# L’uso di un software di Geometria Dinamica nella formazione dei futuri insegnanti

*Aldo Brigaglia, Maria Anna Raspanti, Enrico Rogora*[1]

## Introduzione

L’uso di un software di Geometria Dinamica è ormai pratica piuttosto diffusa tra gli insegnanti della scuola superiore ed è un argomento ampiamente analizzato nei lavori di didattica della matematica.[2] In molte sedi universitarie, le modalità di impiego di tale software nell’insegnamento scolastico è oggetto di studio nei corsi di Didattica. Siamo convinti che possa utilmente trovare spazio, come d’altra parte è stato ampiamente sperimentato in diverse occasioni dagli autori, anche negli insegnamenti di Storia della matematica, Matematiche complementari, Matematiche elementari da un punto di vista superiore, ecc. per facilitare un insegnamento della matematica di carattere “estensivo”, specificamente orientato alla formazione dei futuri insegnanti secondo le idee inizialmente proposte da Guido Castelnuovo e da Felix Klein, che ci sembra utile rivisitare oggi (De Marchis e Rogora, Attualità delle riflessioni di Guido Castelnuovo sulla Formazione dell'insegnante di Matematica 2017) (De Marchis, Menghini e Rogora, The importance of "Extensive teaching" in the education of prospective teachers of Mathematics 2020).

Per insegnamento “estensivo” intendiamo un insegnamento che metta in rilievo i collegamenti piuttosto che gli approfondimenti: tra le varie parti della matematica; tra la matematica e le altre scienze; tra la matematica elementare e le matematiche superiori; tra le diverse fasi che attraversa la genesi storica di un concetto matematico, mirando a un approfondimento concettuale, prima che tecnico, delle conoscenze matematiche. Per contro, l’insegnamento “intensivo”, tipico dei corsi orientati alla formazione dei futuri ricercatori, privilegia l’approfondimento di un tema circoscritto. Nell’insegnamento estensivo della matematica si propone di dare più spazio a ad attività che permettano di sperimentare la logica della scoperta, riducendo quello dedicato all’esposizione assiomatica della teoria.

Lo scopo principale di questo lavoro non è quindi quello di considerare i problemi dell’insegnamento della matematica nella scuola esplorando le possibilità di insegnamento/apprendimento offerte da un ambiente di Geometria Dinamica, ma di esemplificare come sia possibile usare questo software per veicolare l’insegnamento di contenuti superiori secondo modalità appropriate al completamento della formazione matematica del futuro insegnante e di ampliare la prospettiva da cui guardare all'evoluzione di alcuni concetti matematici, recuperando la dimensione geometrica-sintetica caratteristica del pensiero matematico fino al diciannovesimo secolo. La geometria offre numerose occasioni per mettere in risalto l'importanza della matematica nella cultura occidentale, prospettiva che ci sembra molto carente nell'insegnamento scolastico della disciplina. Si pensi ad esempio: all'importanza della geometria nella pittura e nell'architettura; alla suggestione esercitata dalle investigazioni geometriche dello spazio quadridimensionale su pittori e scrittori; all'influenza della geometria non euclidea sull'opera di Escher (Rogora e Tortoriello, Matematica e cultura umanistica 2018). È da questa angolatura, focalizzata sulla formazione del futuro insegnante, che intendiamo presentare e discutere alcuni esempi di uso del software nell'insegnamento della matematica nella formazione universitaria e in quella degli insegnanti in servizio. Questa prospettiva, come vedremo, non è disgiunta da quella di affrontare criticamente l'uso del software in classe ma ne è distinta ed offre anche a quest'ultima un nuovo punto di vista a nostro avviso utile e stimolante.

[1] aldo.brigaglia@gmail.com, marianna.raspanti@gmail.com, rogora@mat.uniroma1.it.

[2] Ci limitiamo a indicare, nella sterminata letteratura a riguardo, alcuni lavori in cui è possibile trovare un’ampia bibliografia, (Laborde, et al. 2006), (Drijvers, Kieran e Mariotti 2010), (Arzarello e Stevenson 2012), (Leung A. 2017), (Mariotti e Baccaglini-Frank, Developing the Mathematical Eye Through Problem Solving in a Dynamic Geometry Environment 2018).

Nel lavoro, dopo alcune considerazioni introduttive sul software, discutiamo quattro percorsi laboratoriali (nel senso di (UMI 2005)). I primi tre (cfr. sez. 2, 3 e 4) sono stati proposti a studentesse e studenti del corso di laurea di matematica mentre il quarto è stato progettato per un nuovo corso di Istituzioni di Matematiche Complementari. Queste attività non sono state pensate come esperimenti didattici ma come semplici laboratori. I risultati ottenuti e le conseguenti riflessioni ci sembrano però rilevanti per sostenere l'importanza di utilizzare un metodo di insegnamento laboratoriale anche nel contesto dell'insegnamento universitario dedicato ai futuri insegnanti.

Nelle sezioni 2 e 3 descriviamo l'uso di GeoGebra per approfondire una teoria già nota agli studenti (la geometria di Euclide) da un nuovo punto di vista. Nelle sezioni 4 e 5 descriviamo l'uso del software per esplorare una teoria non nota. Nella formazione del futuro insegnante crediamo che sia cruciale sperimentare entrambe le modalità. La prima per far sperimentare il ruolo di guida di studenti impegnati in una ricerca matematica. La seconda per far sperimentare il ruolo di studente guidato in tali ricerche. Le due esperienze sono naturalmente collegate, ma è opportuno, a nostro avviso, tenerle distinte e non rinunciare a nessuna di esse.

# 1. Alcune considerazioni introduttive sui software di Geometria Dinamica

Un software di Geometria Dinamica è un software che permette di eseguire costruzioni con riga e compasso e di metterle in movimento, trascinando con l'uso del mouse o facendo muovere lungo una curva prefissata, uno qualsiasi dei punti da cui dipendono le costruzioni stesse. Negli esempi facciamo riferimento al software `GeoGebra` che ha il vantaggio di essere gratuito e ben supportato da una comunità attiva di sviluppatori e utenti. (AAVV s.d.).

Per illustrare, per chi non ne è a conoscenza, la caratteristica dinamica di questo software, consideriamo un esempio semplicissimo. Nella parte sinistra della figura è illustrata la costruzione euclidea del triangolo equilatero (Euclide I.1)[3]. Nella parte destra si mostra il risultato dell'azione di trascinamento, effettuato con il mouse, del punto B (lungo il percorso evidenziato in blu). Si noti come muovendo il punto B e quindi dilatando e ruotando il segmento AB, tutti gli elementi della costruzione (il cerchio di centro A e passante per B, il cerchio di centro B passante per A, il punto C che sta dell'intersezione dei due cerchi e il triangolo ABC) vengono aggiornati di conseguenza, mantenendo le relazioni geometriche che li legano tra loro e, in ultima analisi, agli elementi iniziali. Quindi, in particolare, il triangolo ABC resta sempre equilatero.[4]

[3] Con Euclide X.y intendiamo riferirci alla proposizione y del libro X degli elementi di Euclide, in una qualsiasi edizione, per esempio quella disponibile online in inglese (Euclide s.d.).

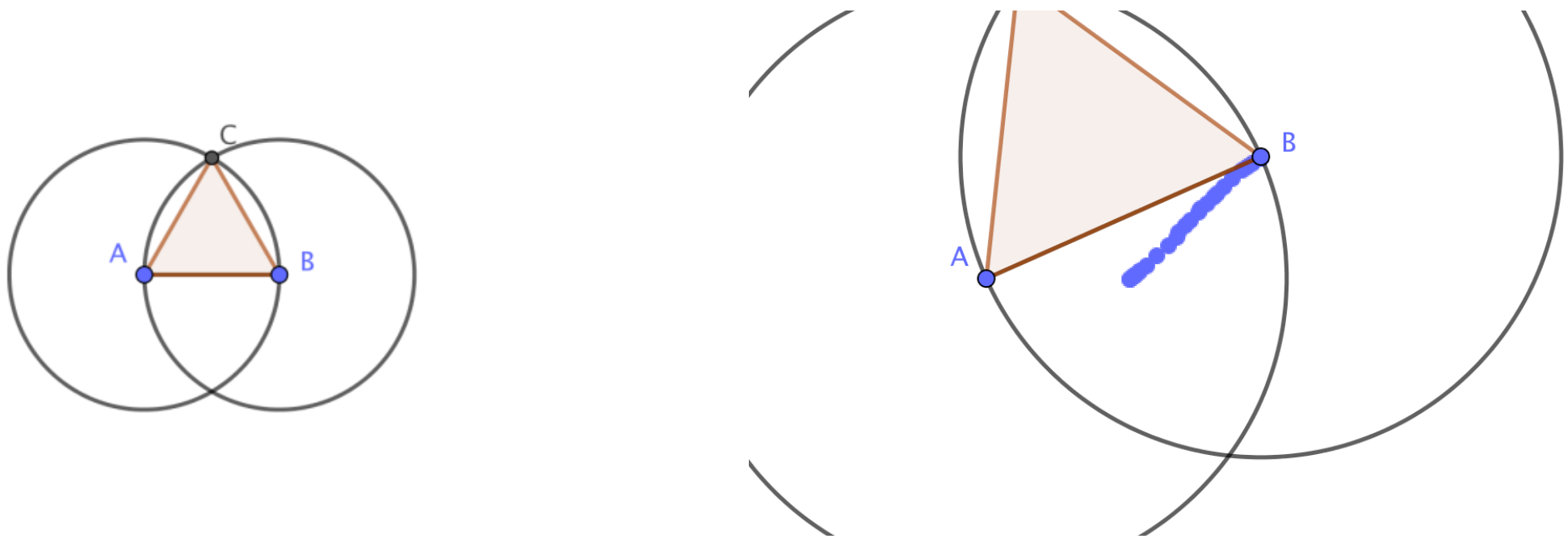

*Figura 1. Trascinando un punto la costruzione si adatta ai nuovi dati iniziali. **

Un software di Geometria Dinamica aiuta ad elaborare un approccio alla geometria in sintonia con il modo di pensare dei matematici antichi e di illustrare il collegamento con approcci più moderni di carattere più analitico e algebrico. Inoltre, permette di illustrare numerosi argomenti senza dover affrontare troppi dettagli tecnici matematici, ma senza rinunciare ad esplorare i collegamenti degli argomenti classici della storia della matematica con argomenti più avanzati e moderni, stimolando un approccio esplorativo e creativo. Infine, l'uso di questo software permette anche di acquisire una competenza informatica utile come formazione professionale, soprattutto, ma non solo, per quanto riguarda la programmazione delle interfacce grafiche.

## 2. Postulati e costruzioni: i primi passi con un software di Geometria Dinamica

Le attività descritte in questa sezione sono state svolte nei corsi di Storia della Matematica e Corso Monografico di Storia della Matematica tenuti tra l'a.a. 2015-2016 e l'anno accademico 2019-2020 presso il corso di Laurea in Matematica della Sapienza, Università di Roma come prima introduzione al software. [5]

Fin dai primi approcci a GeoGebra conviene porre in risalto il collegamento tra il menu delle costruzioni e la presentazione della geometria euclidea negli *Elementi* di Euclide.

Alcune costruzioni corrispondono ai postulati: disegnare un segmento avente gli estremi in due punti dati (Postulato I); prolungare un segmento quanto serve, (Postulato II); disegnare la circonferenza avente centro in un punto e passante per un altro punto (Postulato III). Queste costruzioni sono illustrate nella **Figura 2**.

---

[4] Il lettore interessato può scaricare i file GeoGebra utilizzati nelle figure contrassegnate con il simbolo "*" all'indirizzo web (Rogora, Materiali preparati con GeoGebra 2020). Nei file GeoGebra, i punti indicati in colore blu sono quelli che si possono trascinare con il mouse per muovere la costruzione.

[5] I materiali utilizzati sono raccolti nei siti e-learning dei corsi (Sapienza, Università di Roma s.d.), accessibili con password che può essere richiesta all'indirizzo e-mail del terzo autore. Le attività sono parte integrante dei corsi e sono state svolte presso il laboratorio di Calcolo del Dipartimento di Matematica. Le classi erano composte da un numero di studenti compreso tra 10 e 40.

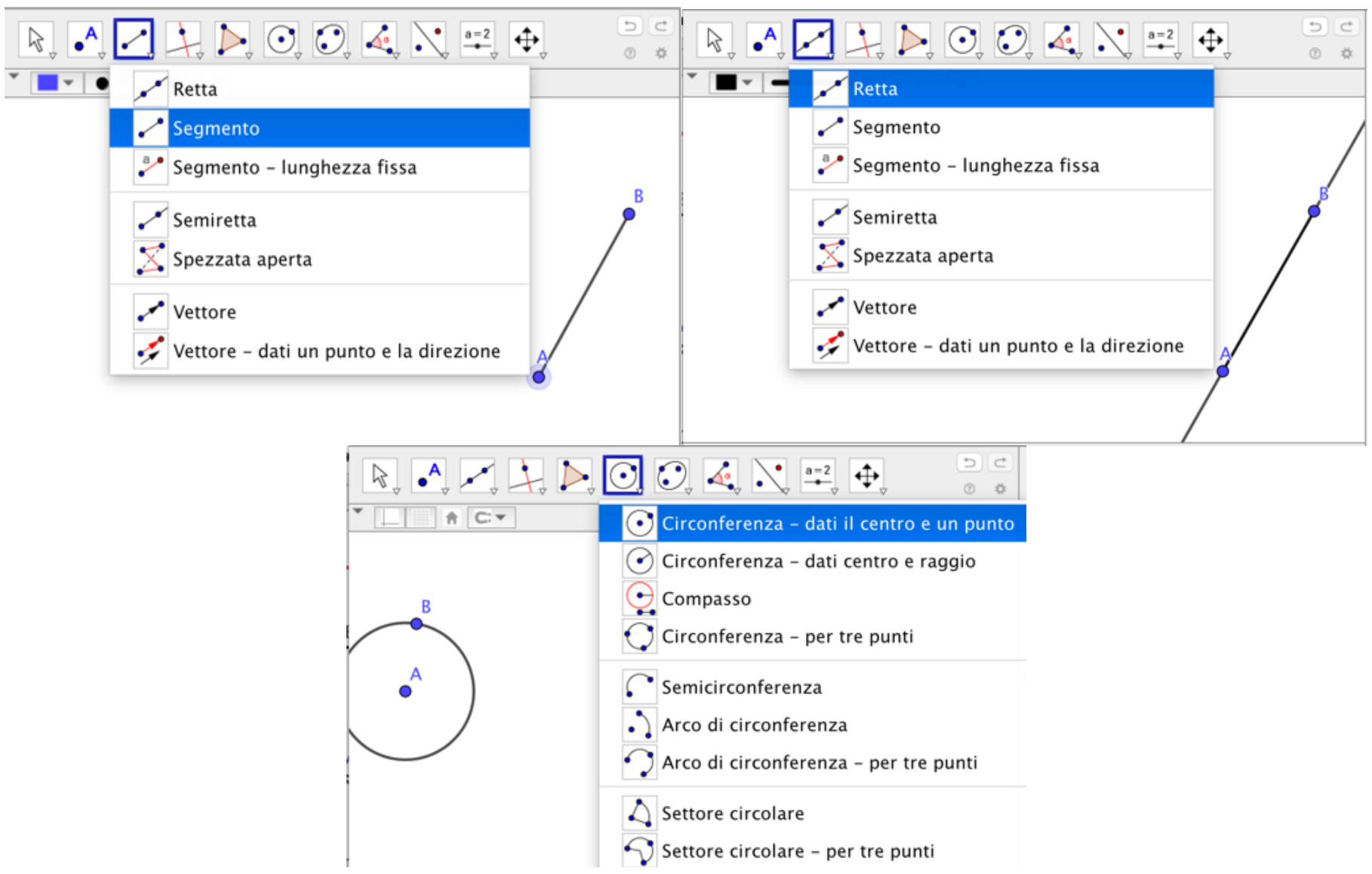

***Figura 2.*** *Le costruzioni corrispondenti ai primi tre postulati di Euclide.*

Un'ulteriore costruzione consiste nel tracciare da un punto la parallela unica a una retta data (Postulato V, nella formulazione di J. Playfair (1748-1819)), cfr. **Figura 3**.

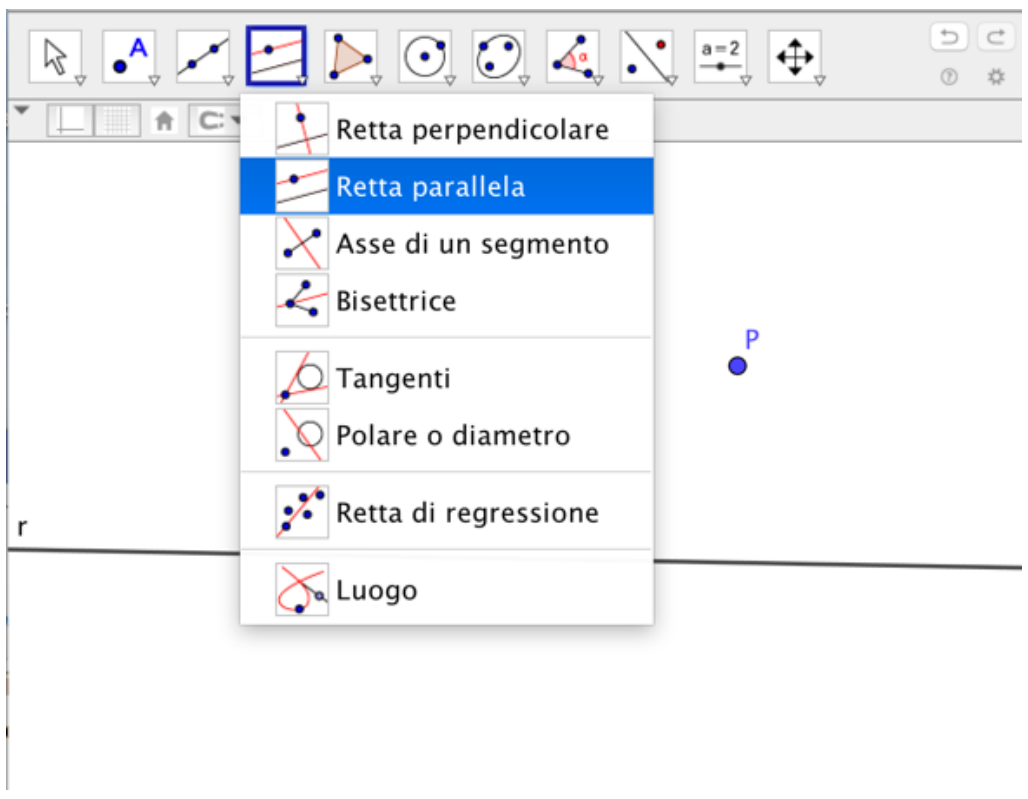

***Figura 3.*** *La costruzione della parallela da un punto.*

Il quarto postulato di Euclide non si identifica con una costruzione del menu ma esprime una particolarità di una delle costruzioni derivate, cioè la costruzione della perpendicolare per un punto a una retta assegnata e offre l'occasione per un'interazione proficua con il software. La perpendicolare (che il software implementa con la costruzione "Retta perpendicolare", sopra alla costruzione "Retta parallela" nella **Figura 3**), viene realizzata da Euclide in (Euclide I.11) e in (Euclide I.12) a partire dalla costruzione del triangolo equilatero, fatta in (Euclide I.1). Il quarto postulato chiede che facendo questa costruzione per due rette distinte che passano per due punti qualsiasi, gli angoli che ottengo siano uguali (o congruenti). L' uguaglianza tra angoli retti, ci può apparire ovvia, ma non discende dagli altri postulati.
Per chiarire questo punto riassumiamo il percorso laboratoriale utilizzato per introdurre gli studenti dei corsi di Storia della matematica Corso Monografico di storia della matematica all'uso di GeoGebra.
Si inizia l'attività chiedendo di costruire con riga e compasso il triangolo equilatero di lato assegnato, utilizzando solo le costruzioni corrispondenti ai primi tre postulati. Il problema è molto semplice, ma serve a mettere in luce la necessità di un postulato ulteriore. Il software è costruito in modo da mascherare la necessità di impiegare un comando per far apparire i punti di intersezione di due circonferenze, ma esaminando il protocollo di costruzione del triangolo si vede che è necessaria una nuova operazione, quella che genera il punto di intersezione di due circonferenze che hanno centro nei due estremi del segmento e passano per l'altro estremo.

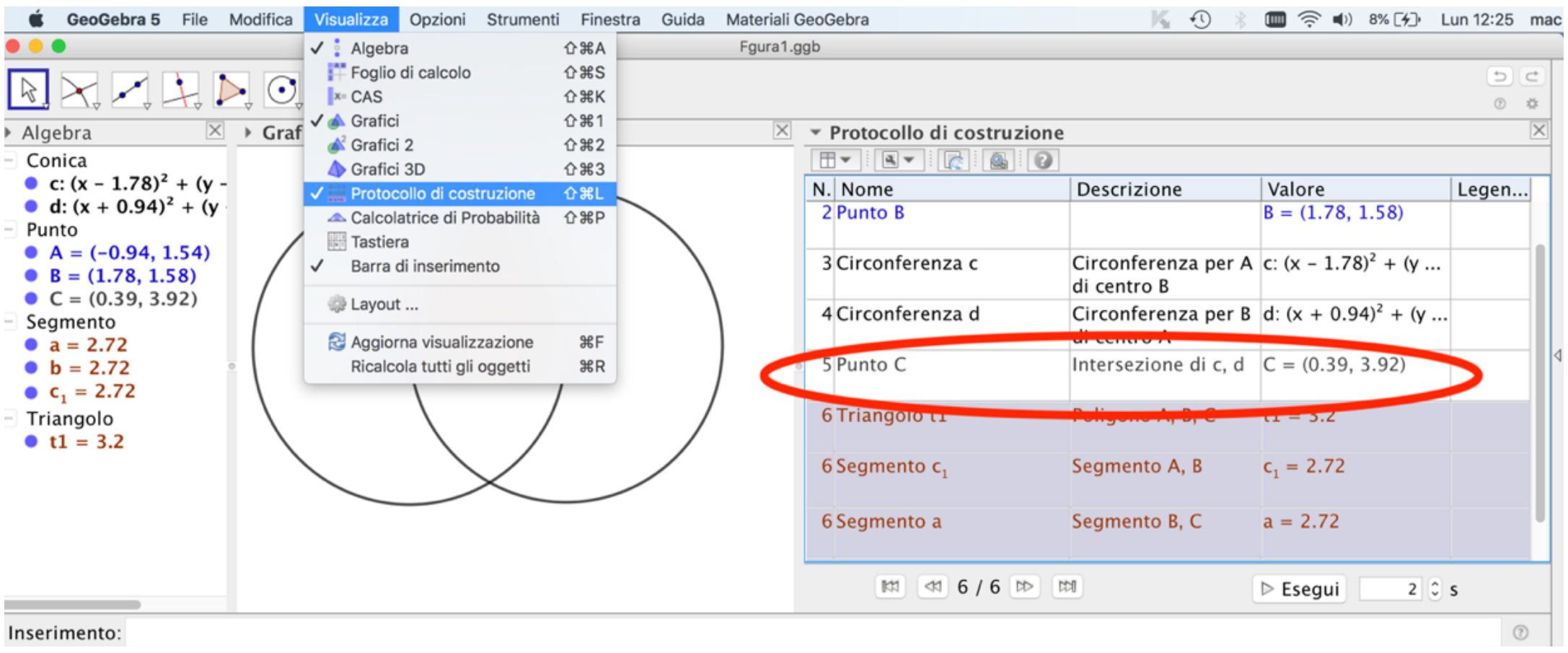

***Figura 4.*** *Dal Menu "Visualizza" è possibile far apparire la finestra "Protocollo di costruzione" dove appare la necessità di impiegare un comando, per intersecare due circonferenze, che non si può definire a partire dai comandi-postulati.*

Il fatto che i punti di intersezione esistano non è conseguenza degli altri postulati. Solo nel XIX secolo ci si è accorti di questa necessità, che, guardando al software, ci appare abbastanza evidente. Bisogna implicitamente o esplicitamente ricorrere a una quarta operazione, corrispondente al postulato nascosto nella Proposizione I.1.

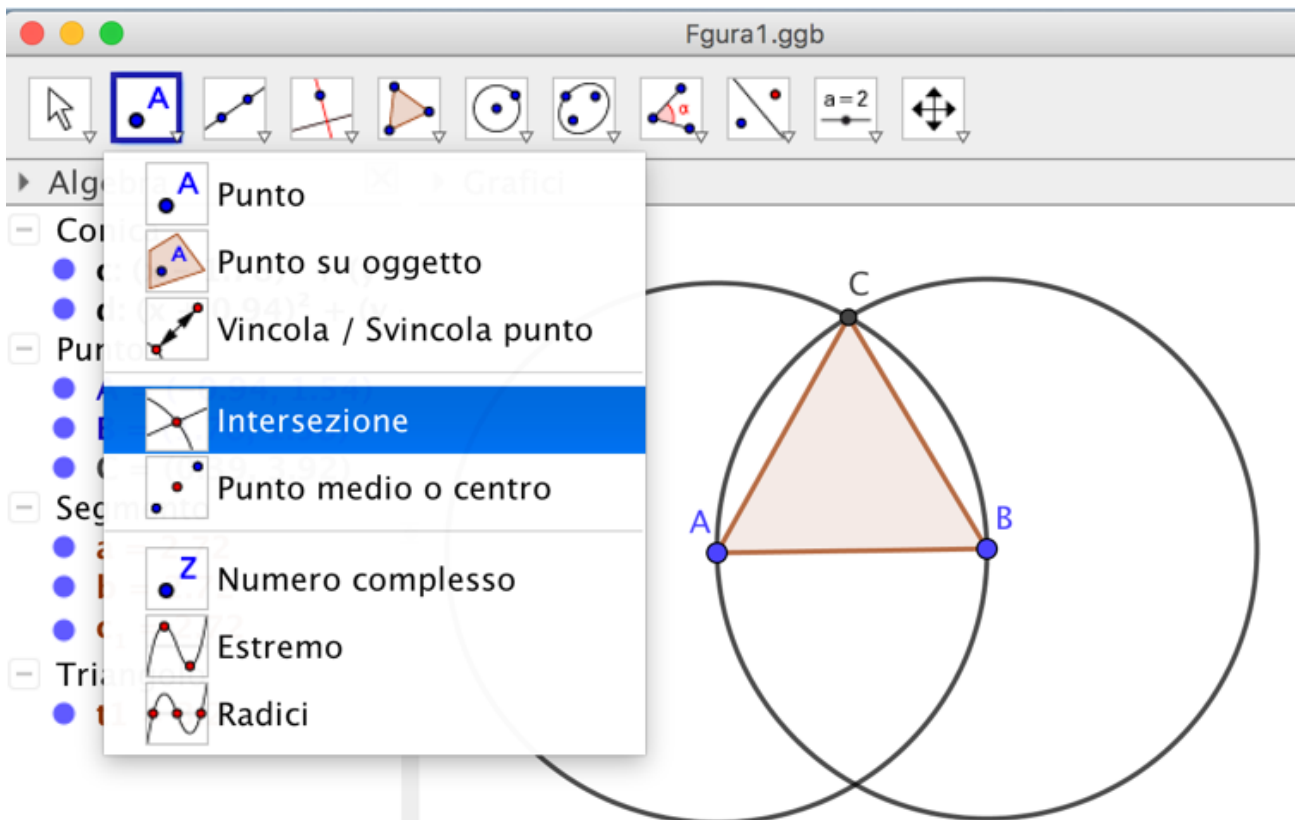

***Figura 5.*** *Nel Menu delle costruzioni appare la costruzione corrispondente ad (una formulazione più generale del) postulato nascosto in Euclide I.1. Il software però è programmato in modo di poterla usare implicitamente in alcune costruzioni. Inoltre, anche se è possibile personalizzare la barra degli strumenti escludendone alcuni, non è possibile inibirne l'uso implicito, almeno nella versione 5. Dal punto di vista didattico questo è un aspetto dove il software potrebbe essere migliorato.*

Il problema successivo che è stato posto ai partecipanti al laboratorio, è quello di affrontare il problema di "trasportare" un segmento da un punto a un altro, lasciando piena libertà di usare tutte le costruzioni disponibili nel menu degli strumenti di GeoGebra (che possono eventualmente essere ridotte grazie allo strumento "Personalizza la barra degli strumenti").

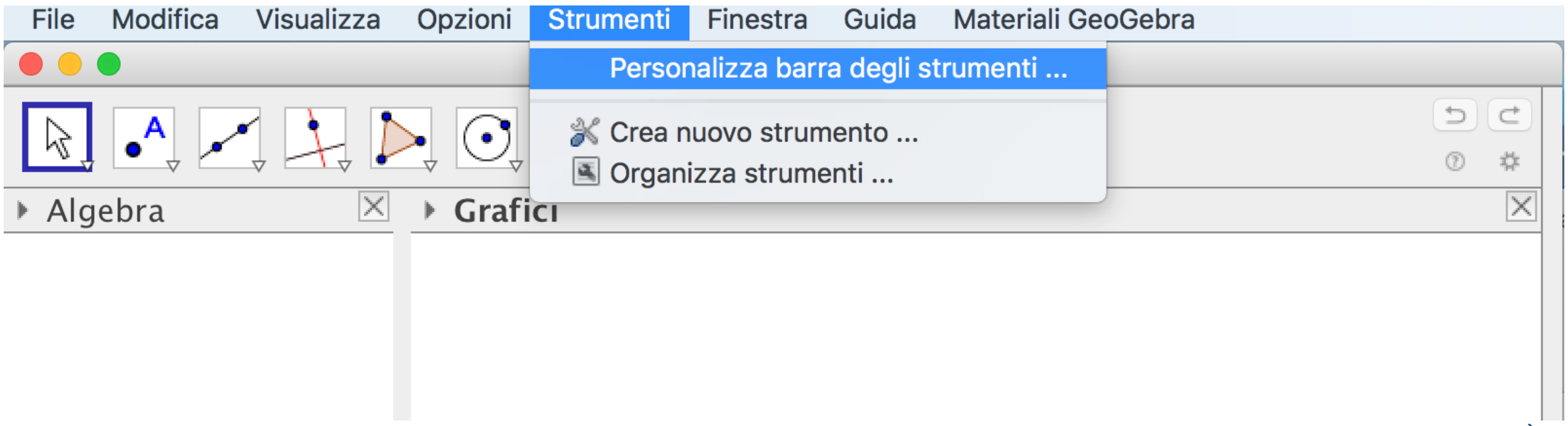

***Figura 6.*** *Nel Menu degli strumenti l'opzione "Personalizza barra degli strumenti" permette di togliere e aggiungere strumenti. È possibile partire con i soli strumenti-postulati e aggiungere gli altri a mano a mano dopo aver letto sugli Elementi o scoperto nel laboratorio come implementarli. Con questa attività la logica degli Elementi di Euclide si lega strettamente alla concretezza di un problema relativo ad un software, con evidenti vantaggi didattici, anche nel senso di far toccare con mano l'utilità di organizzare logicamente e assiomaticamente una teoria. cfr anche (Mariotti, Paola, et al. 2004).*

Le soluzioni al problema del trasporto di un segmento, che abbiamo raccolto nei diversi laboratori, si distinguevano per la scelta di usare una delle costruzioni seguenti:[6]

1. "Retta parallela";
2. "Compasso";
3. "Circonferenza dati centro e raggio";
4. "Segmento – lunghezza fissa".

Successivamente, abbiamo chiesto di trasportare un segmento eliminando tutte le costruzioni del software escluse quelle corrispondenti ai primi tre postulati e al postulato nascosto in Euclide I.1, suggerendo di usare la costruzione del triangolo equilatero.

Dopo aver discusso i tentativi dei partecipanti, si è letta la soluzione di Euclide (Euclide I.2), osservando come questa costruzione fornisce una procedura operativa per confrontare due segmenti qualsiasi e permette, a partire dalle costruzioni/postulato, di implementare la costruzione "Compasso", ma non le altre tre della lista qui sopra, in particolare non la costruzione "Retta parallela". Con il quinto postulato, abbiamo un modo alternativo e rapido per trasportare un segmento rispetto a (Euclide in I.2), ma il metodo di Euclide permette il trasporto anche in contesti in cui il quinto postulato non vale, per esempio nella geometria iperbolica.

Come è emerso dal laboratorio, è possibile risolvere il problema anche con gli strumenti "Circonferenza dati centro e raggio" e "Segmento – lunghezza fissa". C'è però una differenza tra questi strumenti e lo strumento "Compasso". Per utilizzare quest'ultimo si chiede di specificare il raggio in forma "geometrica", cioè selezionando un segmento congruente al raggio. Per utilizzare gli altri, si chiede di specificare il raggio in forma "numerica". Questa identificazione della classe di segmenti congruenti a un segmento dato con un numero non fa parte della geometria di Euclide e il software offre l'occasione di metterlo in luce (ma anche di mascherarlo, quando identifica, nella vista "Algebra", un segmento o un angolo con un numero"). Quindi, solo lo strumento "Compasso" può essere costruibile a partire dagli strumenti-postulati.

---

[6] È possibile infine trasportare un segmento anche facendo uso delle costruzioni del menu trasformazioni. Il problema di costruire queste trasformazioni con le costruzioni postulate da Euclide ovvero di sostituire le costruzioni postulate da Euclide con una opportuna scelta di trasformazioni è un argomento stimolante, ~~che può emergere facilmente dalla proposta laboratoriale che stiamo descrivendo, ma~~ che abbiamo deciso di non approfondire in questa sede.

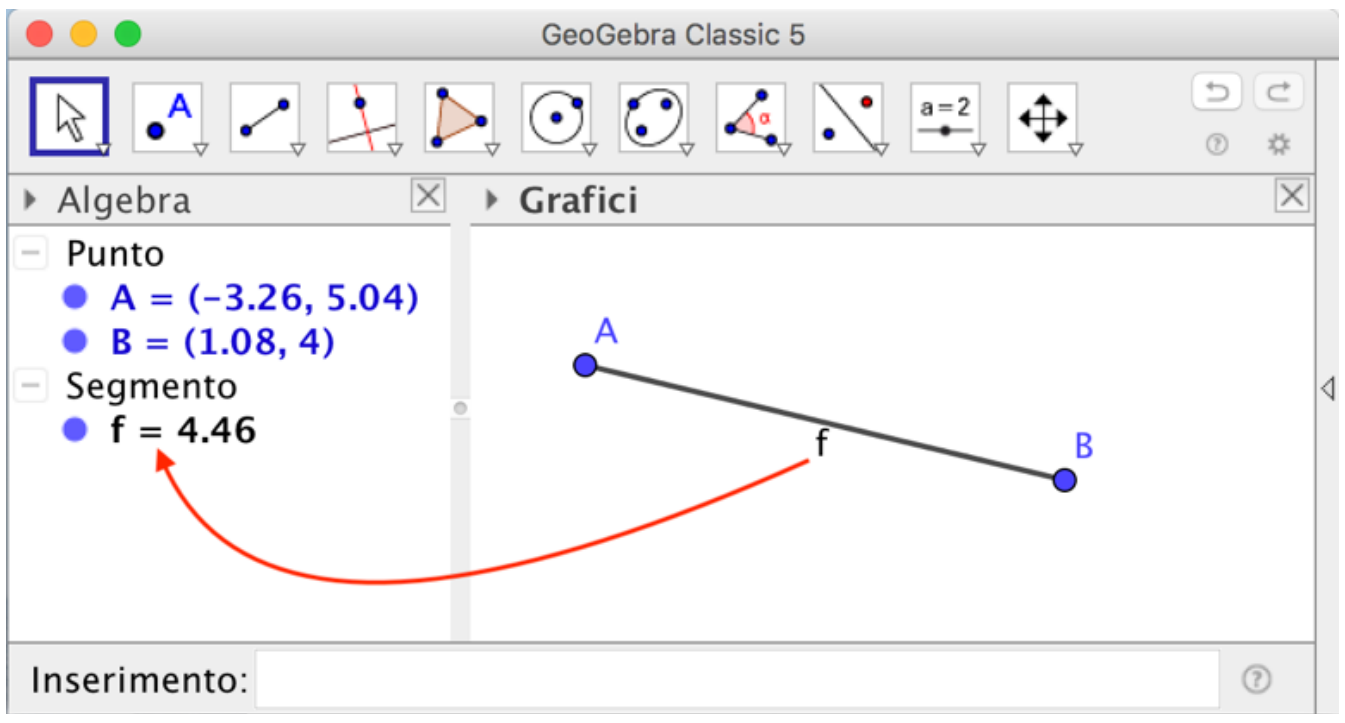

***Figura 7.*** *Il menu Algebra del software facilita l'identificazione di un segmento (o più correttamente di una classe di equivalenza di segmenti congruenti) con un numero, che dal punto di vista didattico è discutibile.*

Torniamo a questo punto al problema di trasportare gli angoli. Nel laboratorio viene richiesto, analogamente a quanto abbiamo fatto per il trasporto del segmento, di trasportare un angolo con gli strumenti di GeoGebra. Per risolvere il problema c'è ora meno scelta di quella che avevamo nel caso dei segmenti, e la scelta ovvia è quella di utilizzare lo strumento "Angolo di data Misura".
Per utilizzare questo strumento si chiede di specificare l'angolo in forma numerica. Non esiste uno strumento "puramente geometrico" analogo al "Compasso" del problema precedente. Questa è una asimmetria facilmente aggirabile,[7] ma a nostro avviso si tratta di una pecca concettuale del software[8].
Poniamo quindi il problema di trasportare un angolo qualsiasi da un punto del piano a una semiretta, con le sole costruzioni relative ai primi tre postulati e al postulato nascosto dell'intersezione.
Una possibile costruzione per trasportare un angolo, che non è stato difficile far riscoprire nel laboratorio, è quella descritta in (Euclide I.23).
La costruzione si può fare con riferimento esclusivamente alle prime quattro costruzioni/postulato, ma per garantire che produca quanto richiesto, cioè un angolo congruente a quello assegnato bisogna far ricorso ad ulteriori postulati. Torniamo per un momento al problema relativo al trasporto di un segmento. La dimostrazione che la costruzione di Euclide funziona, dipende dal fatto che due segmenti con un estremo in comune sono congruenti se e solo se il secondo estremo dell'uno appartiene alla circonferenza, centrata nell'estremo comune, passante per il secondo estremo dell'altro. Non abbiamo ancora nulla invece che ci permetta di confrontare gli angoli e su cui sia possibile basare la dimostrazione che la costruzione proposta da Euclide sia corretta. Esistono due costruzioni, tra le prime 12 proposte da Euclide negli *Elementi*, che producono angoli. Quella del triangolo equilatero (I.1) e quella della perpendicolare ad una retta condotta da un punto (I.12 e I.13). Una possibilità, se vogliamo confrontare gli angoli è quella di concordare sul fatto che queste costruzioni producono angoli congruenti. Si potrebbe assumere, per esempio, che gli angoli di due triangoli equilateri qualsiasi sono congruenti. Questa assunzione, che è formalmente analoga a quella postulata da Euclide nel quarto postulato, è però meno generale. Non vale infatti né per la geometria iperbolica, né per la geometria ellittica[9], per cui valgono invece i primi tre postulati

[7] Quando si definisce un angolo, esso viene registrato nella finestra "Algebra". Lo strumento "Angolo di data misura" accetta come input il nome dell'angolo, oltre che un numero. Nella variabile con tale nome, la misura dell'angolo è registrata con una precisione maggiore di quella che appare nel menu algebra, dando l'illusione di una misura "esatta" e quindi di trattare direttamente con l'oggetto geometrico, come nel caso dei segmenti.
[8] Uno strumento geometrico "trasporta angolo" era presente nelle prime versioni di Cabri [cfr. Mariotti 2002, ZDM].
[9] E neppure per la geometria sulla sfera. In essa, assumendo come rette i cerchi massimi, non vale il primo postulato come viene enunciato solitamente, cioè nella forma: per due punti passa una ed una sola retta. Ma, come osserva Lucio Russo in (Russo e Pirro 2017), la formulazione originale del postulato si deve tradurre piuttosto: si richieda di poter condurre una linea retta da qualsiasi punto a ogni altro punto. In questa forma, senza riferimento all'unicità, vale anche per la geometria sferica.

(modificando il secondo nel senso indicato da Riemann). In queste geometrie l'angolo di un triangolo equilatero dipende dalla lunghezza del lato. L'uguaglianza degli angoli prodotti con la costruzione della perpendicolare, proposta da Euclide in (Euclide I.12) e (Euclide I.13), vale invece anche nelle geometrie non euclidee. Inoltre, dalla uguaglianza tra gli angoli retti e dalla proposizione (Euclide I.9) che insegna a dividere un angolo a metà segue l'uguaglianza di tutti gli angoli $\gamma(h,k)$ che hanno rapporto con l'angolo retto uguale a $\frac{k}{2^h}$, dove h e k sono numeri naturali. Sarebbe quindi possibile basare il confronto tra due angoli qualsiasi $\alpha$ e $\beta$ confrontando l'insieme degli angoli $\gamma(h,k)$ che sono minori/uguali/maggiori di $\alpha$ con l'insieme di quelli che sono minori/uguali/maggiori di $\beta$, analogamente a quanto fa Euclide fa nel libro V per confrontare i rapporti. Esiste però una possibilità più diretta per trasportare e confrontare gli angoli, basata sul terzo criterio di uguaglianza dei triangoli (Euclide I.23) che dipende a sua volta da (Euclide I.4), che appare come proposizione negli *Elementi* ma che deve essere invece aggiunta ai postulati, come nei fondamenti proposti da Hilbert, oppure dedotta da un ulteriore postulato. Come osserva Hilbert nei Grundlagen (Hilbert 2009)[10], l'uguaglianza degli angoli retti si può dimostrare con (Euclide I.23), ed è quindi ridondante quando si aggiungano tra i postulati i criteri di congruenza dei triangoli, che vanno comunque aggiunti ai postulati di Euclide. Resta però a nostro avviso interessante notare lo stretto collegamento del quarto postulato di Euclide con una costruzione fondamentale, quella della perpendicolare. Ciò rende il quarto postulato più omogeneo agli altri di quanto non appaia a prima vista.

Si noti infine che la proposizione "tutti gli angoli retti sono uguali", non ha un analogo per i segmenti. Non esistono infatti, nella geometria euclidea, "segmenti canonici", che possano essere definiti in ogni punto in maniera univoca e non convenzionale, come si possono definire gli angoli retti! Nella geometria iperbolica invece, dove il legame tra angoli e segmenti è più stretto, e dove vale il "quarto criterio di congruenza dei triangoli", cioè due triangoli con angoli congruenti sono congruenti, è possibile definire segmenti canonici, per esempio, il lato di un triangolo equilatero di angoli uguale a $\pi/4$.

L'analisi del menu delle costruzioni di GeoGebra ci offre quindi, come abbiamo cercato di argomentare, uno strumento estremamente efficace per guidare criticamente lo studio degli *Elementi* e per comprendere la loro bellezza e profondità in maniera da coinvolgere il lettore in una esplorazione attiva.

Bisogna osservare però che alcune tra le costruzioni presenti nel menù di GeoGebra non sono sempre eseguibili con riga e compasso. Questo è un fatto che dovrebbe essere messo meglio in rilievo nel software, segnalando in modo chiaro le tante operazioni del menù degli strumenti che non lo sono: per esempio, la costruzione di un poligono regolare qualsiasi; la costruzione di un angolo qualsiasi; le coniche[11]; il luogo generato da un punto al variare di un altro punto su una curva, ecc[12].

---

[10] L'edizione originale di Hilbert, in tedesco, è del 1899. Ad essa seguì immediatamente la traduzione francese a cura di Laugel (1900) e quella inglese a cura di Townsend (1902). Per la prima edizione italiana, a cura di Pietro Canetta, si dovette attendere il 1970.

[11] Con riga e compasso è possibile costruire i punti di una conica che verificano certe condizioni, per esempio quella di stare su una retta. In questa maniera è possibile costruire un numero finito di punti di una conica ma non l'intera conica (salvo naturalmente la circonfrenza!).

[12] Si noti però che nei "luoghi piani", opera smarrita di Apollonio, sono descritti i luoghi generati da punti variabili su una retta o un cerchio nel caso in cui i luoghi risultanti siano ancora una retta o un cerchio~~, p.e, l'inversione circolare~~.

# 3. Pentagono, pentagramma e la logica della scoperta matematica

La logica della scoperta matematica non è quella dell'esposizione assiomatica di una teoria (Lakatos 1979). Anche se un'esposizione assiomatica presenta indubbi vantaggi nell'insegnamento universitario, in quanto permette di raggiungere un'estrema concisione e di massimizzare la quantità di informazioni trasmesse, è sbagliato utilizzare esclusivamente questo registro nella formazione del futuro insegnante, che dovrà essere in grado di gestire nelle proprie classi un lavoro di carattere laboratoriale che ricalchi più da vicino l'iter della scoperta matematica invece di quello dell'esposizione asettica di teorie cristallizzate. A tal fine, risulta molto utile e stimolante l'uso di GeoGebra nell'insegnamento dei corsi universitari per i futuri insegnanti per esplorare in un contesto laboratoriale alcuni argomenti, sia di carattere elementare, come in questa sezione, sia di carattere più elevato, come nelle successive.

Le attività descritte in questa sezione sono state proposte in un corso di Storia della Matematica di 30 studenti tenuto alla Sapienza, Università di Roma, nell'a.a. 2019-20.

Dopo aver letto e discusso il brano del Menone in cui Socrate fa ricordare allo schiavo come raddoppiare un quadrato e dopo aver letto e discusso il capitolo dell'opera di Lakatos "Dimostrazioni e confutazioni", relativo alla formula di Eulero per i poliedri (Lakatos 1979), è stato chiesto alle studentesse e agli studenti del corso di scrivere un dialogo ispirato al Menone e al brano di Lakatos in presentare la costruzione di Euclide del pentagono regolare iscritto in una data circonferenza e la dimostrazione che la figura ottenuta è effettivamente equilatera ed equiangola. Successivamente è stato chiesto di «Realizzare con GeoGebra un percorso *storicamente plausibile*[13] che porti alla scoperta della costruzione del pentagono regolare nella maniera che la studentessa/lo studente giudica naturale». A supporto della realizzazione di questo progetto sono state dedicate alcune sessioni di Laboratorio in cui i partecipanti potevano discutere tra loro e con l'insegnante le idee per realizzare il loro percorso. La richiesta di elaborare un percorso storicamente plausibile serve a mettere il futuro insegnante nelle condizioni di guardare al problema con occhi vicini a quelli dei suoi futuri studenti e studentesse, mettendo momentaneamente da parte le conoscenze matematiche che gli derivano dagli studi universitari.

L'uso di GeoGebra si è rivelato molto utile, sia in fase esplorativa, come ambiente in grado di stimolare l'ideazione di ipotetici percorsi, sia in fase comunicativa, per illustrare e discutere il percorso con gli altri partecipanti e con l'insegnante.

Ed ecco in sintesi il "percorso storicamente plausibile" costruito dall'insegnante sulla base degli spunti proposti dai lavori degli studenti. Il percorso costruito dall'insegnante ha la forma di un laboratorio per una scuola superiore ed è stato discusso in laboratorio con gli studenti del corso di Storia della Matematica, cui è stato chiesto di sperimentare le attività laboratoriali proposte nel percorso e di riflettere su di esse. Presentiamo il percorso nella forma di un racconto che un ipotetico insegnante di scuola superiore fa alla sua classe in un laboratorio. Durante il racconto l'insegnante pone alla classe domande, cui la classe può rispondere usando GeoGebra. Queste stesse domande sono state poste agli studenti universitari durante la discussione del percorso. Cominciamo il racconto.

Possiamo facilmente immaginare che giocando con riga e compasso gli antichi matematici abbiano scoperto molto presto come costruire un esagono regolare di lato assegnato. Dato un segmento AB e la circonferenza centrata in A e passante per B basta riportare il segmento AB su una corda BC della circonferenza, congruente ad AB e iterare la costruzione di corde adiacenti e congruenti, (CD, DE, EF, FG) finché, all'ultimo passo, l'estremo H della corda GH viene a coincide con B. [Si chiede agli studenti di ripetere la costruzione.]
I punti B,C,D,E,F e G sono quindi i vertici di un esagono equilatero di lati BC, CD, DE, EF, FG e GB e equiangolo.

---

[13] Per percorso storicamente plausibile intendiamo un percorso che usi le tecniche e le conoscenze dei matematici greci e si sviluppi conformemente alla loro (plausibile) intuizione. Si tratta quindi di evitare l'algebra dei numeri, di usare l'algebra delle proporzioni, di evitare, per quanto possibile, l'identificazione di una misura con un numero. L'uso critico di GeoGebra e la conoscenza approfondita delle sue funzionalità, è di grande aiuto in questo.

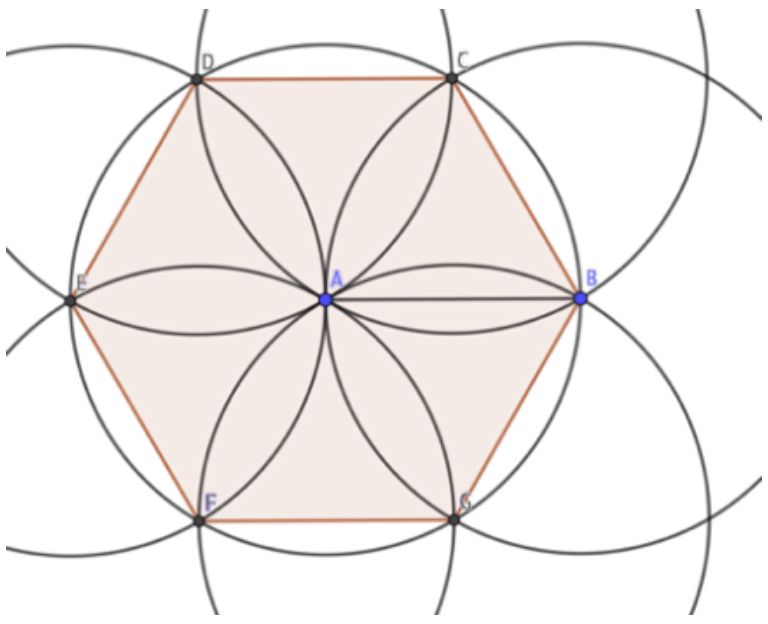

***Figura 8.*** *La costruzione dell'esagono regolare.*

È possibile immaginare la curiosità degli antichi matematici di esplorare cosa succede, costruendo corde congruenti a segmenti diversi dal raggio. Per esempio qualcuno potrebbe essersi chiesto: se invece di riportare sulla circonferenza il raggio AB riportiamo il segmento AH pari alla metà di AB e ripetiamo cosa succede? [Questa domanda viene posta alla classe]. La spezzata non si richiude esattamente in B.

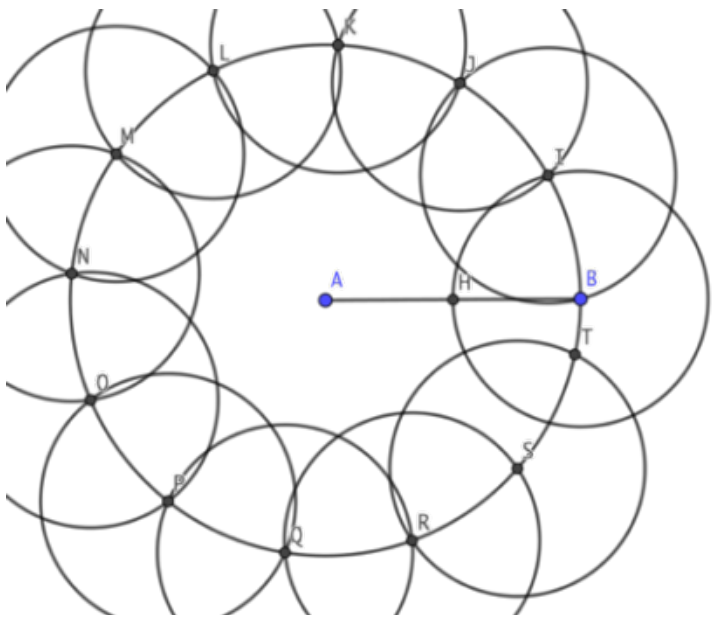

***Figura 9.*** *La costruzione della spezzata di corde congruenti.*

A questo punto viene posta alla classe la seguente domanda: esistono altre divisioni del segmento AB in due parti AX e BX che possiamo usare per costruire, con lo stesso procedimento, poligoni regolari? Possiamo anche ipotizzare che questo "programma di ricerca" sia stato subito precisato nella seguente maniera: cerchiamo una "costruzione" del punto X e non una scelta casuale, in modo da poter replicare la costruzione e condividere la scoperta. Viene chiesto alla classe: quali punti volete costruire? ce n'è qualcuno che permette di chiudere in B la costruzione? Un'idea proposta dagli studenti universitari è quella di iterare la costruzione del punto medio. In questa maniera si possono costruire i punti che individuano i segmenti che hanno con AB rapporto uguale a $\frac{h}{2^k}$. Poiché questi punti sono "densi" in AB possiamo con questo metodo approssimare tanto bene quanto vogliamo la costruzione di un poligono regolare con un numero di lati qualsiasi (ovviamente maggiore di 6), ma si tratta di un approccio poco efficiente e, naturalmente, neppure esatto.[14]

[14] Non era chiaro ai primi matematici greci, come non lo è per la maggioranza degli studenti che iniziano a riflettere su questi problemi, che il procedimento di approssimazione di un punto di AB con un insieme denso di punti, non debba necessariamente arrestarsi dopo un numero finito di passi, come prevedeva, per esempio, la teoria delle monadi geometriche di Pitagora.

La possibilità di una facile costruzione approssimata dell'ettagono regolare, scoperta nel corso dell’attività laboratoriale da uno studente del corso di Storia, si ottiene scegliendo $BX/AB=\frac{7}{8}$.

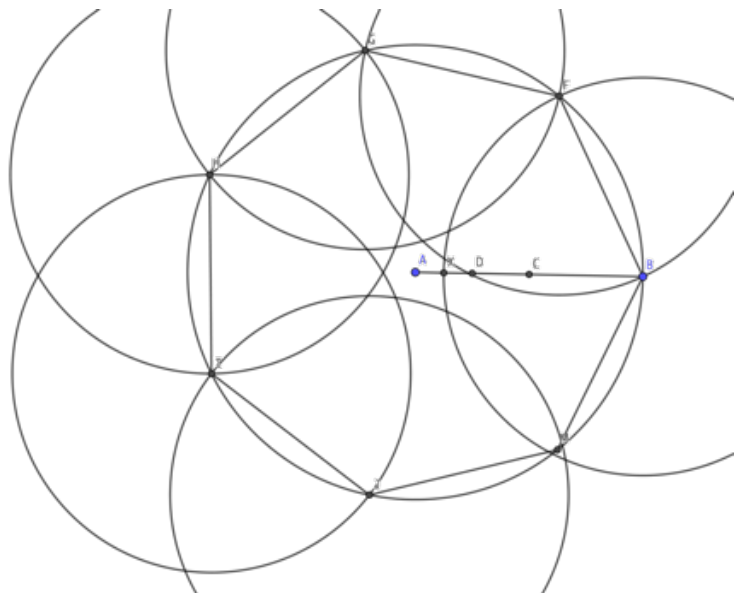

***Figura 10.*** *Una costruzione approssimata dell’ettagono regolare “scoperta” in laboratorio. **

Possiamo anche ipotizzare che gli antichi matematici greci si fossero posti il problema di costruire i punti X tali che AX abbia rapporto razionale con AB. La costruzione di tali punti offre spunti interessanti per un laboratorio scolastico, rivolto a trattare il problema di suddividere un segmento in parti razionali, che però non approfondiamo in questa sede.
Il passo successivo può essere facilmente motivato con un linguaggio moderno, ma non può far parte del percorso storicamente plausibile che stiamo costruendo: costruiamo punti X tali che AX abbia rapporti irrazionali con il segmento AB. Come possiamo arrivare ai punti irrazionali in modo elementare? Uscendo dal segmento! Per esempio, costruendo il quadrato che ha come diagonale AB e riportando il suo lato.

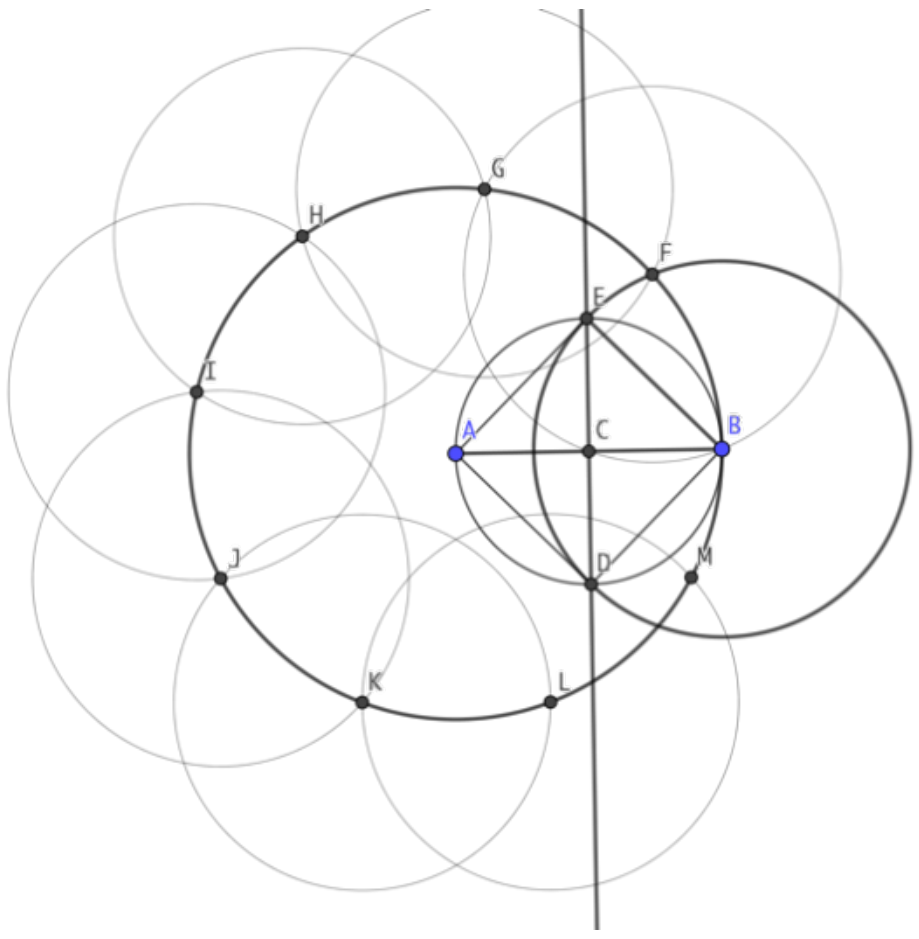

***Figura 11*** *Anche riportando le corde congruenti al lato del quadrato di diagonale AB non riusciamo a chiudere la spezzata ma abbiamo acquisito l’idea di uscire dal segmento AB per costruire una divisione AX, BX del segmento AB. **

L’idea che per costruire la divisione bisogna considerare punti esterni alla retta AB è fondamentale per proseguire la costruzione del percorso. Questo cambiamento di prospettiva, se non viene in mente ai partecipanti al laboratorio (per gli studenti universitari era ovvio, vista l’attività che avevano svolto sulle costruzioni di Euclide, ma non è affatto scontato per gli studenti di scuola) deve essere suggerito dall’insegnante cercando però sempre di:

> *“Non rivelare immediatamente tutto quello che sai e devi spiegare agli studenti - fallo congetturare dagli studenti prima di dirlo - fai loro scoprire, da soli, quanto più è possibile.”*
> (Polya 1971), vol. 2.

Una volta accettata l'idea di uscire dal segmento per cercare una divisione che ci permetta di costruire un poligono regolare, possiamo guidare il laboratorio alla scoperta della seguente costruzione[15].

Costruiamo la perpendicolare in B al segmento AB. Riportiamo sulla perpendicolare che abbiamo appena costruito e con un estremo in A un segmento AC congruente ad AB. Sul segmento AC costruiamo il punto medio (senza bisogno di suggerirlo, è questo il punto che viene scelto generalmente dagli studenti) che indichiamo con la lettera D. Intersechiamo in J la circonferenza centrata in D e passante per A con la retta AC. Riportiamo il segmento JB su AB con la circonferenza centrata in B e passante per J e individuando il punto K.

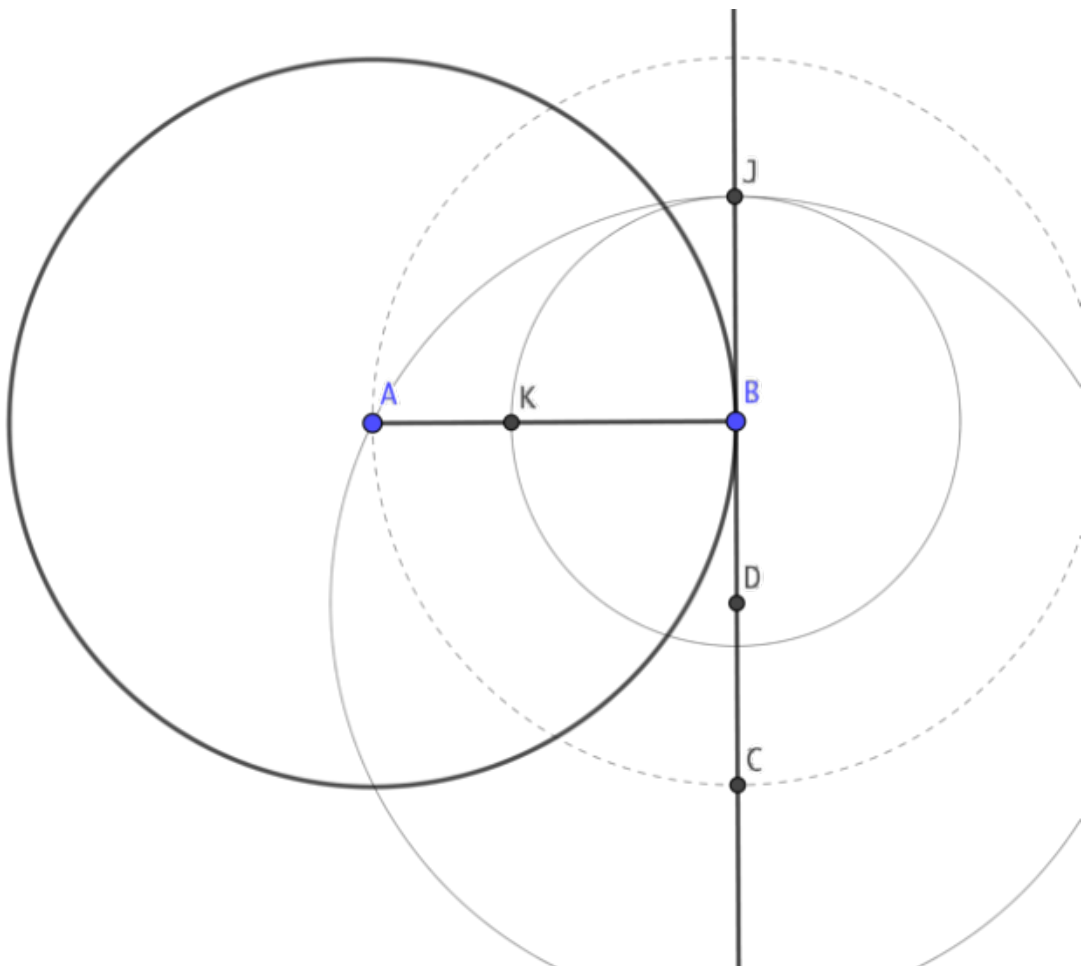

***Figura 12.*** *Una sezione "irrazionale".*

La costruzione delle corde adiacenti congruenti al segmento BK sembra chiudersi perfettamente, suggerendo la costruzione del decagono regolare e, congiungendo un vertice sì e uno no, anche del pentagono regolare.

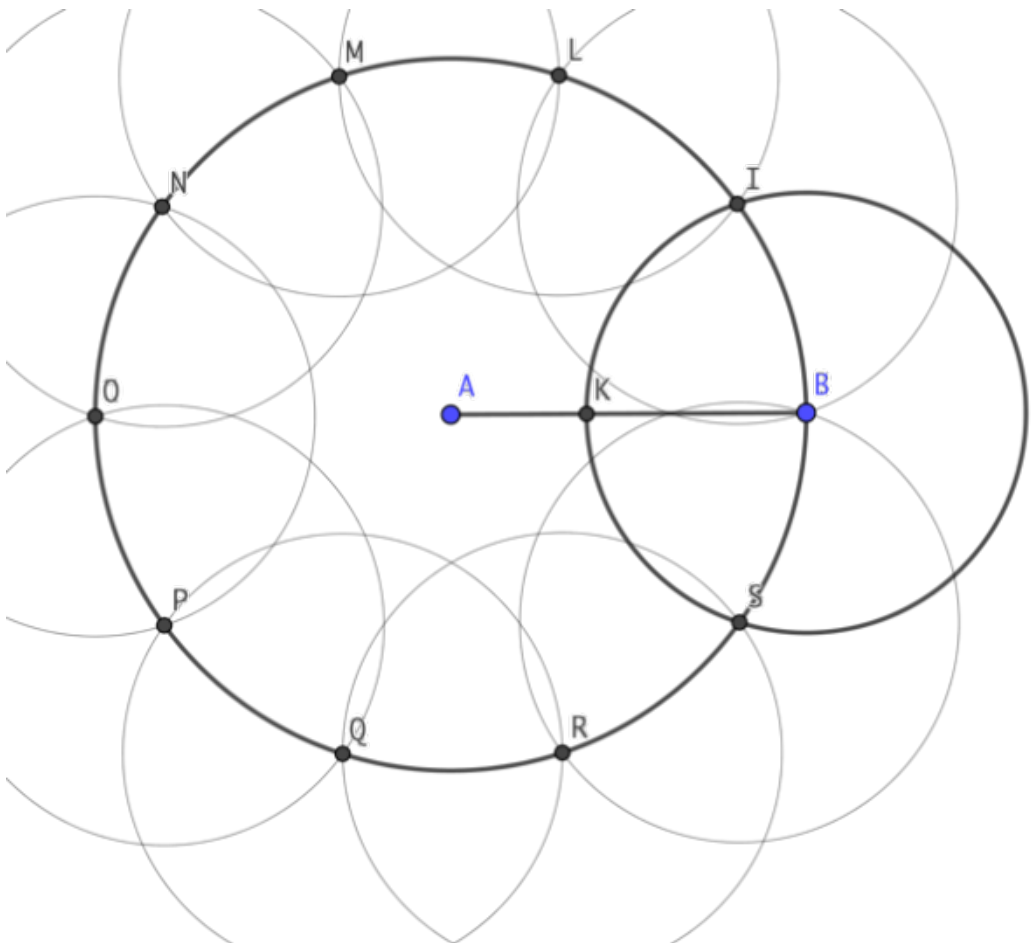

***Figura 13.*** *La scelta di K permette di costruire il decagono regolare.*

In una attività laboratoriale fatta a scuola, quanto tempo ci aspettiamo che i partecipanti impieghino per scoprire la costruzione del punto K e osservino che la spezzata di segmenti adiacenti congruenti a BK si chiude? Dipende molto dai suggerimenti che un insegnante darà, ma un insegnante ben preparato è in grado di portare la classe a scoprire questa costruzione in un tempo ragionevole. Vale la pena? Crediamo di si, perché crediamo che insegnare a congetturare e a porsi le domande giuste quando si affronta un problema è almeno altrettanto importante che insegnare a dimostrare[16]! È chiaro che non è possibile dimostrare solo ciò che congettura una classe, ma dedicare tempo ad attività specifiche che stimolino la capacità di congetturare è necessario. Il percorso universitario dei futuri insegnanti li prepara sufficientemente a questo compito? Crediamo di no. L'uso di un software di Geometria Dinamica (ma anche altri software, come per esempio *Mathematica*[17]) può dare efficaci stimoli in questo senso.
Il punto K che abbiamo costruito nella **Figura 12** fornisce una *divisione* o *sezione* del segmento AB in due segmenti AK e KB. La proprietà caratteristica della divisione determinata dal punto medio D è banale: BD e DC sono congruenti. Nel laboratorio a l'insegnante pone il problema di scoprire la proprietà caratteristica della divisione di AB determinata da K. Gli studenti cercano di trovare una relazione tra le lunghezze, presumibilmente senza riuscirci e l'insegnante, dopo aver suggerito di guardare ai "rapporti" conduce la classe a verificare (utilizzando gli strumenti messi a disposizione da GeoGebra) che la proprietà caratteristica è

$$AB{:}BK = BK{:}KA.$$ [18]

[15] Nel laboratorio con gli studenti universitari si è discusso sul come guidare gli studenti di scuola a scoprire questa costruzione. Non ci interessa in questa sede riportare la discussione ma solo di dare un'idea del tipo di lavoro che ci sembra utile fare in tutti corsi universitari dedicati ai futuri insegnanti e non solo in quelli di Didattica della matematica.

La divisione di AB determinata da K è considerata da Euclide in (Euclide II.11) (dove la costruzione di K è sostanzialmente quella che viene fatta riscoprire nel laboratorio) e in (Euclide VI.30), dove Euclide, a differenza di quanto fatto nel libro II, impiega la teoria delle proporzioni per dividere un segmento in *media ed estrema ragione*. La sezione determinata da K prende il nome di *sezione aurea* e il rapporto $AB/BK = BK/KA$, *rapporto aureo*.[19] Quindi, il segmento BK è il lato del decagono regolare. Una possibile dimostrazione di questo fatto segue dalle proprietà del triangolo ABI. Attorno all'esplorazione delle proprietà di questo triangolo viene quindi incentrata l'attività del laboratorio, in cui lo studente, usando gli strumenti di misura di GeoGebra viene condotto a congetturare il fatto, dimostrato in (Euclide IV.11), che gli angoli alla base sono doppi rispetto all'angolo al vertice.
Il laboratorio prosegue quindi con l'esplorazione delle proprietà del pentagramma regolare e dei suoi rapporti con il pentagono regolare, ma crediamo che quanto detto sia sufficiente per dare l'idea dell'attività e dei suoi scopi.
Le proprietà principali del pentagono e del pentagramma regolare sono studiate da Euclide in una serie di proposizioni che, pur essendo molto ben organizzata, non segue necessariamente la strada più naturale e intuitiva. Euclide sa già dove vuol arrivare mentre, nel laboratorio, abbiamo cercato di tracciare una via plausibile per la scoperta. D'altra parte, come è noto, negli *Elementi*, Euclide presenta solo la sintesi dei suoi teoremi, cioè la dimostrazione che ciò che si afferma è vero, non la via per scoprire una verità matematica. Secondo quanto sappiamo da Pappo i matematici greci avevano elaborato una strategia per ritrovare quanto cercato, una sorta di problem solving. Tutto questo era presentato in quei lavori che Pappo presenta come "Il tesoro dell'Analisi" una serie di testi per la maggior parte perduti. Anche se il metodo di Analisi è assai diverso da quello proposto dalla didattica contemporanea, la differenza tra i due momenti, quello della scoperta e quello della dimostrazione, resta centrale[20].

---

[16] Il tema è ampiamente trattato nella letteratura della didattica della matematica, cfr. p.e. (Mariotti e Maffia 2018).
[17] https://www.wolfram.com
[18] Dal punto di vista della teoria delle proporzioni, si tratta di una caratterizzazione molto semplice delle proprietà del punto K. Crediamo che sia opportuno dedicare spazio alla formazione di un "pensiero proporzionale" nei futuri insegnanti, dedicando allo scopo, per esempio, alcune lezioni del corso di Storia della Matematica, perché, in diversi contesti di geometria elementare, come questo relativo alla caratterizzazione del punto K, la ricerca della proporzionalità è senza dubbio la via più semplice e naturale (una volta acquisita) per motivare un risultato (teorema o costruzione). Altrettanto importante, a nostro avviso, è la discussione dei limiti del pensiero proporzionale, della necessità di completarlo con un pensiero algebrico più raffinato e delle difficoltà (storiche e didattiche) che si devono superare per costruire il "pensiero algebrico".
[19] Ripercorrendo la soluzione della proporzione con linguaggio algebrico avremmo, detto x=BK=BJ e a=AB, l'equazione $x^2 + ax - a^2 = 0$ e quindi $x = BJ = \sqrt{a^2 + \left(\frac{a}{2}\right)^2} - \frac{a}{2}$; che corrisponde esattamente alla costruzione euclidea. La manipolazione algebrica ulteriore, svolta solitamene nei manuali, che porta a $x = a\frac{\sqrt{5}-1}{2}$ rende del tutto opaca l'interpretazione geometrica.

Nella dimostrazione costruita in laboratorio della correttezza della costruzione del pentagono regolate si parte dal triangolo isoscele ABG, si osserva che gli angoli alla base sono doppi rispetto a quello al vertice e, ammettendo che la somma degli angoli in ogni triangolo sia pari a due angoli retti[21], si deduce che la costruzione del decagono regolare, e quindi quella del pentagono regolare, è esatta.

L'uso di GeoGebra ha tratto fin qui solo parziale vantaggio dal carattere "dinamico" del software (per esempio lasciando variabile il punto K nella figura 13 per esplorare il problema della chiusura della spezzata). Risulta invece cruciale per formare e sviluppare l'intuizione geometrica nell'affrontare i problemi (ad esempio attraverso l'uso critico dei diversi strumenti) e nel comunicare le scoperte o i controesempi ai partecipanti il laboratorio. L'aspetto dinamico sarà invece pienamente sfruttato nelle attività descritte nelle sezioni seguenti.

# 4. L'inversione circolare

Le attività che presentiamo in questa sezione riguardano le trasformazioni geometriche: non le isometrie e le dilatazioni, che si considerano nella geometria elementare ma l'inversione circolare e le sue generalizzazioni. Il percorso ci è stato in parte suggerito dallo studio dei lavori di Giusto Bellavitis (1803-1880) (Bellavitis 1838) e di Thomas Archer Hirst (1830—1892) (Hirst 1865), che sono stati analizzati da uno degli autori nella sua tesi di dottorato (Raspanti, Tesi di dottorato 2017) e in due successivi lavori (Raspanti, Dall'inversione circolare all'inversione quadratica: aspetti storici e potenzialità didattiche 2016) (Raspanti, Giusto Bellavitis e la sua geometria di derivazione s.d.). Le attività sono state proposte nei laboratori del Corso monografico di Storia della Matematica (Sapienza, Università di Roma), in particolare nell'a.a. 2018-19 (15 studenti) e 2019-20 (10 studenti).

L'inversione circolare può collegare diverse parti della geometria, sia come punto d'arrivo di un percorso di studio della geometria elementare che si potrebbe svolgere alle scuole secondarie di secondo grado, sia come punto di partenza di un percorso di formazione per insegnanti e futuri insegnanti, che mette in luce i collegamenti con la geometria proiettiva e le geometrie non euclidee. Tratteremo l'argomento da questo secondo punto di vista, senza tralasciare alcuni commenti relativi alla possibilità di costruire percorsi di approfondimento nella scuola.

[20] Il possibile uso della visualizzazione nell'elaborazione della matematica antica è attualmente oggetto di studio da parte di numerosi storici e filosofi della scienza. Cfr. ad esempio (Saito 2012) (Netz 1999) (Mumma e Panza 2012).

[21] Questa proprietà è facile da dedursi in maniera convincente a partire da fatti semplici e quindi con ogni probabilità era ben nota fin dagli albori della matematica greca. Si osservi però che le deduzioni antiche, quelle che si attribuiscono, per esempio alla scuola pitagorica, non sono completamente assimilabili alle dimostrazioni euclidee. I postulati, posti da Euclide a fondamento della geometria erano sconosciuti ai primi matematici greci. Ogni dimostrazione dipendeva da un insieme di fatti evidenti che potevano essere continuamente modificati e criticati. Si tratta di una forma di deduzione più vicina al metodo dialettico di Platone, che si affermava nella ricerca filosofica, piuttosto che a quello assiomatico di Euclide. Questo modo di argomentare è anche più vicino alla logica della scoperta matematica secondo Lakatos ed è questo il modello che vogliamo ricreare nelle attività laboratoriali proposte, perché è più naturale e stimolante e in definitiva più "comprensibile" per chi si avvicina alla matematica. Crediamo che per poter impiegare un siffatto metodo di insegnamento in maniera efficace sia necessaria una formazione del futuro insegnante che preveda momenti di ricostruzione di percorsi matematici standard, come quello suggerito in questo paragrafo, e momenti di esplorazione di problemi matematici più avanzati, come diremo nel successivo paragrafo

## 4.1 La proiezione stereografica

Nello spazio tridimensionale la *proiezione* da un punto P di un piano $\pi$ ad un piano $\pi'$ (entrambi non contenenti P) si definisce associando ad un punto $X \in \pi$ il punto $X' =< P, X > \cap \pi'$, dove $< P, X >$ indica la retta passante per i punti P e X.

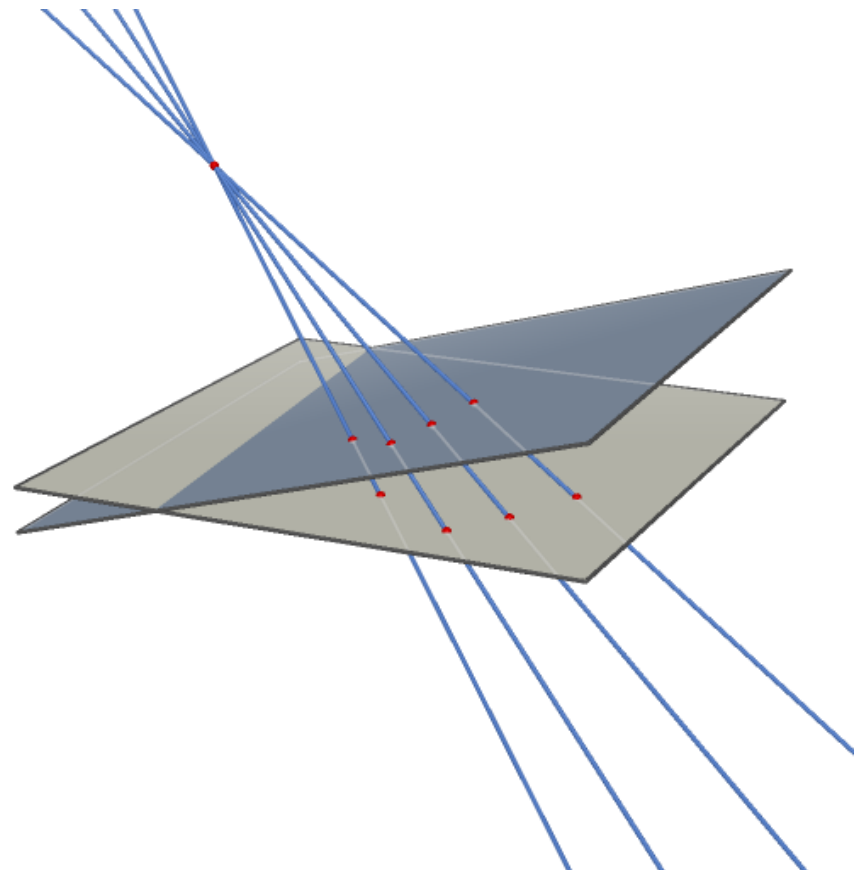

***Figura 14.*** *Proiezione da un punto di un piano su un altro.*

Nello stesso modo è possibile definire la proiezione tra figure più generali di una coppia di piani. Per esempio, la *proiezione stereografica*, proietta una sfera da un suo punto N, sul piano per il centro O della sfera, ortogonale al vettore NP. GeoGebra permette una visualizzazione tridimensionale della proiezione stereografica $X' \rightarrow X$ e della sua inversa $X \rightarrow X'$ (parte destra della **Figura 15**) e la contemporanea visualizzazione dei punti immagine sul piano su cui proiettiamo (parte sinistra della **Figura 15**).

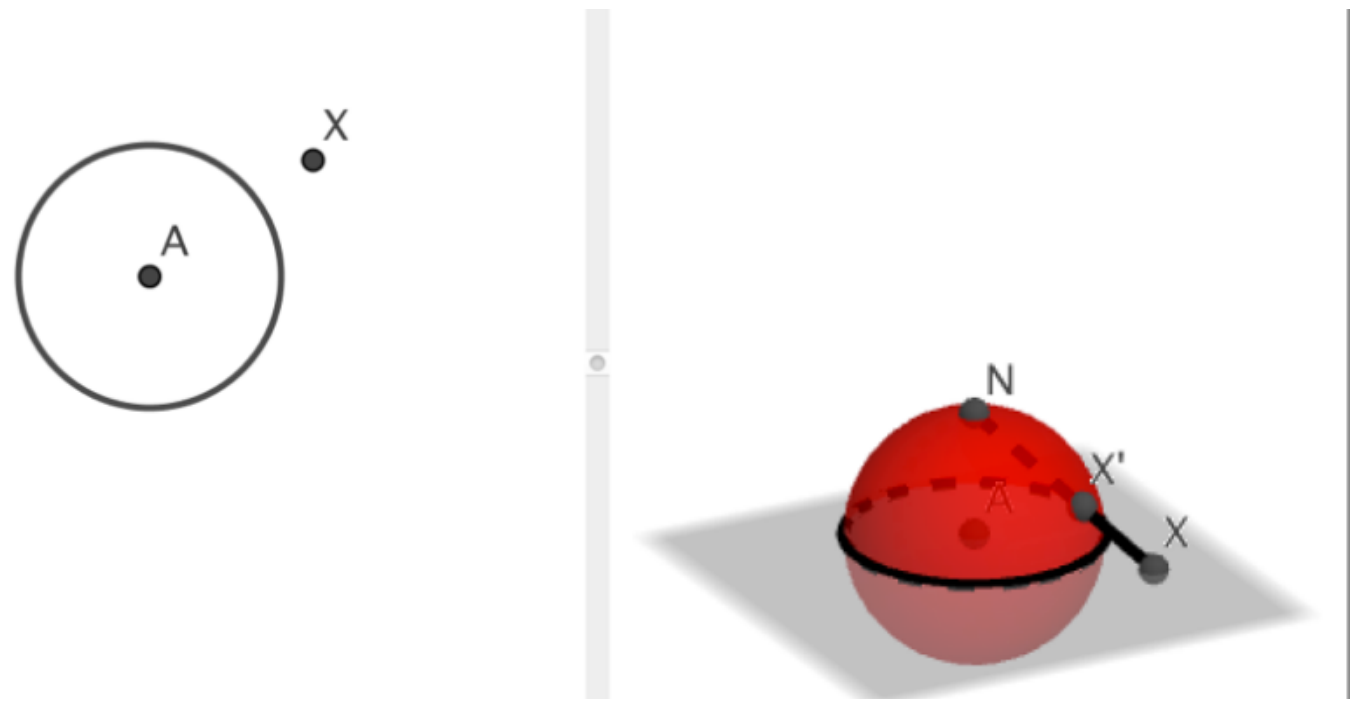

***Figura 15.*** *In GeoGebra è possibile avere la vista contemporanea della proiezione stereografica (a destra) e della sua traccia sul piano orizzontale (a sinistra). Abbiamo osservato l'utilità di mantenere entrambe le viste nella fase iniziale dell'esplorazione e in particolare nella considerazione della composizione di una proiezione stereografica con la sua inversa, illustrata nella* ***Figura 21****. La costruzione tridimensionale a destra viene realizzata in laboratorio dai partecipanti mentre la restrizione al piano orizzontale che GeoGebra mostra automaticamente nella vista bidimensionale.*

La proiezione manda rette in rette, trasforma cerchi in coniche e deforma distanze e angoli.

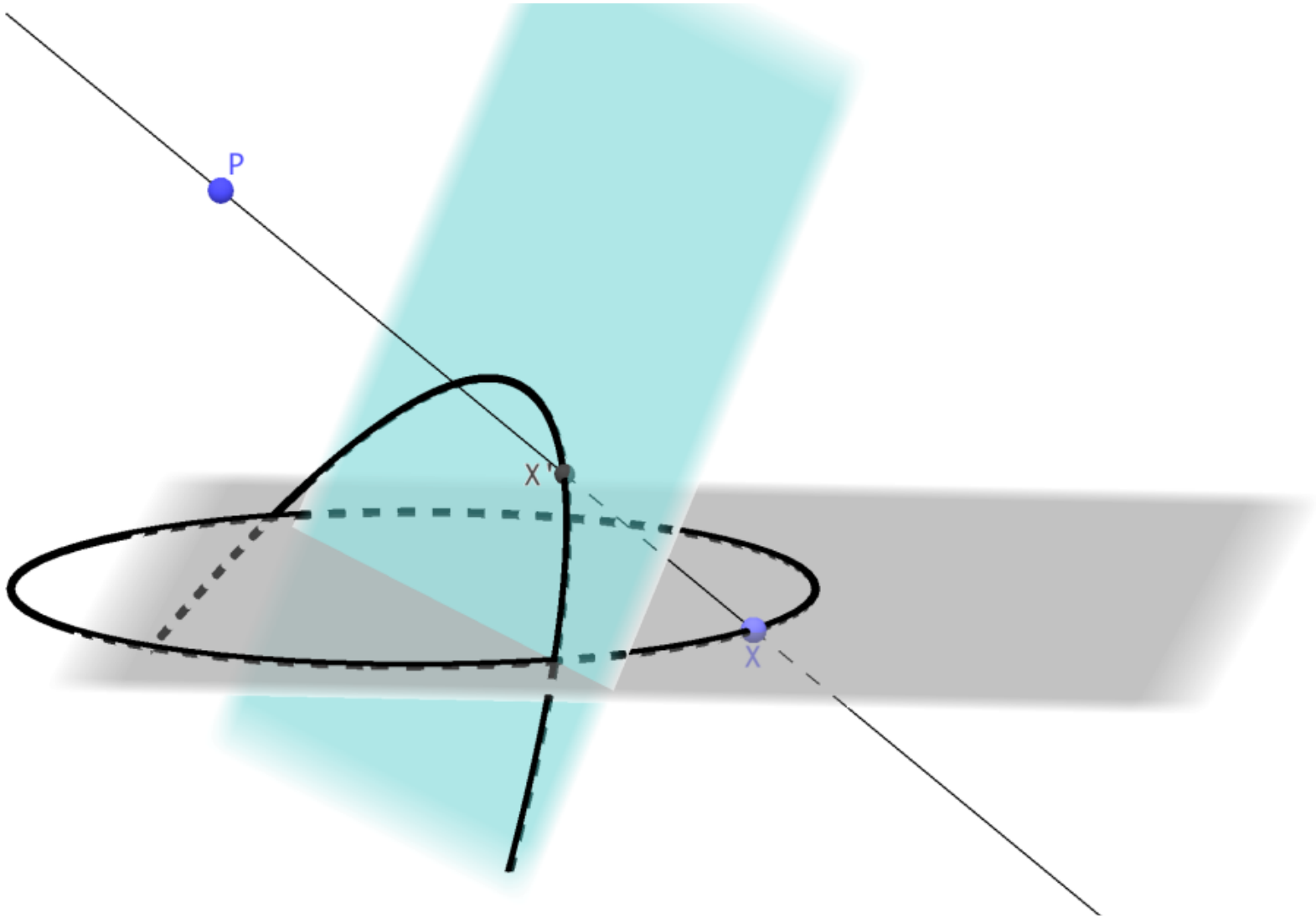

***Figura 16.*** *Ogni conica si può ottenere come proiezione di una circonferenza.*

Quali sono le proprietà della proiezione stereografica? Nel laboratorio sulle trasformazioni, abbiamo cominciato con il porre questa questione nella forma più precisa: quali proprietà della proiezione stereografica potete osservare studiando con GeoGebra come si trasformano le rette e le circonferenze?

Qual è il valore formativo di un'attività del genere per i futuri insegnanti?

1. La possibilità di costruire facilmente esempi, porta ad apprezzare il valore nell'apprendimento e nell'insegnamento della matematica dell'elaborazione di schemi di organizzazione sistematica delle osservazioni.
2. Il contesto laboratoriale in cui si svolgono le attività permette di porre domande aperte molto generali, quali "stabilisci le proprietà di una certa costruzione" e di apprezzarne il valore nell'apprendimento/insegnamento della matematica come strumento per stimolare un approccio creativo nell'affrontare i problemi.
3. I momenti di condivisione e di discussione degli esiti delle esplorazioni permettono di apprezzare il valore della discussione e del confronto critico nella costruzione della conoscenza matematica.
4. La possibilità di costruire con il software meccanismi digitali sofisticati[22] per realizzare oggetti geometrici dinamici permette di aprire diversi "punti di vista" su un problema di carattere geometrico e di sviluppare l'*intuizione geometrica*.

---

[22] Meccanismi digitali sono, per esempio, quelli considerati in (Bartolini Bussi e Maschietto 2006), che realizzazione con un software di Geometria Dinamica le macchine matematiche descritte nel libro. È possibile progettare e costruire "Macchine digitali" più generali anche senza che debbano esistere corrispettive macchine reali. Un primo esempio elementare è quello del meccanismo digitale per illustrare la costruzione tridimensionale di una trasformazione bidimensionale e contemporaneamente la sua visualizzazione bidimensionale, come nella **Figura 15** e più significativamente nella **Figura 21**. Nella costruzione delle macchine digitali è fondamentale la progettazione dell'interfaccia con cui l'utente può interagire con la Macchina. La progettazione delle interfacce offre numerose opportunità di attività importanti sia dal punto di vista didattico, sia da quello della formazione professionale.

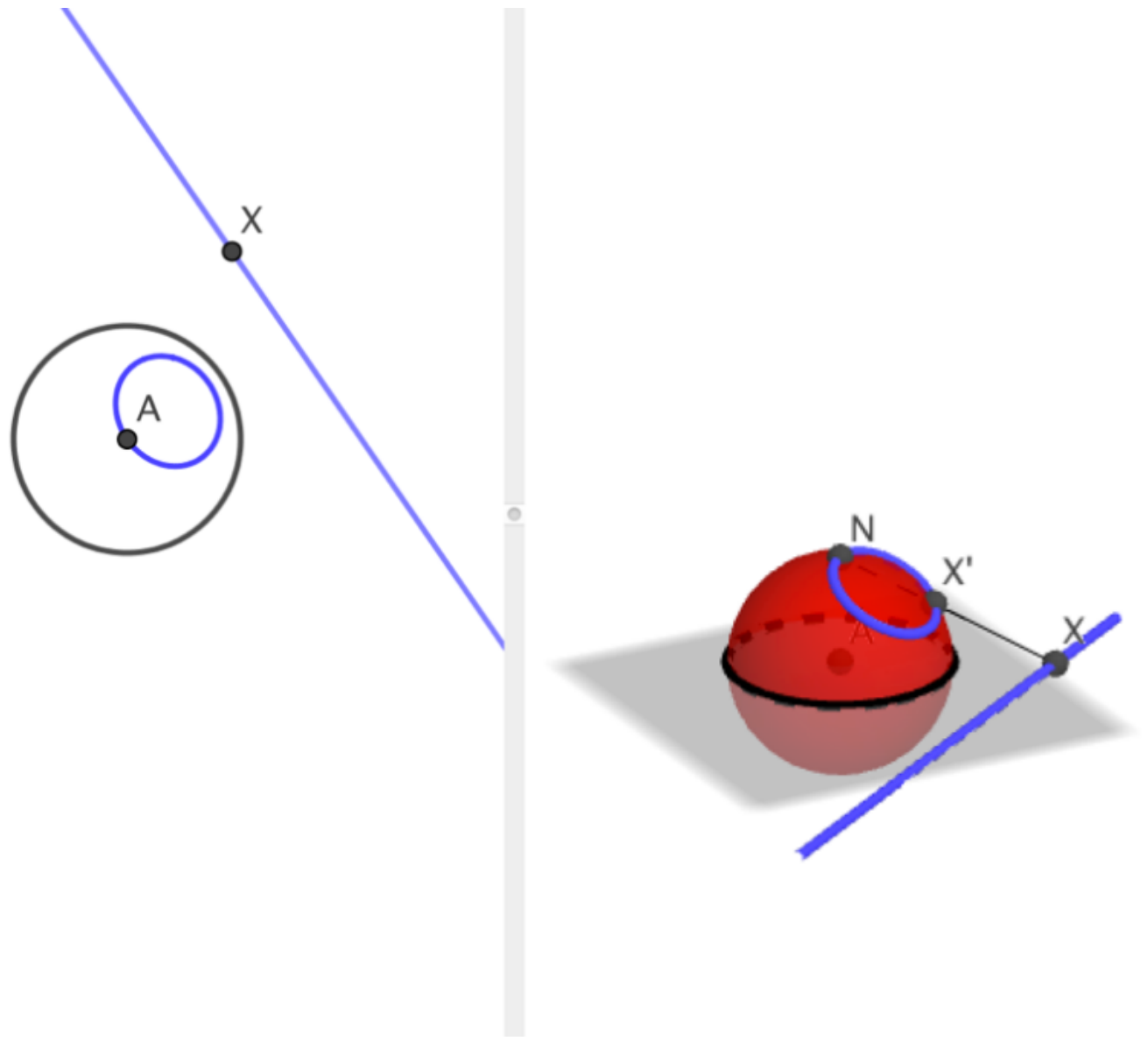

*Figura 17. La proiezione stereografica manda cerchi per N in rette.*

Ecco alcune delle proprietà osservate nel laboratorio.
i) La proiezione stereografica manda cerchi per N in rette e manda cerchi che non passano per N in cerchi.

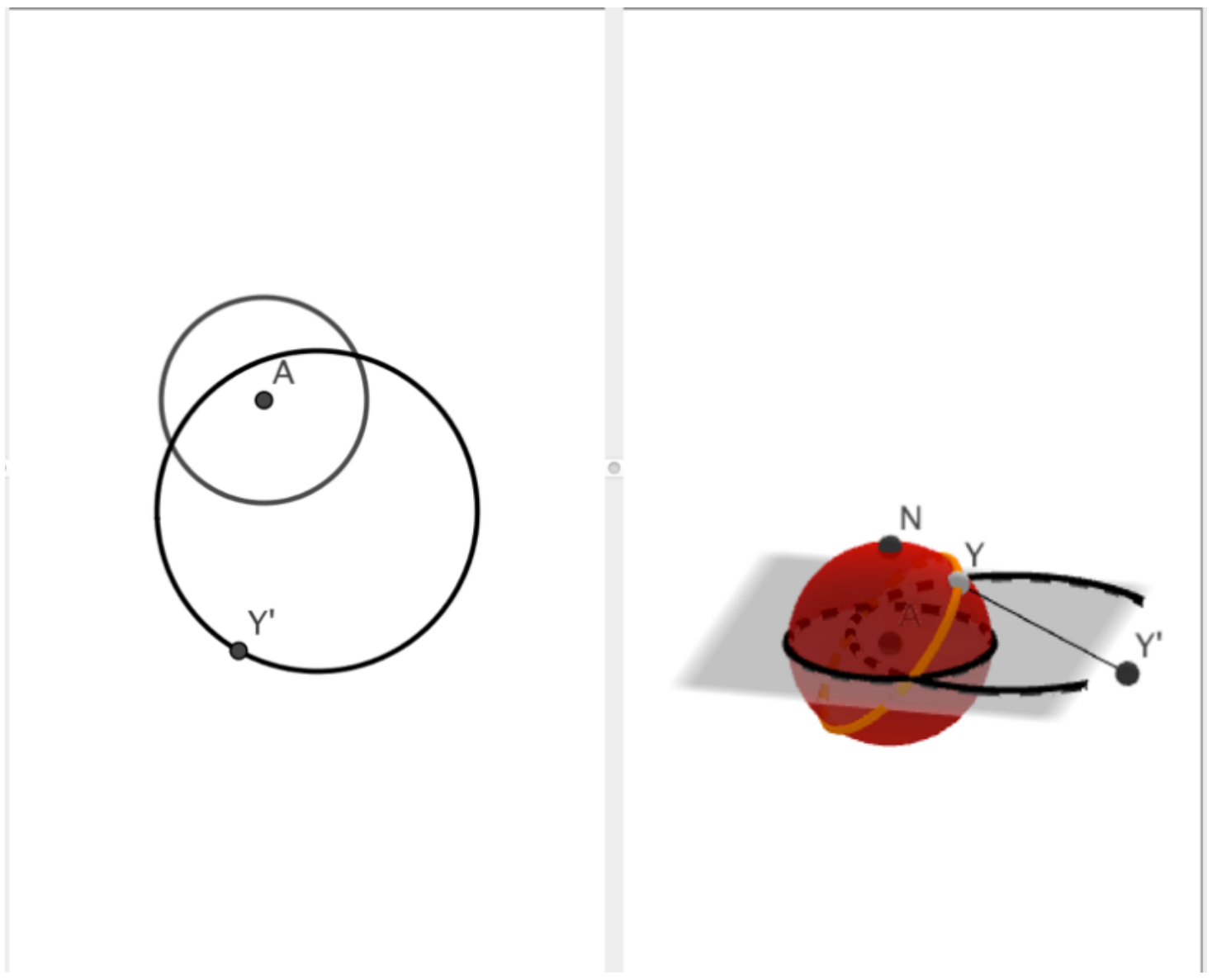

*Figura 18. La proiezione stereografica manda cerchi che non passano per N in cerchi.*

ii) Le immagini dei cerchi massimi sono cerchi che hanno centro in un punto I qualsiasi del piano e passano per le intersezioni H e G della perpendicolare condotta dal centro A del cerchio "base" alla congiungente AI.[23]

[23] Usando l'inversione circolare che introdurremo nel paragrafo seguente, possiamo caratterizzare questi cerchi come quelli che se passano per un punto P passano anche per il punto Q simmetrico rispetto ad A dell'inverso di P rispetto all'equatore. Prendendo come rette le immagini dei cerchi massimi e come cerchi i cerchi euclidei (ma con centro "sferico" diverso dal centro euclideo) otteniamo un modello piano della geometria sferica.

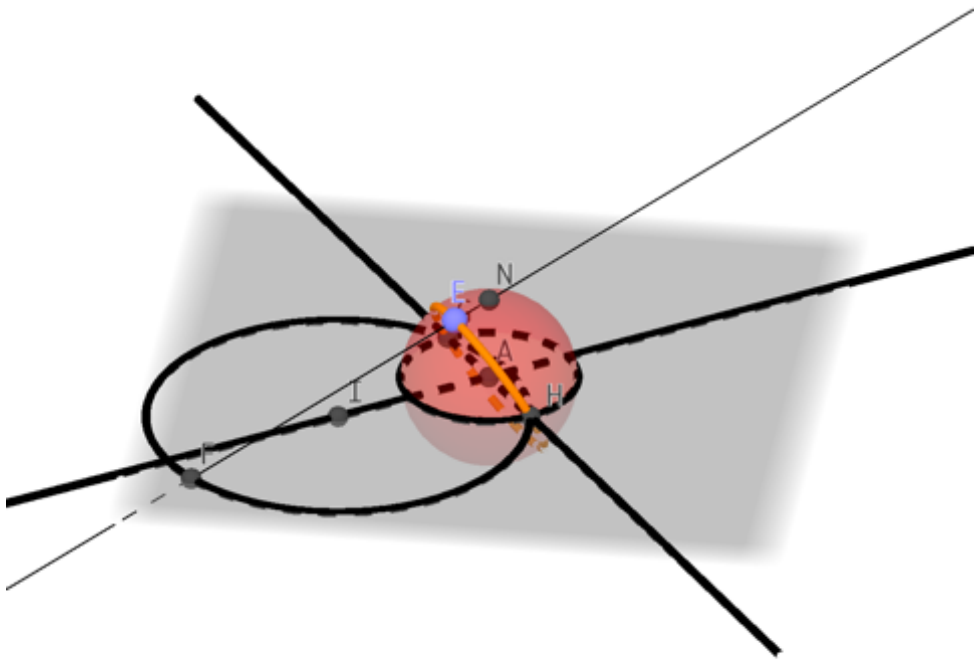

*Figura 19. La proiezione stereografica di un cerchio massimo.*

iii) La proiezione stereografica deforma le distanze, ma preserva gli angoli.

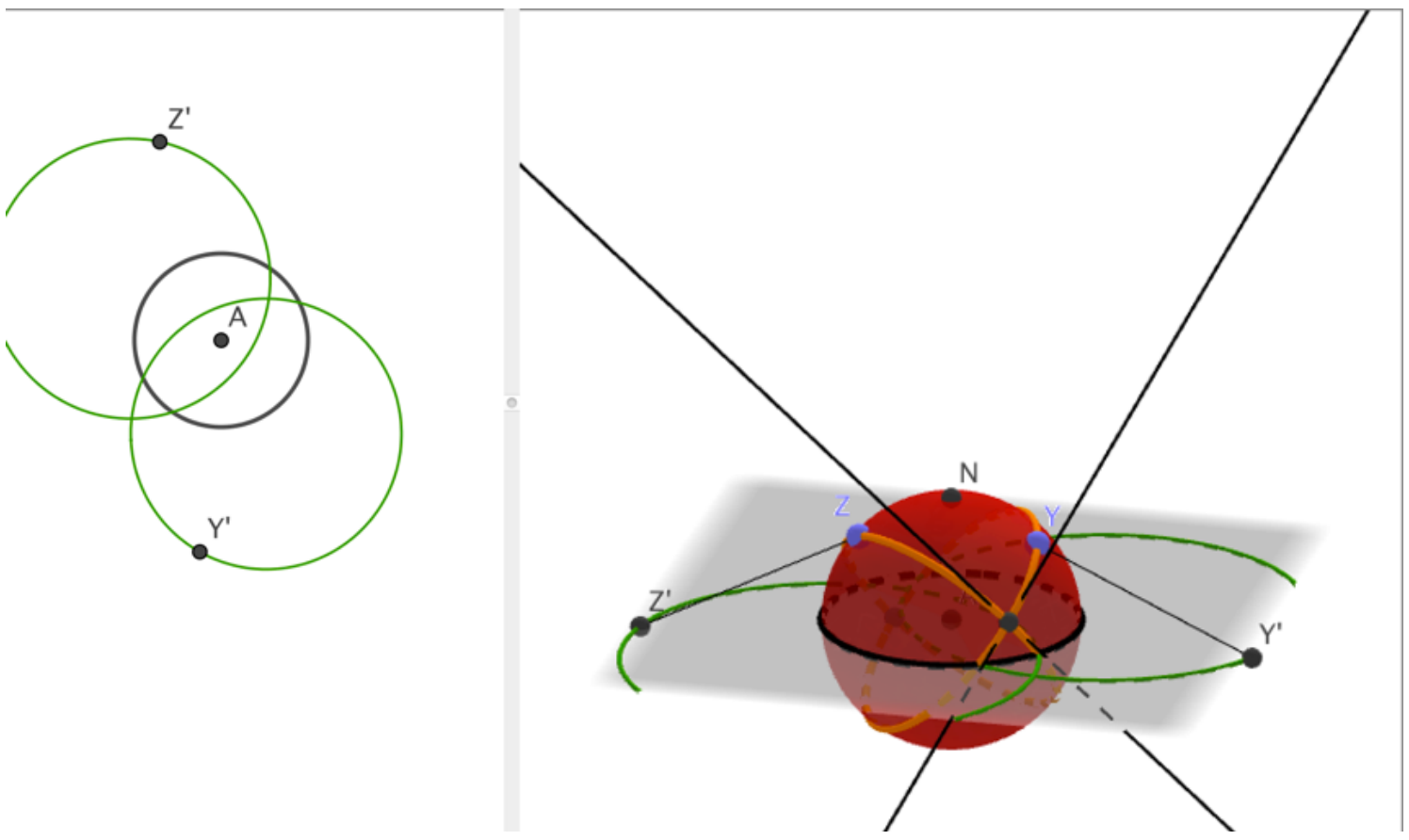

*Figura 20. Proprietà conforme della proiezione stereografica.*

Queste proprietà, osservate durante le attività di laboratorio, non erano precedentemente note ai partecipanti, L'insegnante ha lasciato liberi i partecipanti di esplorare le proprietà della proiezione stereografica Ogni gruppo è stato in grado di osservare almeno alcune delle proprietà che abbiamo riportato.

## 4.2 L'omologia e l'inversione circolare

Componendo due proiezioni (proiezioni da un punto o proiezioni stereografiche) si ottengono trasformazioni più generali.

Scelti due punti P e Q e due piani $\pi$ e $\pi'$, la composizione della proiezione da P di $\pi$ su $\pi'$ con quella da Q di $\pi'$ su $\pi$ è un'*omologia*. È interessante far scoprire con GeoGebra le principali proprietà dell'omologia da questa descrizione, ma non ci soffermeremo su questo punto. Componendo l'inversa della proiezione stereografica da N con la proiezione stereografica dal punto antipodale S otteniamo l'*inversione circolare*.[24]

[24] Consideriamo sulla sfera la riflessione f rispetto a un piano passante per il centro A e Q l'immagine stereografica di un punto P e Q' quella di f(P). Q' è l'inverso di Q rispetto alla circonferenza massima tagliata dal piano nella sfera. Questa corrispondenza significa che le isometrie euclidee corrispondono stereograficamente nel piano a prodotti di inversioni rispetto a circonferenze massime.

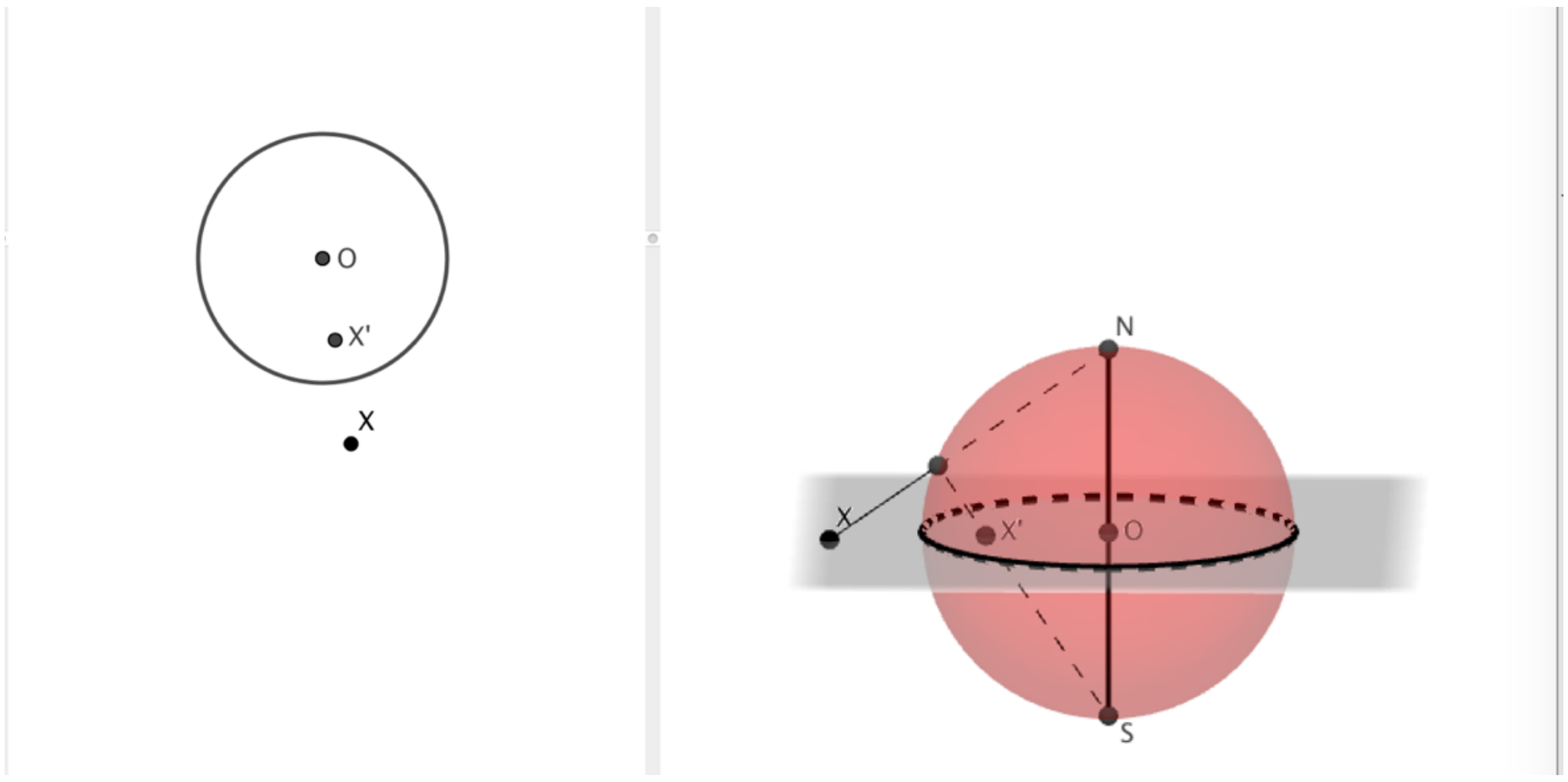

***Figura 21.*** *L'inversione circolare (a sinistra), illustrata come composizione di una proiezione stereografica con l'inversa della projezione stereografica dal punto antipodale. Il meccanismo digitale raffigurato nell'illustrazione permette di collegare la trasformazione piana alla sua costruzione tridimensionale.*

L'inversione circolare suggerisce esplorazioni laboratoriali ricche di spunti di approfondimento e di collegamenti. Si tratta di una trasformazione importante in vari ambiti della matematica, anche in ambiti elementari, tanto da essere implementata nel menu delle costruzioni di GeoGebra. Nel laboratorio useremo senz'altro questa implementazione, abbandonando la costruzione tridimensionale che abbiamo usato fin qui. Per convincersi, senza scrivere la forma analitica e senza argomentare geometricamente, che quella implementata in GeoGebra e quella costruita componendo le due proiezioni stereografiche coincidono, basta considerare, come suggerito nel laboratorio, l'immagine X'' di X con la prima nella vista bidimensionale del meccanismo digitale illustrato nella **Figura 21** e osservare che, muovendo il punto X, i punti X' e X'' risultano sempre sovrapposti.

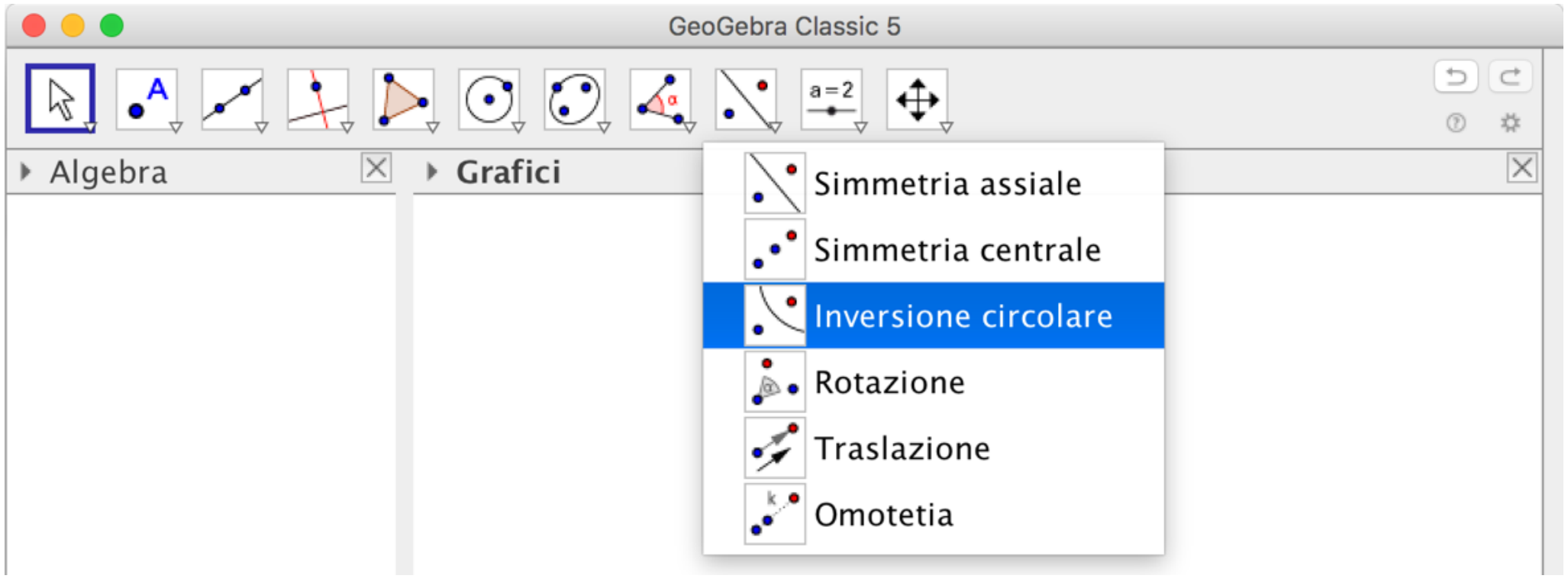

***Figura 22.*** *L'inversione circolare fa parte del menu delle trasformazioni di GeoGebra.*

Per cercare le proprietà dell'inversione circolare abbiamo posto agli studenti le seguenti domande: i) Con riferimento alle trasformazioni euclidee (Traslazioni; Rotazioni; Riflessioni; Omotetie) elenca, per ogni classe di trasformazioni, una lista di proprietà che giudichi interessanti e caratteristiche. ii) Indaga, utilizzando GeoGebra, se e quali delle proprietà che hai osservato al punto precedente sono verificate dell'inversione circolare. iii) È possibile modificare alcune delle proprietà osservate in maniera che siano verificate dell'inversione circolare? iv) Puoi osservare ulteriori proprietà dell'inversione circolare che ti sembrano degne di nota?

Tra le diverse proprietà che sono state osservate e successivamente studiate nei laboratori citiamo: la determinazione dell'immagine di rette e cerchi[25]; la determinazione della forma analitica dell'inversione; il fatto che ogni cerchio attraverso $P$ e $P'$ interseca ortogonalmente il cerchio di inversione e che un cerchio che interseca ortogonalmente il cerchio di inversione resta invariato per inversione; l'inversione è una mappa conforme.
È stato posto infine il problema di caratterizzare l'immagine delle coniche, senza la pretesa che si arrivasse a risultati esaurienti ma solo per rendere conto delle possibilità del software. Si è poi introdotta la nozione di birapporto e si è fatto notare, usando gli strumenti metrici del software, che se $C$ e $C'$ sono inversi rispetto ad un cerchio $\Gamma$, i punti di intersezione $\{A, B\}$ di $\Gamma$ con la retta $CC'$ sono coniugati armonici di $\{C, C'\}$. Questo ci ha permesso di esplorare i collegamenti con la geometria proiettiva, come diremo nel successivo paragrafo.

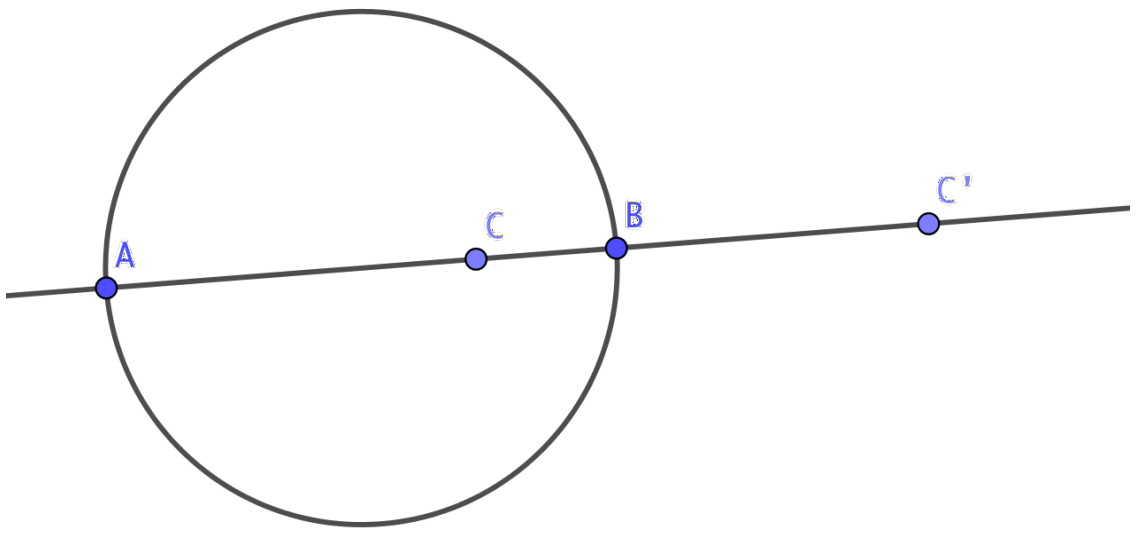

***Figura 23.*** *Se C' è l'inverso circolare di C, la quaterna ABCC' è armonica.*

# 4.3 La costruzione del quarto armonico e la geometria proiettiva

La costruzione del quarto armonico attraverso l'involuzione circolare, che abbiamo discusso al termine del paragrafo precedente, permette di esplorare questioni interessanti e stimolanti di geometria proiettiva.
La geometria proiettiva è quasi completamente scomparsa dal curriculum dei corsi di matematica, ma crediamo che possa ancora giocare un ruolo importante nella formazione del futuro insegnante, offrendo numerosi spunti per il trattamento "estensivo" di diversi temi di grande rilevanza storica e pedagogica. Essendo particolarmente adatta ad essere esplorata tramite un software di Geometria Dinamica, è possibile affrontare in maniera efficace diversi argomenti, senza la necessità di una trattazione sistematica ma permettendo una trattazione atta a metterne in luce i collegamenti. È inoltre possibile affrontare questi argomenti mettendo lo studente nelle condizioni del "ricercatore matematico". Si tratta di un'esperienza spesso preclusa ai futuri insegnanti. Vogliamo illustrare alcuni esempi di attività mirate a coltivare questa attitudine alla ricerca matematica con le modalità e nella misura in cui è utile al futuro insegnante.
Usando la costruzione del quarto armonico, che può essere memorizzata come un nuovo strumento nel menu degli strumenti[26], è naturale porre il seguente problema: *data una conica C e un punto P che non gli appartiene, ogni retta r per P che incontra la conica la sega in due punti A e B. Al variare di r tra le rette per P, qual è il luogo dei quarti armonici dopo A, B e P?* Si tratta di un segmento che identifica una retta. Abbiamo chiesto se esiste una costruzione, presente nel menu di GeoGebra che restituisce la stessa retta. Non è stato difficile scoprire che tale costruzione esiste: si tratta della *polare* di un punto rispetto alla conica.

[25] Le immagini di rette e cerchi sono rette e cerchi, ma alcune rette vengono trasformate in cerchi e alcuni cerchi in rette.
[26] Una volta terminata una costruzione è possibile creare un nuovo strumento dal Menu Strumenti, specificando gli elementi di input, quelli di output e il nome del nuovo strumento, che potrà essere utilizzato esattamente come gli altri già presenti nel software.

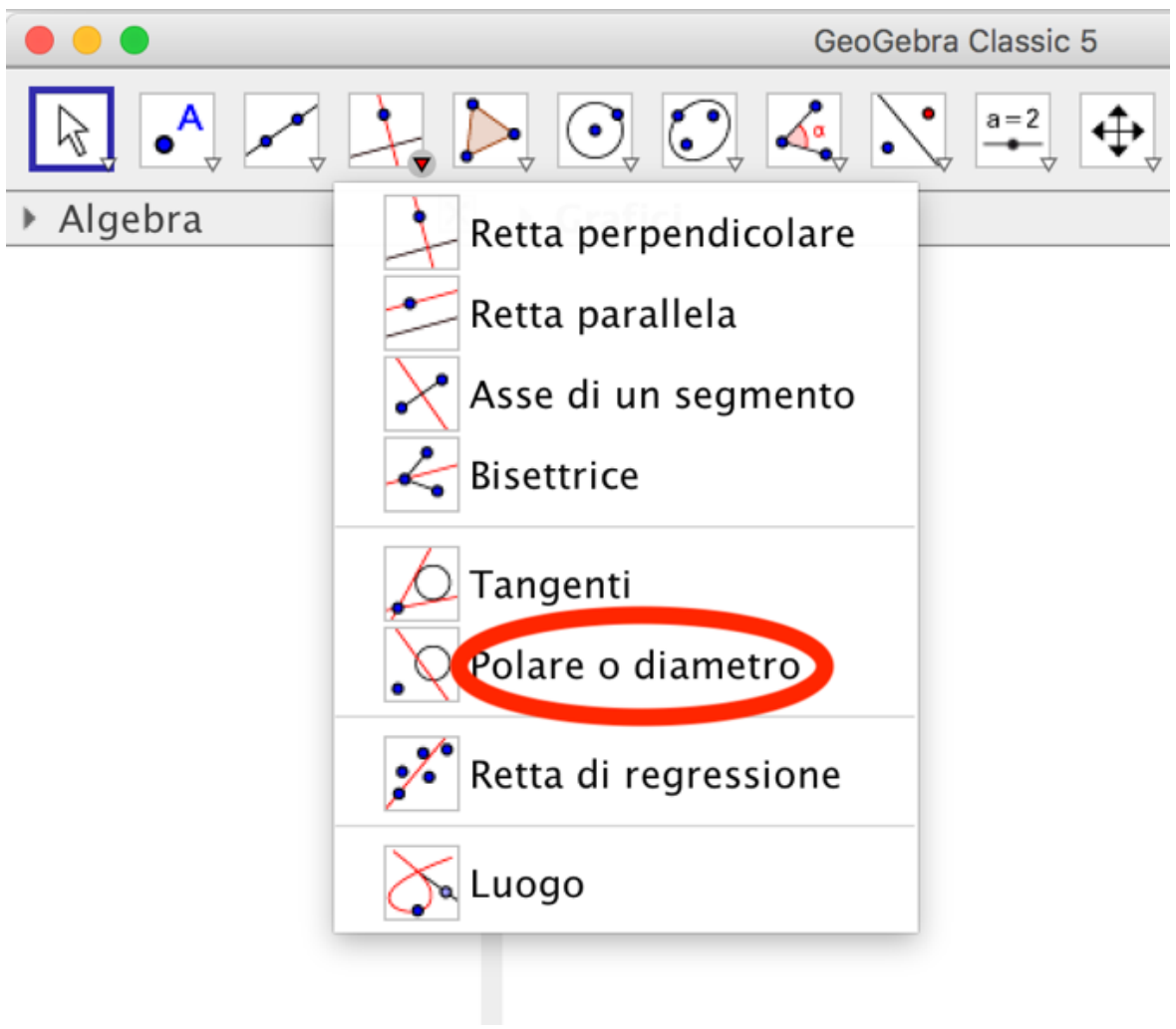

***Figura 24.*** *La costruzione della polare di un punto si trova nel menu delle costruzioni di GeoGebra.*

Chiedendo di disegnare, al variare della retta per P, il luogo del coniugato armonico di P rispetto all'intersezione della retta con la conica abbiamo così condotto i partecipanti alla scoperta di una costruzione della polare in termini di quarto armonico, che ~~ci~~ permette, volendo, di trovarne la forma analitica. Ma non ci interessa tanto questo spunto, quanto quello di studiarne le proprietà. Con GeoGebra abbiamo guidato la scoperta di queste proprietà in laboratorio, ponendo alcune domande: qual è la proprietà delle congiungenti di P con i punti di intersezione tra la polare e la conica? Se P varia su una retta, le polari di P cosa inviluppano? Se P varia su un cerchio, le polari di P cosa inviluppano? Quest'ultima domanda risulta particolarmente interessante da esplorare con GeoGebra e stimolante anche dal punto di vista teorico.

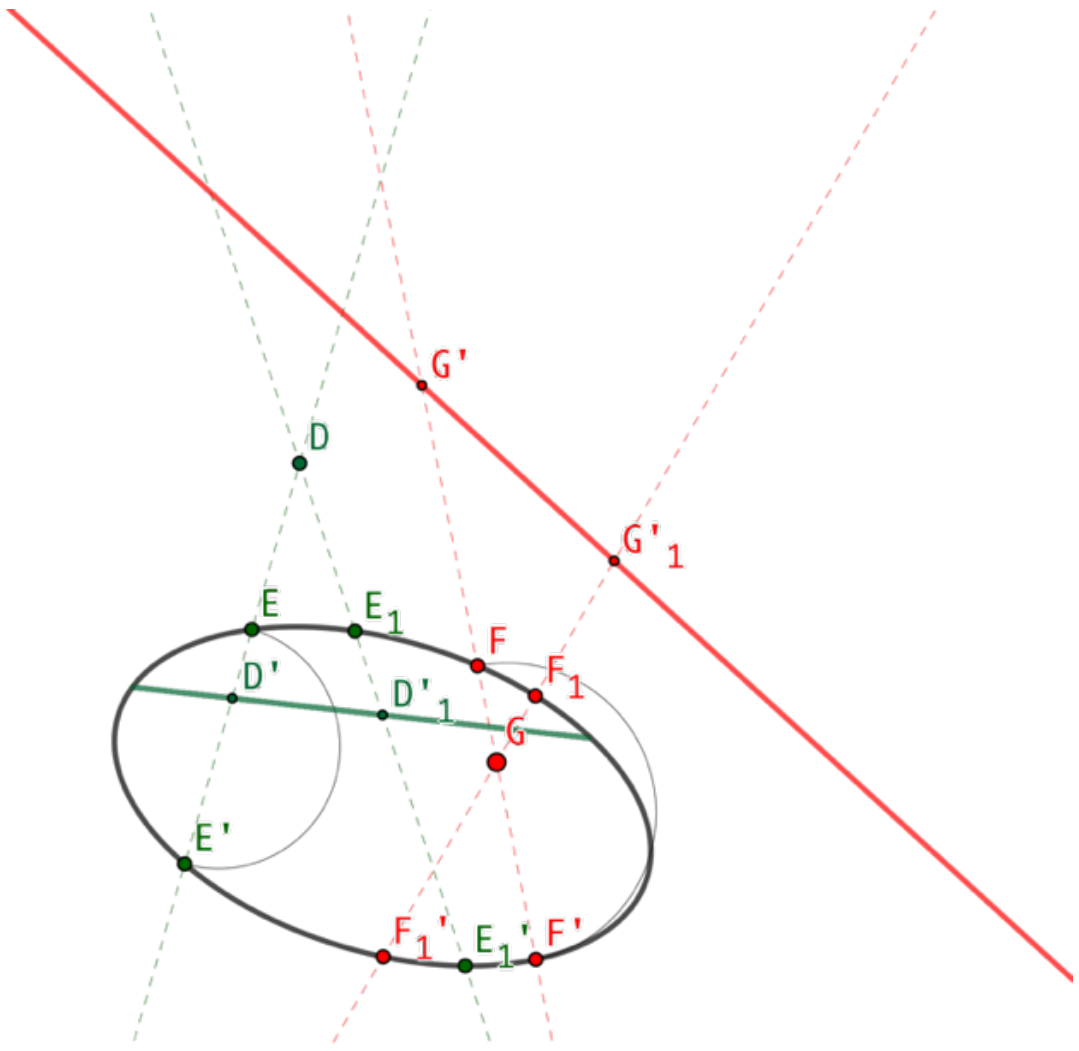

***Figura 25.*** *La costruzione tramite quarto armonico permette di definire la polare sia quando il punto è esterno all'ellisse, sia quando è interno ad essa. **

## 4.4 La costruzione della conica coniugata e delle secanti immaginarie

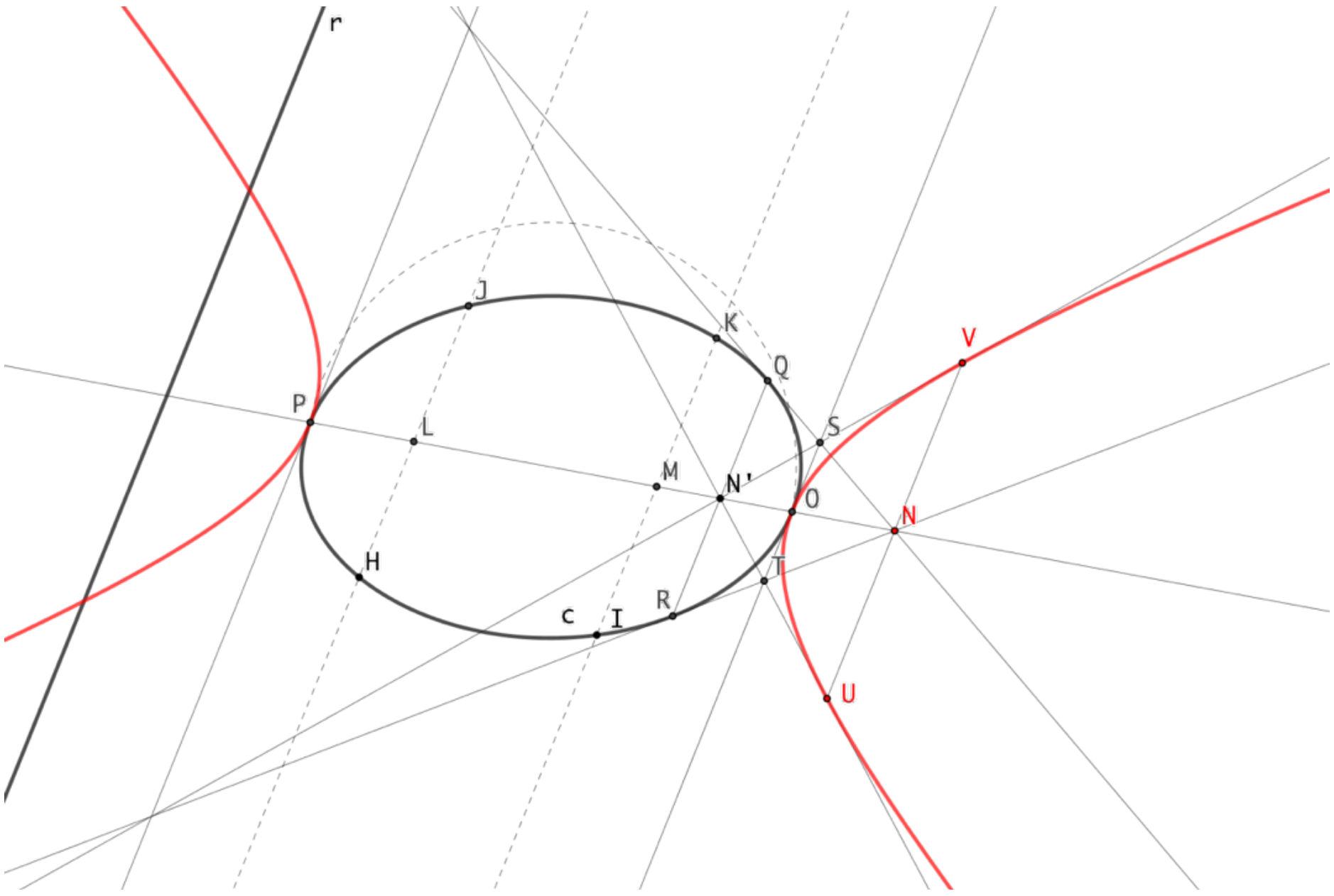

***Figura 26.*** *In rosso è disegnata la conica coniugata dell'ellisse HIJKQ rispetto al diametro LM. **

Una costruzione che riteniamo molto utile discutere ed esplorare è quella della *conica coniugata ad una conica data rispetto a una direzione*. Essa permette di visualizzare in maniera completamente geometrica gli elementi immaginari che intervengono nella teoria proiettiva delle coniche. La costruzione è descritta in (Bellavitis 1838) e crediamo che la implementazione della costruzione in GeoGebra attraverso la lettura delle parole di Bellavitis sia estremamente formativa per il futuro insegnante, in quanto si tratta di interpretare una costruzione descritta in un linguaggio diverso a quello a cui è abituato, in parziale analogia, pur con evidenti differenze, a quando l'insegnante dovrà interpretare il linguaggio impreciso dei suoi studenti. In questo genere di attività un software di Geometria Dinamica risulta particolarmente efficace in quanto implementa in maniera geometrica, e quindi vicina alla concezione originale, i principali strumenti necessari. Riportiamo per una costruzione della conica coniugata, ottenuta interpretando lo scritto di Bellavitis, cui abbiamo dedicato un laboratorio del Corso Monografico di storia della matematica.[27]

Data una conica c e una direzione r vogliamo costruire la conica coniugata a c rispetto alla direzione r. Scelti due punti H ed I di c tracciamo le corde HJ e IK parallele ad r. Siano L ed M i rispettivi punti medi. La congiungente LM è il diametro coniugato ad r. Siano O e P i punti di intersezione di tale diametro con c. Prendiamo un punto N esterno alla conica, cioè tale che la parallela ad r non intersechi la conica. Vogliamo costruire la *secante immaginaria* per N alla conica c. Sia N' il quarto armonico di N dopo O e P. Per costruirlo, come detto precedentemente, basta prendere l'inverso circolare di N rispetto al cerchio di diametro OP. Sia N' tale inverso. Si tracci da N' la corda QR parallela alla direzione r. Le congiungenti NQ ed NR intersecano la parallela a r per O (che è anche tangente in O alla conica c) nei punti S e T rispettivamente. Le congiungenti N'S ed N'T intersecano la parallela ad r passante per N nei punti V ed U rispettivamente. Il segmento UV è la corda ideale di c passante per N e parallela a r. Al variare di N sul diametro gli estremi della corda descrivono la conica coniugata a c rispetto alla direzione r.

Si osservi che la polare di N rispetto a c interseca la conica nei punti Q ed R che, congiunti a N determinano le tangenti a c in Q ed R. Analogamente, la polare di N' rispetto a c interseca la conica coniugata (rispetto alla direzione coniugata al diametro su cui sta N') nei punti U e V che, congiunti a N' determinano le tangenti alla conica coniugata in U e V.

[27] Sull'uso delle fonti nell'insegnamento della matematica, vedi anche (Jahnke, et al. 2000).

Guardando le cose dal punto di vista della geometria analitica possiamo osservare quanto segue. Ogni conica a centro ha come equazione, rispetto a un diametro e al suo diametro coniugato, $y^2 = p(r-x)(r+x)$ (ellisse) e $y^2 = p(x-r)(x+r)$ (iperbole) con $p$ e $r$ parametri ~~costanti~~. Due coniche con gli stessi parametri rispetto agli stessi diametri sono dette coniugate. È evidente che le ordinate di un punto variabile su uno dei diametri saranno reali se e solo se quelle della conica coniugata sono immaginarie. Questo suggerisce di usare la conica coniugata per definire, nel piano reale, oggetti immaginari (*ideali* nella terminologia di Poncelet, ripresa da Bellavitis).[28] Per esempio, dette $C$ e $C'$ una conica e la sua coniugata (supplementare), la corda ideale tagliata su $C$ da una retta che non la interseca, sarà quella tagliata su $C'$.[29]
Non proponiamo certo di tornare a usare teorie superate che permettono la rappresentazione reale degli enti complessi in sostituzione all'approccio algebrico, ma crediamo che possa essere ancora utile riflettere sulla possibilità di una visualizzazione geometrica di questi enti come stimolo didattico che permette di avvicinare alcune teorie alla loro origine storica per poterle meglio comprendere nelle loro relazioni con altri aspetti, non necessariamente matematici, della cultura dei tempi in cui sono stati prodotti, in modo da riconoscere e mettere in evidenza l'importanza culturale e non solo tecnica della matematica.[30] Vedi per esempio, con riguardo alla teoria degli spazi di dimensione maggiore di tre e alla geometria non euclidea (Rogora e Tortoriello, Matematica e cultura umanistica 2018).
La costruzione della conica coniugata permette di vedere, per esempio, la corda ideale comune a due coniche che non si intersecano. Poncelet, ad esempio, applica il concetto di *corda ideale* al caso di due circonferenze, ottenendo ciò che sarà poi detto *asse radicale*. Dice Poncelet: *Il résulte d'ailleurs … que les sécantes réelles et idéales communes doivent jouir absolument des mêmes*

---

[28] A tal proposito, è interessante, a nostro avviso, leggere come Bellavitis deduce le stesse relazioni qui espresse in forma analitica e come fa ricorso alla conica coniugata per definire la *secante ideale comune* a due coniche: "Se *AB* è un diametro della conica *AMB*, e parallelamente al suo diametro conjugato è condotta l'ordinata *PM*, è noto che $PM^2$ = *m.AP.PB*, essendo *m* un numero costante per tutta l'estensione del diametro *AB*. Risulta da ciò, che le condizioni necessarie perché due coniche abbiano i due punti comuni *M N* sono, che i diametri conjugati colla direzione della retta *MN* la incontrino in uno stesso punto *P*, e che inoltre sia *m.AP.PB=m'.A'P.PB'*. Ora queste due condizioni possono realizzarsi rispetto alla retta *TU* anche nella figura 14., nella quale quei due prodotti essendo negativi, rendono immaginaria la ordinata *PM*; e ad onta di quest'ultima circostanza la retta *TU* conserverà nella figura 14. Molte proprietà di cui essa godeva nella figura 13. Con questo scopo il Poncelet chiamandola *secante comune* nella figura 13., la disse poi *secante-comune ideale* nel caso della figura 14., indicando coll'epiteto *ideale* che sono immaginarii i punti *M N*, nei quali la retta *TU* è supposta tagliare ambedue le coniche" (Bellavitis 1838, §. 63., p. 266).

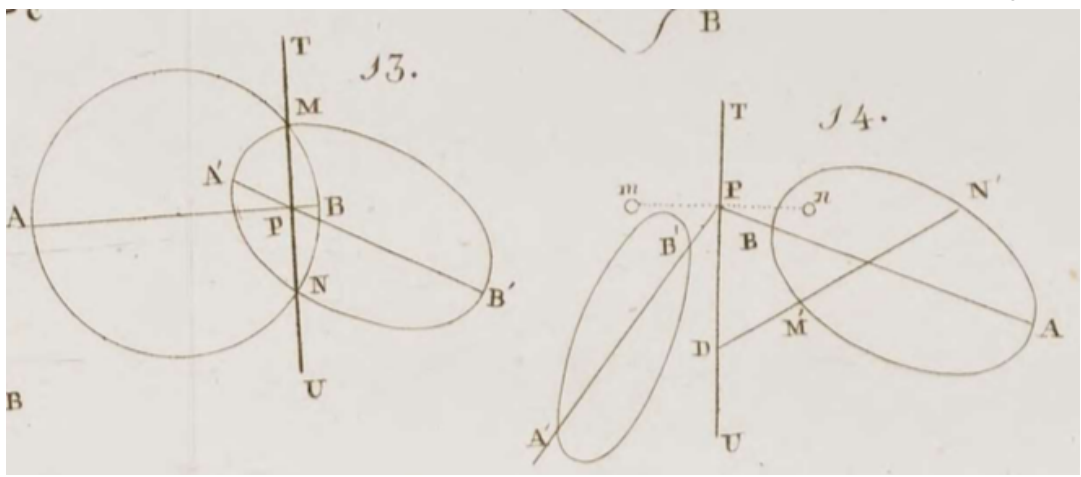

È interessante notare che Bellavitis esprime queste relazioni in un linguaggio non immediatamente comprensibile a chi è abituato ad un linguaggio più contemporaneo, ma che adeguatamente interpretato coincide con quanto espresso qui in forma più analitica; inoltre egli stesso deduce la definizione e le proprietà della *secante comune ideale* dal caso reale ricorrendo alla realizzazione grafica; in laboratorio il processo continuo di cui le figure 13 e 14 fissano staticamente due momenti è stato realizzato dinamicamente tramite GeoGebra, "muovendo" le due coniche in modo da passare dal caso in cui le loro intersezioni sono reali a quello in cui le loro intersezioni sono immaginarie.

[29] Si noti che, per determinare i punti V ed U nella parallela a r per N basta intersecarla con le congiungenti di Q ed R con O. Ciò dipende dal fatto che Q, O, U (come anche R, O, V) sono allineati.

[30] Senza contare che un approccio esclusivamente tecnico spesso disorienta gli studenti o una parte di essi. Mentre un approccio laboratoriale che consenta di *visualizzare* l'oggetto di studio, possibilmente in una prima fase di scoperta che preceda la fase di formalizzazione dei concetti, può servire ad abbattere la ritrosia di alcuni studenti davanti ad un approccio esclusivamente formale.

*propriétés ; mais c'est ce qu'on peut aussi démontrer pour le cas particulière de deux cercles, d'une manière fort tout à fait directe et forte simple*. Con Geogebra questa conservazione delle proprietà può essere verificata movendo uno dei cerchi. Si vedrà che non solo quando si intersecano la corda ideale si muta in quella reale, ma che la proprietà indicata da Poncelet (l'uguale lunghezza di due segmenti di tangente tratte da un punto della corda, non dipende dal fatto che essa sia reale o ideale, e tutto senza far uso alcuno dei numeri complessi.
Seguendo la descrizione della costruzione fatta da Poncelet in (J. Poncelet 1822) ci possiamo rendere conto che il software ci mette a disposizione uno strumento potentissimo, cioè la possibilità di determinare una costruzione "muovendo" con continuità alcuni oggetti. Si noti come questo modo di vedere è per Poncelet una straordinaria dote intuitiva che gli permette di "vedere" con la mente il movimento delle costruzioni, mentre, grazie al software, noi possiamo realmente vederle muoversi. Siamo convinti che questa pratica di "vedere il movimento" in contesti di ricerca sia molto utile per preparare l'insegnante a far uso efficace delle possibilità espressive del software nel momento in cui lo impiegherà a lezione, molto più della pratica relativa a casi elementari, in cui ha già raggiunto con l'intuizione la comprensione di un fatto che si può anche illustrare dinamicamente. L'utilizzo del software consente di "vedere il movimento" anche a chi non possiede la dote innata di "immaginare il movimento". In questo aspetto del software risiede un'enorme potenzialità didattica: allenare la mente a "vedere il movimento", per sviluppare capacità come l'intuizione e la generalizzazione

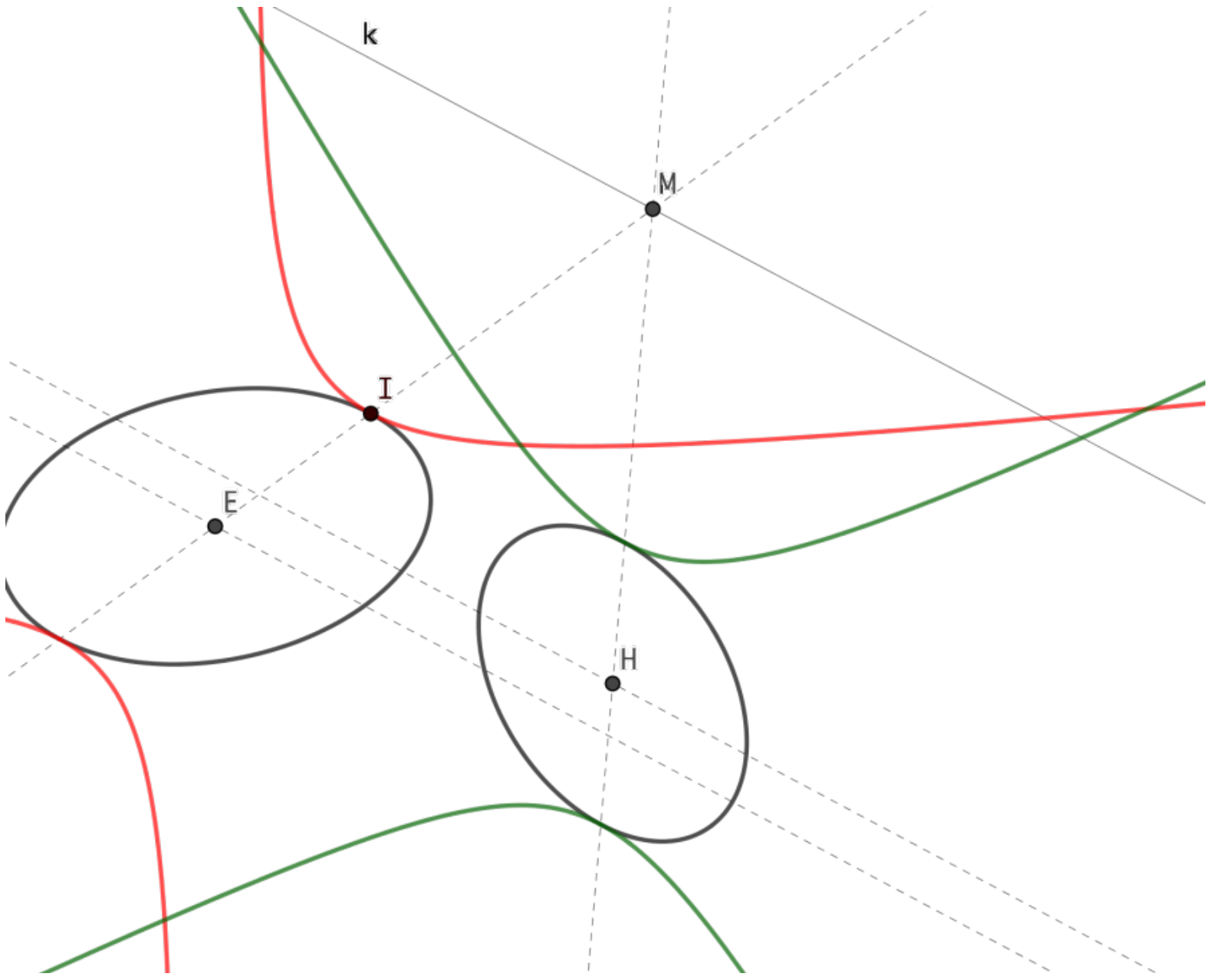

***Figura 27.** La "costruzione" della corda immaginaria comune a due coniche prive di punti di intersezioni reali. **

# 5. L'inversione di Bellavitis-Hirst

Le attività presentate in questa sezione sono state progettate per un laboratorio del corso di Istituzioni di Matematiche Complementari e non sono ancora state provate in classe. Questo laboratorio completa quello della sezione precedente.
L'inversione circolare si può anche definire a partire dalla polarità. La cosa si può far scoprire in una attività laboratoriale, chiedendo di costruire con GeoGebra l'inversione circolare utilizzando lo strumento "polare". La costruzione è la seguente: dato un cerchio c di centro A, l'inverso circolare di un punto P coincide con l'intersezione tra la polare di P e la congiungente AP. Questa descrizione suggerisce una possibile doppia generalizzazione, che anch'essa può essere utilmente suggerita ed esplorata attraverso un approccio laboratoriale: intersecare la polare con la congiungente AP con A diverso dal centro, e considerare una conica, invece del cerchio, rispetto a cui tracciare la polare. La trasformazione così ottenuta è *l'inversione di Bellavitis-Hirst* (Bellavitis 1838) (Hirst 1865) (Raspanti, Tesi di dottorato 2017) (Raspanti, Giusto Bellavitis e la sua geometria di derivazione s.d.) e la ricerca delle sue proprietà offre uno spunto per un laboratorio molto interessante.
Chiedendo di esplorare come si trasformano le curve (a partire dalle più semplici, cerchi e rette) applicando l'inversione di Bellavitis Hirst, riteniamo che non sia difficile guidare gli studenti a scoprire diverse proprietà e a formulare diverse congetture, tra cui le seguenti:
i) Trasforma le rette nelle coniche per il punto A e per i punti B e C in cui la polare di A interseca la conica.

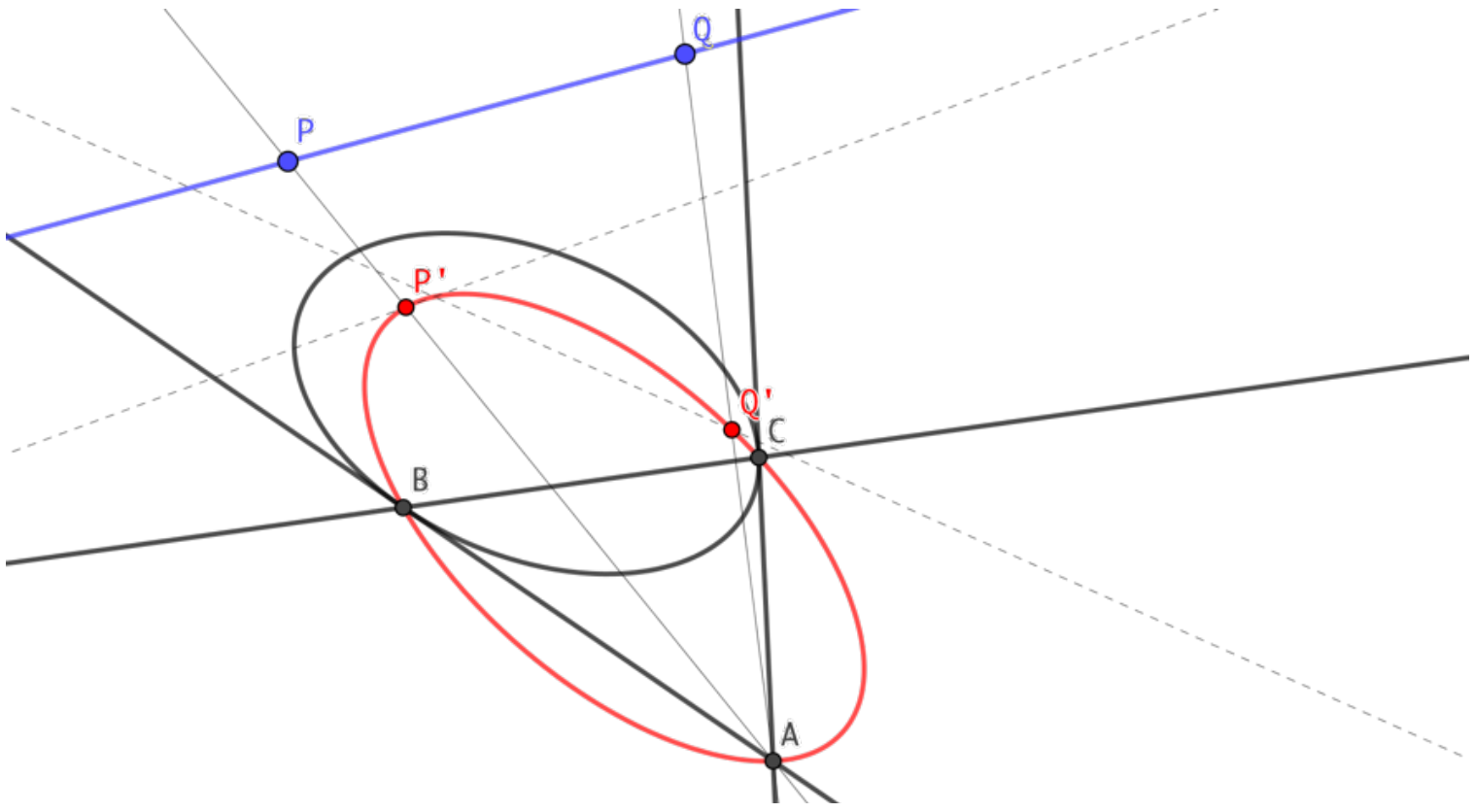

***Figura 28.*** *La retta per P e Q è trasformata nella conica per A, B, C, P' e Q'.*

ii) L'immagine di una retta per C è una conica riducibile che si spezza nella retta AC e in una retta per B. In particolare, l'immagine della retta CB si spezza nelle rette AB e AC.

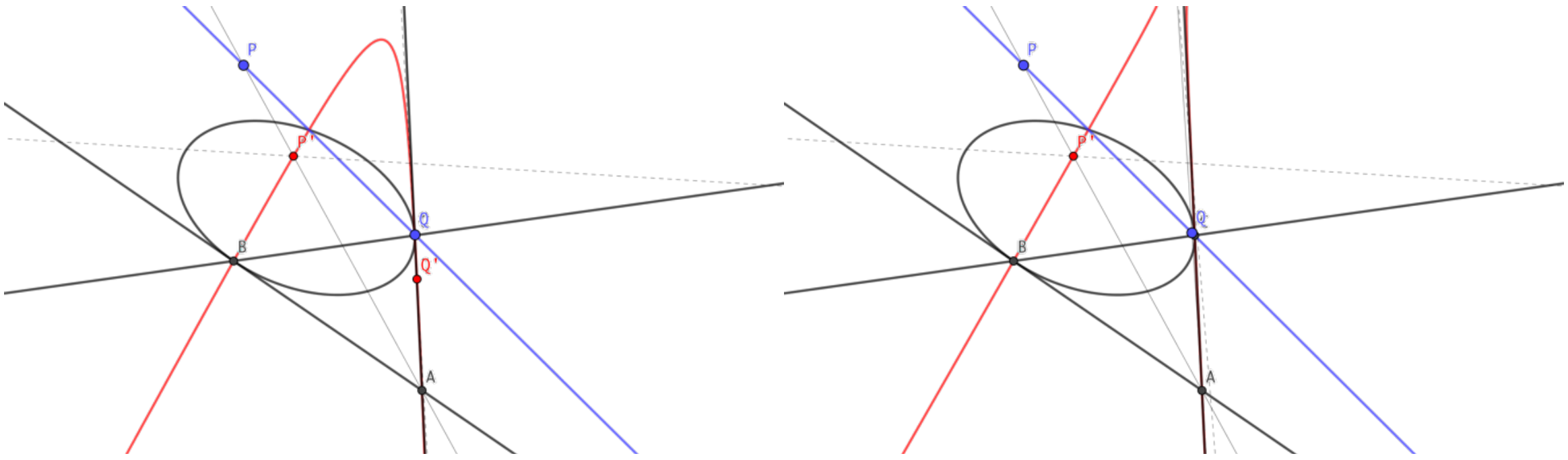

***Figura 29.*** *Avvicinando il punto Q al punto C si passa con continuità dalla figura a sinistra a quella a destra. Al "limite" l'immagine si spezza nella conica riducibile che passa per A, B P' e "doppiamente" per C, stando con ciò a significare che deve essere tangente alla conica base in C e quindi si spezza nelle rette BP' e AC.*

iii) L'immagine di una retta per A è una conica riducibile che si spezza nella retta BC e in una retta per A. In particolare, l'immagine della retta AC si spezza nelle rette AC e BC

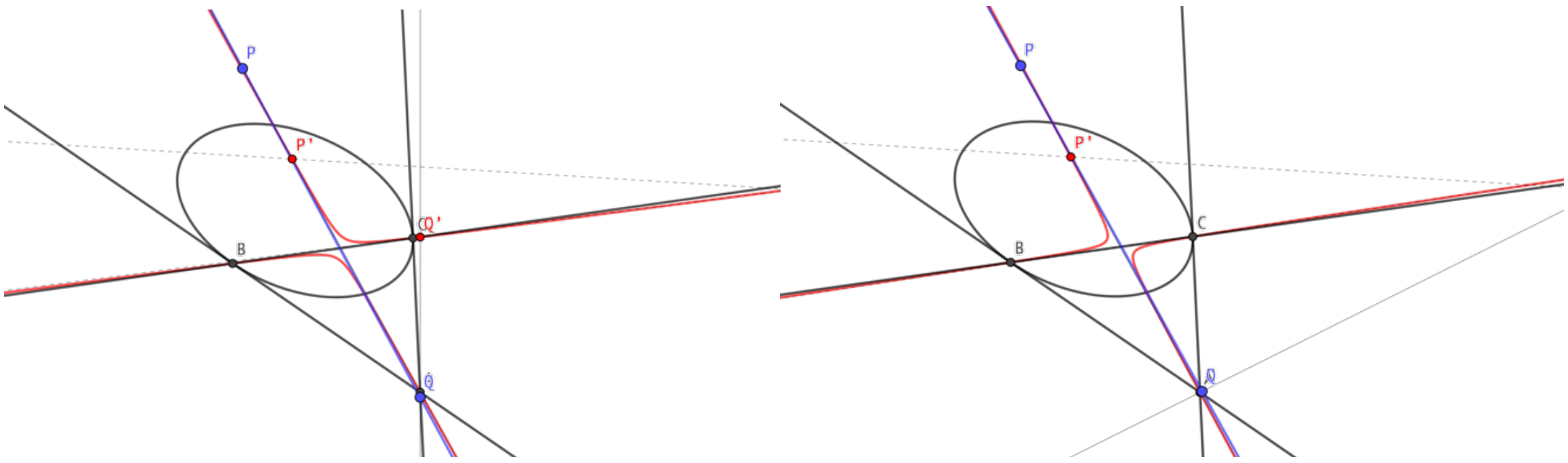

***Figura 30.*** *Avvicinando il punto Q al punto A si passa con continuità dalla figura a sinistra a quella a destra. Al "limite" l'immagine si spezza nella conica riducibile che passa per B, C, P' e "doppiamente" per A, stando con ciò a significare che deve spezzarsi nella retta AP' e nella retta BC. **

Queste immagini chiariscono come sia possibile con GeoGebra avere una visione dinamica delle proprietà di una figura. La possibilità di "vedere" le proprietà per deformazione incarna uno dei precetti didattici di Emma Castelnuovo: "Muovere le figure fa muovere le idee", che grazie a un software di Geometria Dinamica può essere sperimentato dai futuri insegnanti, che ne possono comprendere l'utilità didattica.

L'inversione di Bellavitis-Hirst si specializza nell'inversione circolare quando la conica è il cerchio di inversione, A è il centro del cerchio e C e B sono i punti di intersezione della polare del centro (retta all'infinito) con la circonferenza. Questi punti sono immaginari e prendono il nome di *punti ciclici.* Le immagini delle rette sono quindi coniche per i punti ciclici, quindi cerchi per il centro di inversione, come abbiamo già notato.[31]

L'inversione di Bellavitis-Hirst trasforma: il fascio di rette per C nel fascio di rette per B, in modo da trasformare la retta CB nella retta AB; il fascio di rette per A in sé, trasformando la retta AC in sé e la retta AB in sé.

È naturale a questo punto interrogarsi sulla natura delle trasformazioni indotte da proiettività tra fasci di rette che godono delle proprietà appena messe in evidenza. Il risultato fondamentale da tenere in mente per seguire il percorso è che se esiste una proiettività $\sigma$ tra il fascio di rette per un punto $P$ (che indichiamo con $P^*$) e il fascio di rette per un punto $Q$ (che indichiamo con $Q^*$) che fissa la retta $PQ$, allora intersecando ogni retta $r \in P^*$ con la sua corrispondente $\sigma(r) \in Q^*$otteniamo una conica passante per $P$ e $Q$ rispettivamente.[32] Questo ci permette di mostrare che la mappa di Bellavitis-Hirst può essere descritta anche così. Siano date tre proiettività $\sigma: B^* \rightarrow C^*, \tau: C^* \rightarrow B^*$ $\phi: A^* \rightarrow A^*$ tali $\sigma(BC) = CB$ e $\sigma(BA) = CA$; $\tau(CB) = BC$ e $\tau(CA) = BA$; $\phi(AB) = AB$ e $\phi(AC) = AC$. Per ogni $P$ che non appartiene al triangolo fondamentale, siano $r_B(P), r_C(P), r_A(P)$ le rette dei fasci per B, C, e A rispettivamente che passano per il punto P. Le immagini $\sigma(r_B(P))$ $\tau(r_C(P))$ $\phi(r_A(P))$ con le proiettività di cui sopra, si intersecano in un punto $P'$ che è l'immagine della mappa di Bellavitis-Hirst.[33]

Le inversioni di Bellavitis-Hirst sono trasformazioni da $P^2 \rightarrow P^2$ mentre le inversioni sono trasformazioni definite da $R^2 \cup \{\infty\} \rightarrow R^2 \cup \{\infty\}$. Riteniamo istruttivo discutere come è possibile interpretare le inversioni come trasformazioni quadratiche associate a un triangolo fondamentale e come l'inversione circolare può essere vista come (la restrizione di) una trasformazione quadratica ma anche come la composizione di una proiettività della retta complessa in sé con la mappa di coniugio.
Il percorso che parte dalla proiezione e arriva alle trasformazioni quadratiche illustra la nascita della geometria birazionale secondo un'angolatura interessante per il futuro insegnante che permette di esplorare i collegamenti con la geometria non euclidea (le inversioni circolari generano le isometrie del modello di Riemann-Beltrami-Poincaré (Arcozzi 2012)) e con il programma di Eralngen di Klein.

---

[31] È utile usare GeoGebra per vedere, oltre alle figure in movimento, anche le figure immaginarie.

[32] La costruzione è attribuita generalmente a Steiner ma già Newton la conosceva con il nome di *Generazione organica* (Lemma XXI dei Principia) è tra le più gradite nei laboratori dedicati alle costruzioni proiettive con GeoGebra nell'ambito del Corso Monografico di Matematica.

[33] Queste cose sono presentate in modo chiaro in (Beltrami 1862). Vedi anche, per un quadro storico, (Vaccaro 2020).

## 5.1 Esplorando con Cremona e Hirst

Una delle questioni naturali che possono essere sollevate dai partecipanti di un laboratorio in cui si tratti dell'inversione di Bellavitis-Hirst è: qual è l'immagine di una circonferenza? Il problema non è semplice da trattare in generale ma offre l'occasione per un altro genere di sperimentazione, che ha a nostro avviso un grande valore formativo, quella cioè di partecipare al dialogo tra grandi matematici del passato potendo ripercorrere, grazie all'uso del software, l'iter di una scoperta, fatta grazie a una straordinaria intuizione geometrica, cioè grazie alla capacità di vedere con gli occhi della mente le relazioni geometriche tra gli oggetti di studio, che un software di Geometria Dinamica permette di mostrare concretamente sullo schermo di un computer. Questo percorso è stato descritto in (Raspanti, Dall'inversione circolare all'inversione quadratica: aspetti storici e potenzialità didattiche 2016) ma pensiamo che valga la pena di presentarne una parte in conclusione del nostro lavoro.
Iniziamo con la lettura del seguente estratto di una lettera di Hirst a Cremona del 16 ottobre 1864.

> *Forse ricorderai la nostra conversazione (vicino Porta S. Stefano) sulla quartica di Steiner di terza classe, che ha tre cuspidi, e che tocca la retta all'infinito nei punti circolari;*[34] *poi ho notato che mi ero imbattuto sulla stessa curva, ma non riuscivo a ricordare come. Ecco un altro metodo per generarla da aggiungere ai tanti già proposti da Steiner, Schaefli, Schroeter e te stesso.*
> *È l'inviluppo delle iperboli, circoscritte ad un triangolo equilatero, che hanno i loro angoli asintotici uguali l'uno l'altro e ad ogni angolo del triangolo.*
> *Una definizione equivalente (mediante inversione quadrica) è la seguente.*
> *Preso un cerchio (B), ed un cerchio (C) passante per il centro del primo e avente raggio doppio. Sia A il centro esterno di similitudine dei due cerchi.*
> *Allora l'inversa quadrica di (B) relativamente all'origine A e al cerchio fondamentale (C) sarà esattamente la quartica di Steiner. In altre parole, se p e p' costituiscono una coppia di punti coniugati, relativamente a (C), su una retta passante per A, allora mentre p descrive il cerchio (B), il suo inverso (quadrico) p' descriverà la quartica di Steiner. Le tre cuspidi di quest'ultima saranno A, e i due punti $A_1$, $A_2$ di contatto delle due tangenti condotte da A a (C).*[35]*(…) Osservo che l'inverso quadrico del cerchio, concentrico con (B)*

[34] Questa curva, chiamata anche *deltoide* per la forma che ricorda la lettera maiuscola greca Delta, fu considerata da Eulero nel 1745 in relazione a un problema di ottica. Fu successivamente studiata da Steiner nel 1856 e per questo è anche nota come quartica di Steiner. Può essere generata come traiettoria di un punto sul bordo di un cerchio di raggio b che rotola senza strisciare internamente a un cerchio di raggio a nell'ipotesi che a/b=3, da cui il nome di ipocicloide tricuspide. L'identità tra la quartica di Steiner e l'ipocicloide fu dimostrata solo alcuni anni più tardi, e per la prima volta per via puramente geometrica, da Luigi Cremona in (Cremona 1865). Per informazioni ulteriori su questa curva, rimandiamo a (Weisstein s.d.).

[35] (Nurzia 1999), pp. 47-50: «You perhaps remember our conversation (near the Porta S. Stefano) on Steiner's quartic of the third class, which has three cusps, and touches the line at infinity at the circular points; I then remarked that I myself had stumbled on the same curve but could not recall how. Here is another method of generating it to be added to the many already proposed by Steiner, Schaefli, Schroeter and yourself. It is the envelope of the hyperbolas,

*e il cui diametro è AC, è la retta all'infinito e tutte le proprietà dei cerchi concentrici (AC) e (B) si possono immediatamente trasferire per inversione quadrica alla quartica di Steiner.*

Alla lettera è allegato un disegno (**Figura 31**, a sinistra). La prima richiesta da cui prende le mosse il nostro laboratorio è quella di ricostruire il disegno con GeoGebra (**Figura 31**, a destra, dove è stata aggiunta la polare di *p* per visualizzare meglio l'inversione di Bellavitis-Hirst rispetto alla circonferenza centrata in C e al polo A).

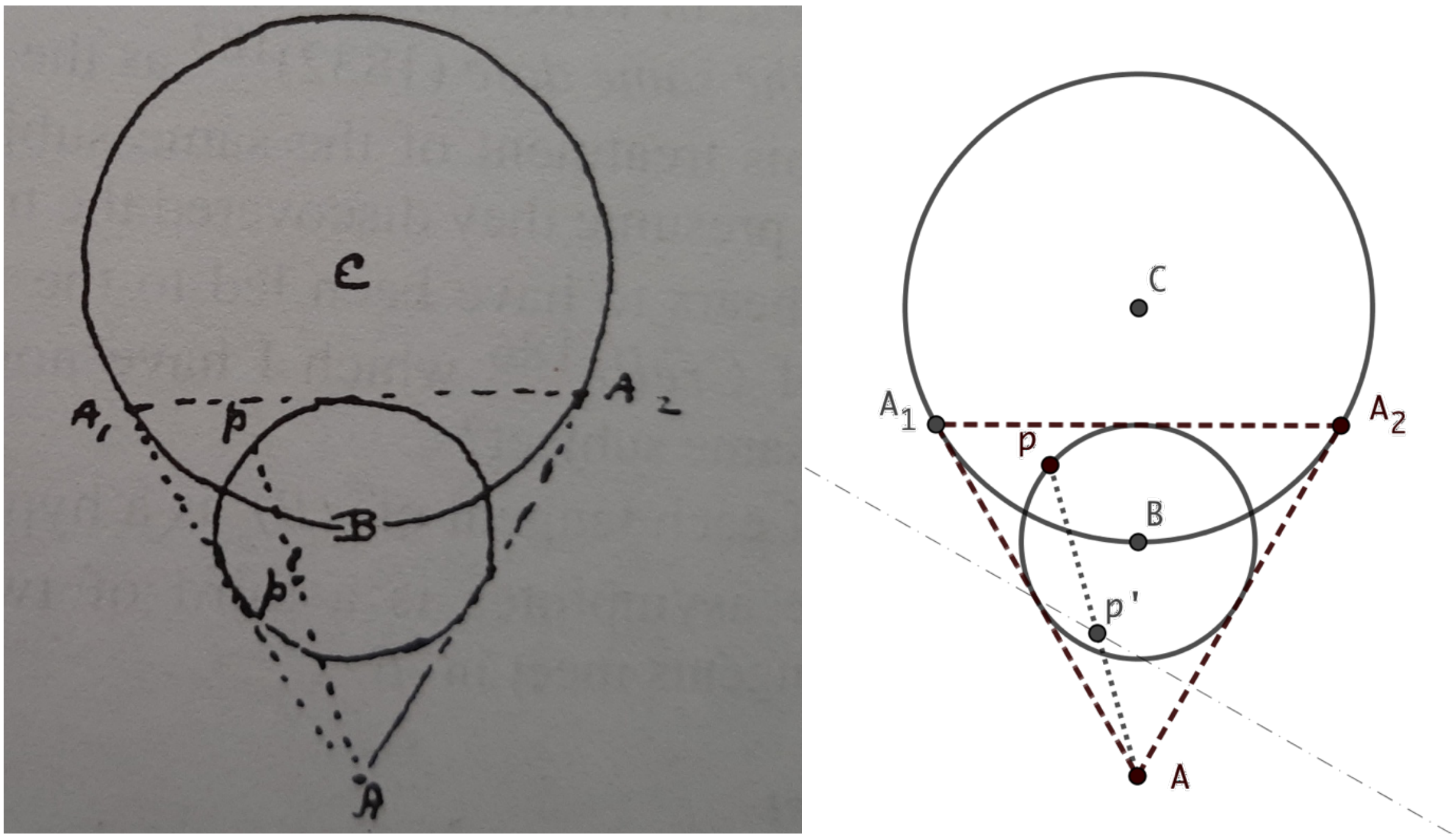

***Figura 31.*** *A sinistra il disegno di Hirst, a destra la sua ricostruzione con GeoGebra. La polare di p rispetto alla circonferenza centrata in C è rappresentata con una linea tratteggiata.*

Una volta costruita l'immagine di un punto *p*, l'immagine del cerchio centrato in B si ottiene con lo strumento "Luogo".

---

circumscribed to an equilateral triangle, which have their asymptotic angles equal to one another, and to each angle of the triangle. An equivalent definition (by Quadric Inversion) is the following. Take any circle (B), and a circle (C) passing through the centre of the first and having twice as great a radius. Let A be the external centre of similitude of the two circles. Then the quadric inverse of (B) relative to the origin A, and the fundamental circle (C) will be precisely Steiner's quartic. In other words, if p and p' be a pair of conjugate points, relative to (C), on a line passing through A then as p describes the circle (B), its quadric inverse p' will describe Steiner's quartic. The three cusps of the latter will be at A, and the points $A_1$, $A_2$ of contact of the two tangents from A to (C). (…) I may observe that the quadric inverse of the circle, concentric with B and whose diameter is AC is the line at infinity, and all properties of the concentric circles (AC) and B can at once be transferred by quadric inversion to Steiner's quartic. »

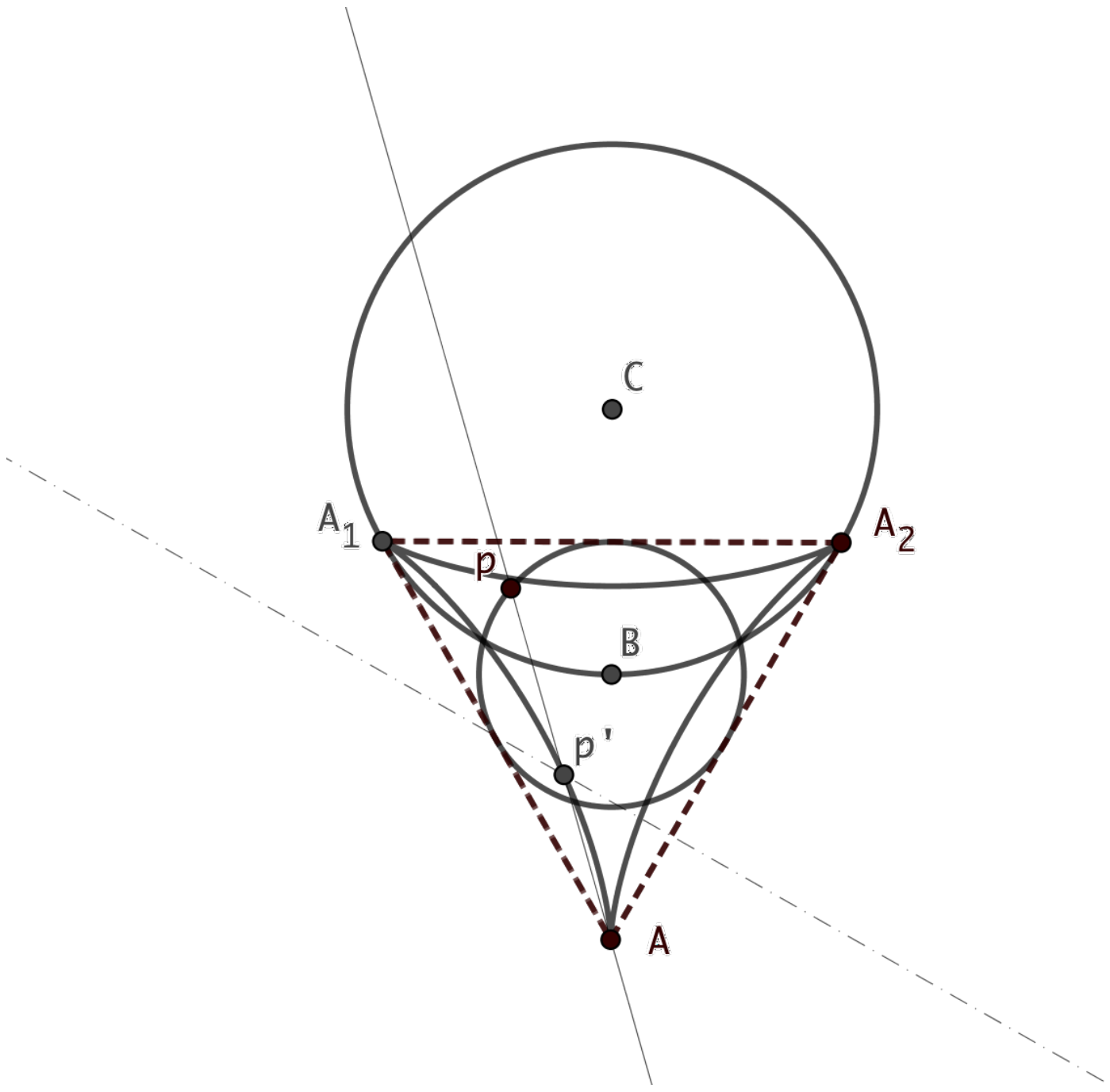

*Figura 32. La forma del luogo è quella dell'ipocicloide di Steiner. **

Cominciamo l'esplorazione delle proprietà del luogo con lo scopo di dimostrare la coincidenza con l'ipocicloide di Steiner. Nella lettera di Hirst viene ricordato che la curva di Steiner è una quartica con tre cuspidi; come far discendere questo dalle proprietà dell'inversione di Bellavitis-Hirst? Per esplorare la questione, cominciamo a studiare l'iniettività e la suriettività della trasformazione e successivamente come si trasformano le rette e le coniche.
Nella **Figura 33**, p' è l'immagine di p e quindi p, p' e A sono allineati (retta rossa). La retta blu è la polare di p rispetto alla circonferenza (C) mentre la retta verde è la polare di p' rispetto alla circonferenza (C). L'involutorietà della trasformazione segue da una proprietà involutiva della polarità, secondo cui se P e p sono una coppia polo e polare e Q e q son un'altra coppia polo e polare, allora Q appartiene a p se e solo se P appartiene a q.
Usando la possibilità offerta dal software di muovere il punto p trascinandosi dietro l'intera costruzione, si vede facilmente che l'applicazione contrae la retta $AA_1$ in $A_1$ e la retta $AA_2$ in $A_2$, indicando anche chiaramente la linea della dimostrazione.[36] Si vede anche che la retta $A_1\,A_2$ viene invece contratta sul punto $A$, e si vede anche che la ragione sta nel fatto che la polare di A è proprio la retta $A_1A_2$.

[36] Invitiamo a sperimentare con il file Figura34.ggb, disponibile all'indirizzo web http://matstor.wikidot.com/materialiggb.

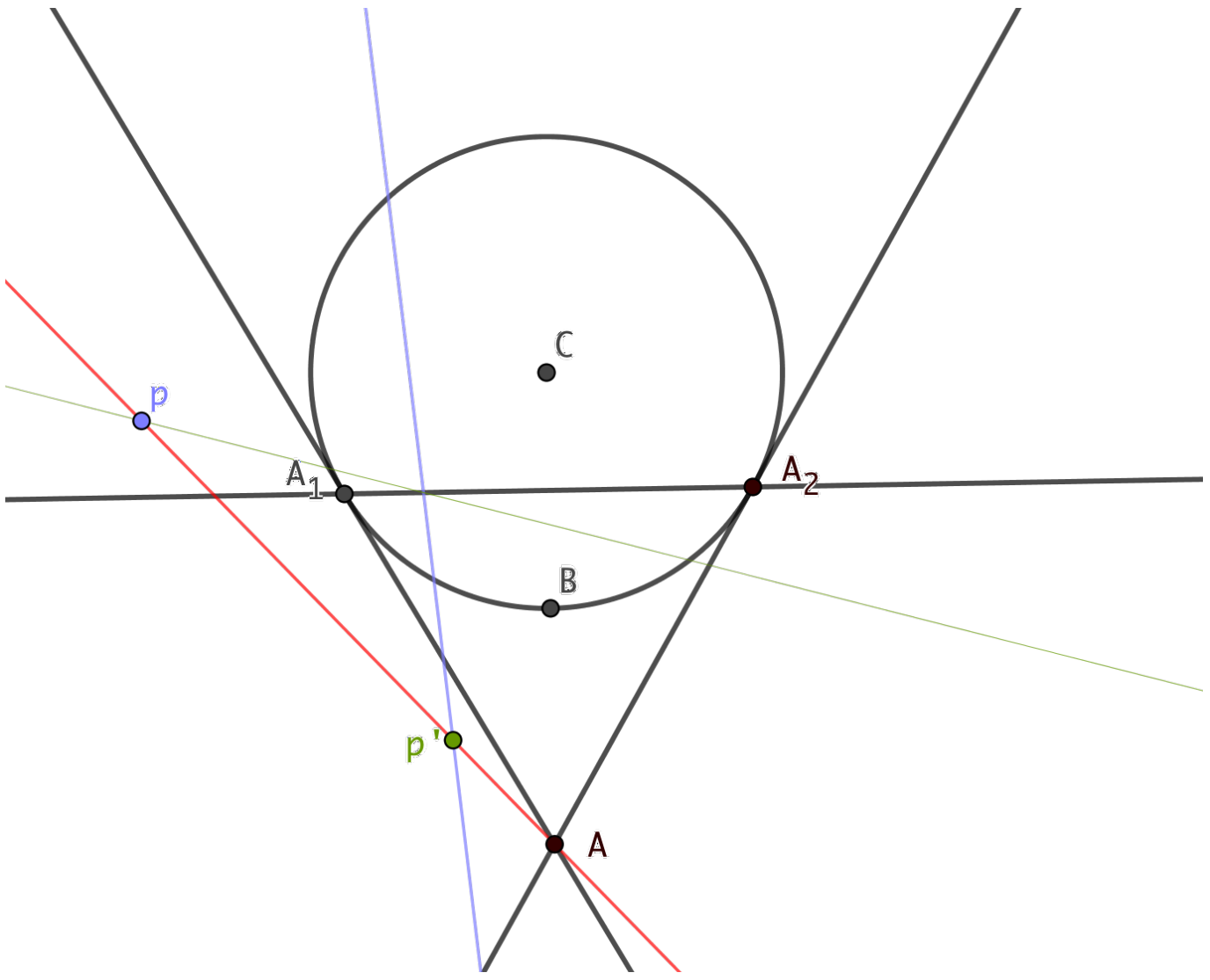

***Figura 33.*** *L'inversione di Bellavitis-Hirst è involutiva fuori dal triangolo fondamentale.*

Con GeoGebra è possibile anche illustrare come il punto $A_1$ ($A_2$) viene *scoppiato* sull'intera retta $AA_1$ ($AA_2$), mentre il punto A viene scoppiato sulla retta $A_1A_2$. Per far ciò, come è illustrato nella figura seguente, si vincoli p a muoversi su una retta per $A_1$ (o per un altro dei punti base). Al tendere di p ad $A_1$ su questa retta, il punto p' tende ad un particolare punto della retta $AA_1$, che dipende dalla retta per $A_1$ che abbiamo scelto.

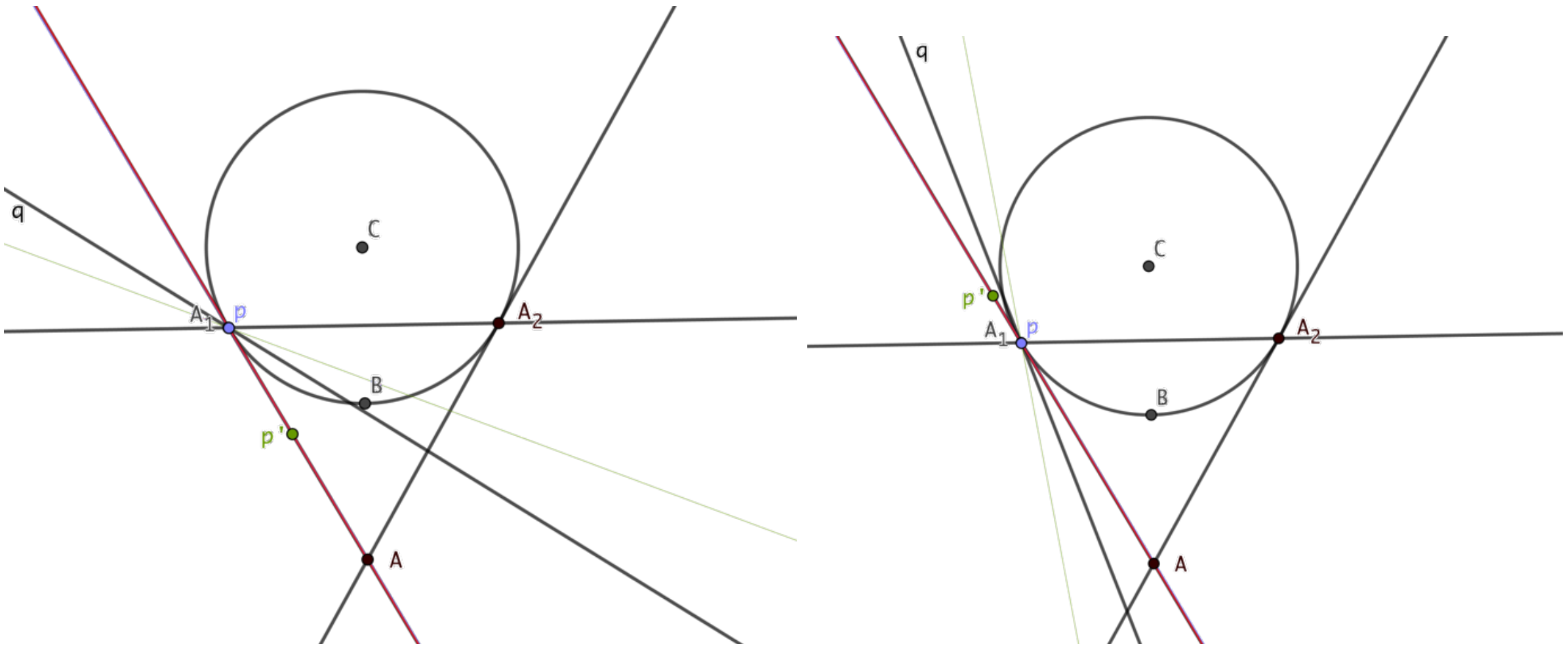

***Figura 34.*** *Quando p si avvicina ad $A_1$ lungo una retta, il limite di p' dipende dalla retta lungo cui si muove p. La retta $AA_1$ rappresenta quindi l'intorno del primo ordine di $A_1$. **

Chiediamoci ora qual è l'immagine di una retta. Abbiamo già visto che i vertici del triangolo fondamentale sono scoppiati, cominciamo quindi con il considerare una retta che non passa per i vertici. Trasformando un punto di una retta (rappresentata in azzurro) e costruendo il luogo dell'immagine di tale punto con GeoGebra (rappresentato in rosso), otteniamo figure come le seguenti.

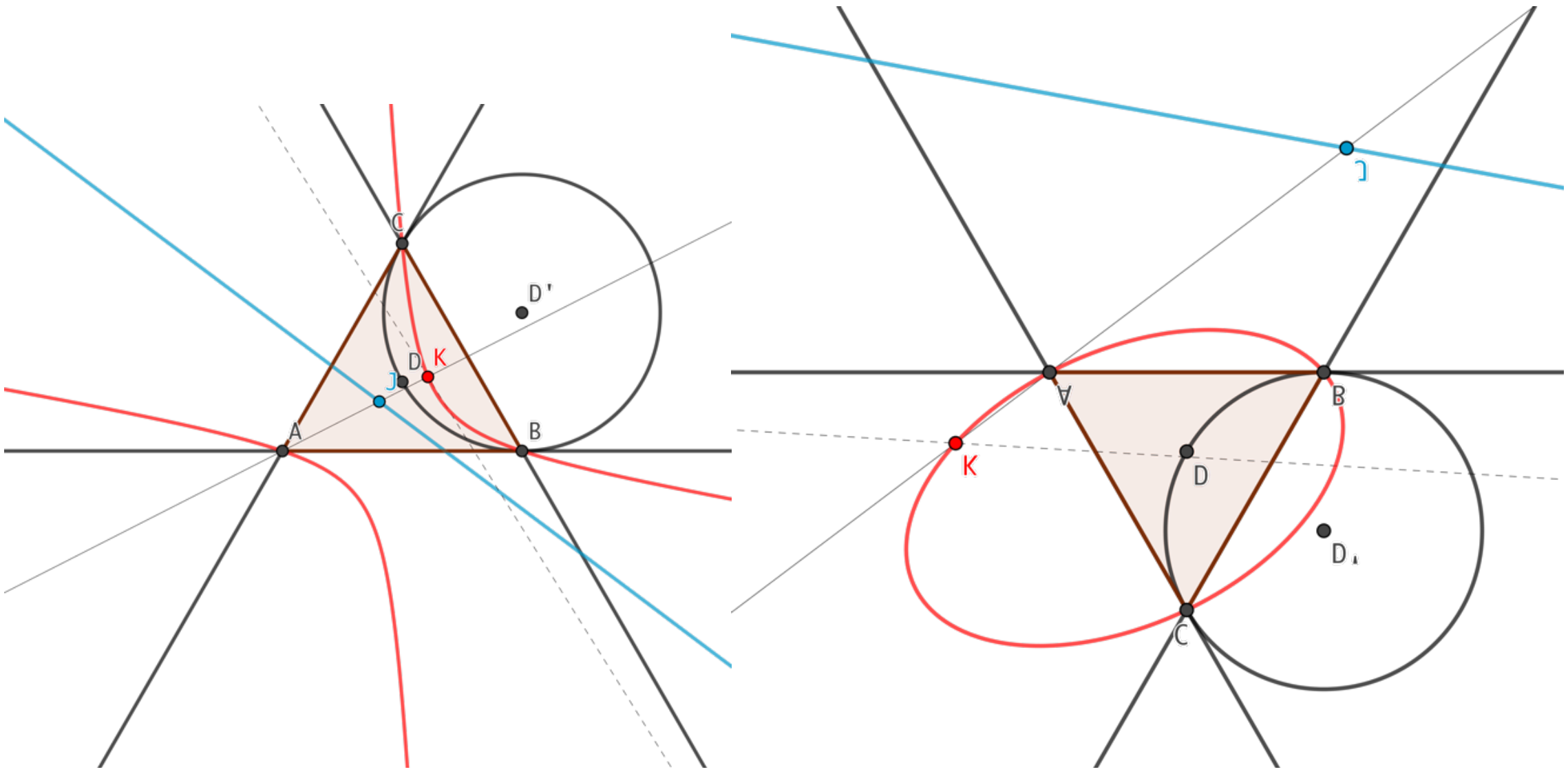

***Figura 35.*** *L'immagine di una retta che non contiene vertici è una conica passante per i vertici. A sinistra, un'iperbole, a destra un'ellisse.*

La curva immagine passa per i vertici. Infatti, come abbiamo già osservato, l'intersezione della retta originale con la retta $AA_1$ viene contratto in $A_1$, ecc.  La curva immagine sembra essere un'iperbole. Cambiando la retta, l'immagine continua a sembrare una conica. Possiamo ipotizzare quindi che il sistema lineare delle rette viene trasformato per inversione nel sistema lineare delle coniche per i tre vertici del triangolo fondamentale.
Cosa succede se la retta passa per un vertice, per esempio per il vertice A della Figura 1. Abbiamo già detto che A viene scoppiato nella retta per gli altri due vertici (CB in figura). Togliendo questa componente resta la cosiddetta *trasformata stretta* che, come vediamo nella figura, e come si può dimostrare geometricamente e facilmente, è la retta stessa. Quindi, le rette per A, che non contengono altri vertici, si trasformano in una conica riducibile per i vertici di cui una componente è la trasformata stretta e l'altra è lo scoppiamento del punto A.

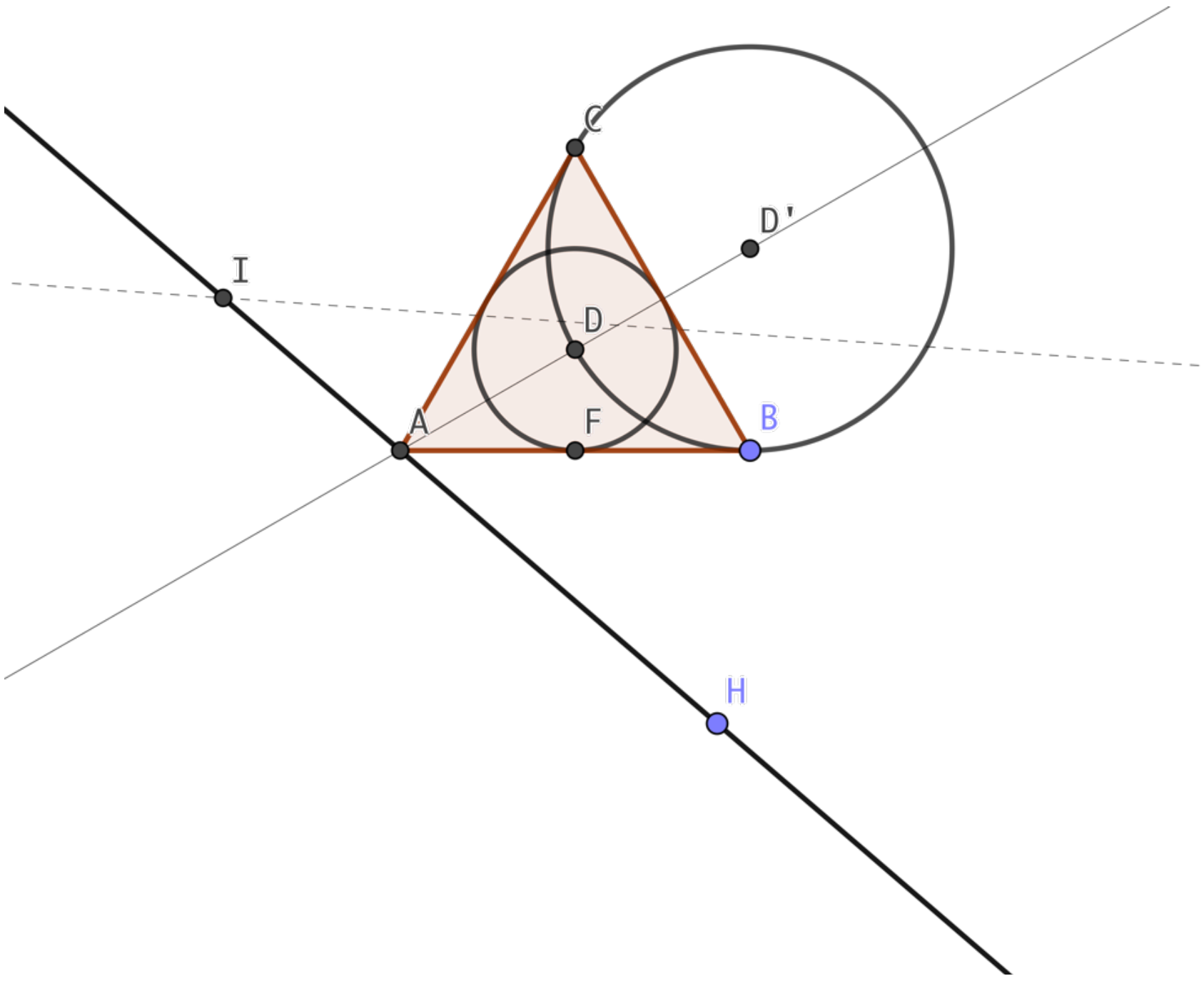

***Figura 36.*** *Una retta per A viene trasformata in sé stessa (escluse le rette AC e BC)*

Si può fare la stessa analisi per una retta passante per C (con riferimento alla figura 2).

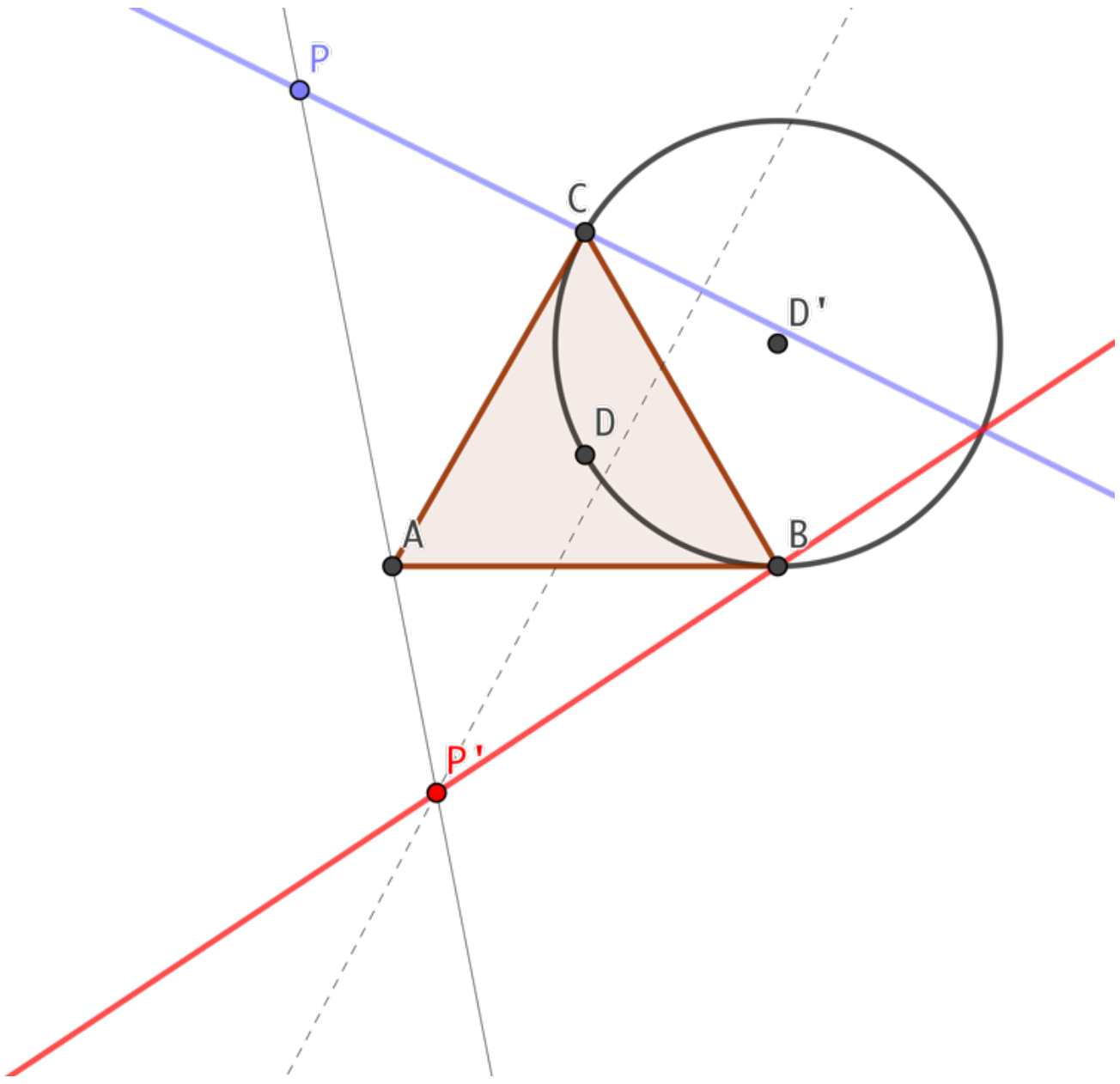

*__Figura 37.__ Una retta per C (escluse le rette CB e CA) viene trasformata in una conica riducibile che ha una componente (trasformata stretta) che è una retta per B (che incontra la prima in un punto della circonferenza fondamentale) e una seconda componente, che viene dallo scoppiamento di C, che coincide con la retta CA.*

Vediamo ora cosa succede di un cerchio che interseca le tre rette fondamentali in due punti che appartengono ad una sola di tali rette.

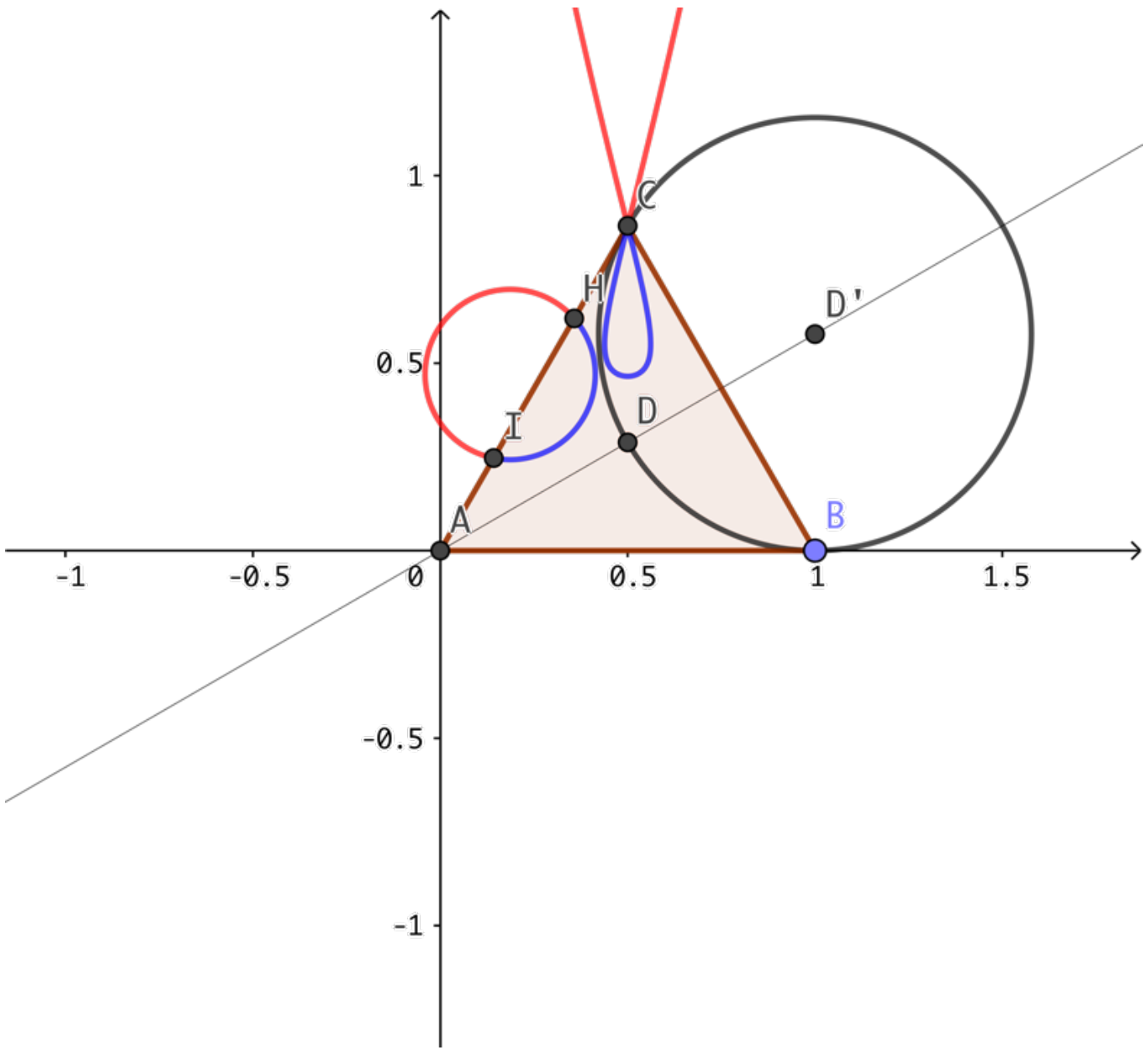

*__Figura 38.__ La figura illustra la seguente proprietà fondamentale: La retta AC viene contratta nel punto C e una curva che interseca AC in h punti semplici viene trasformata in una curva con un punto h-plo in C. **

Non è difficile immaginare cosa succede quando una curva che non passa per i vertici interseca la retta CB in h punti semplici

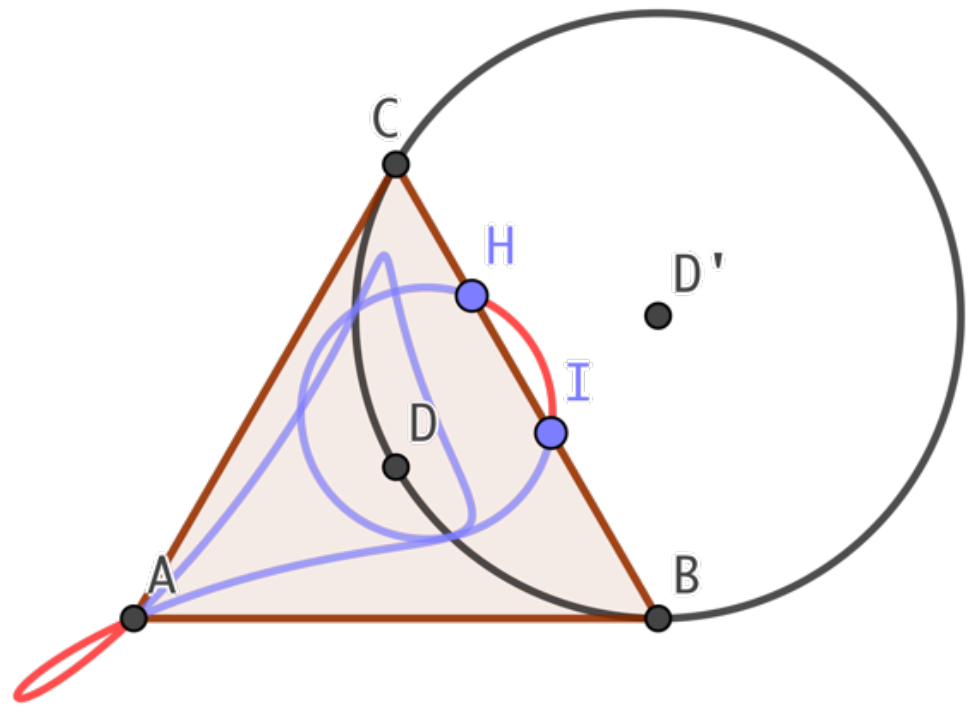

*__Figura 39.__ La retta BC viene contratta nel punto A e una curva che interseca BC in h punti semplici viene trasformata in una curva con un punto h-plo in A. **

Chiediamoci ora in quale curva viene trasformata una conica passante (semplicemente) per i tre vertici. Per l'involutorietà della trasformazione e per quanto abbiamo congetturato considerando l'immagine di una retta, ci aspettiamo che vengano trasformate in rette, cosa che viene confermata da GeoGebra.[37]

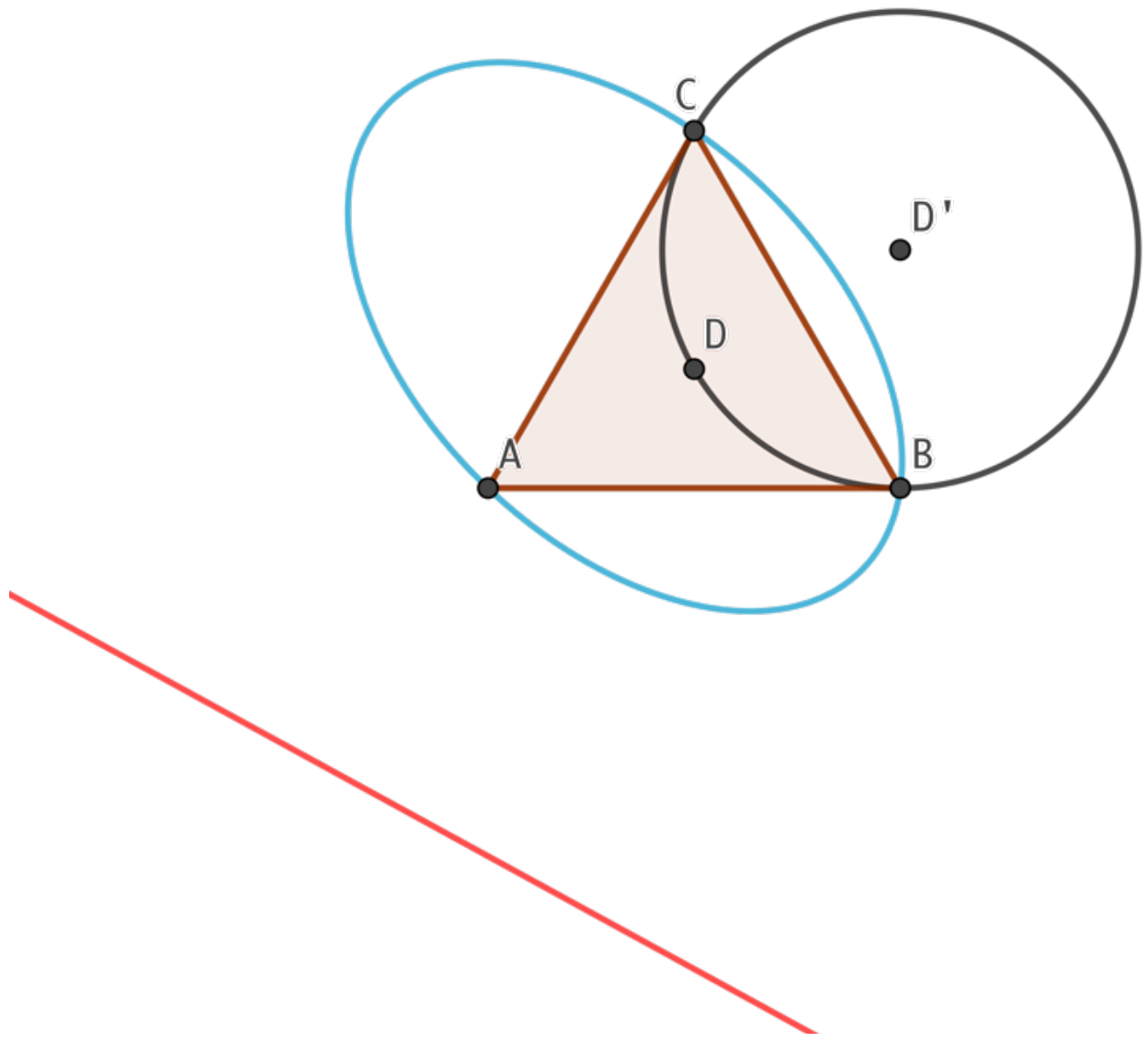

*__Figura 40.__ Una conica (azzurra) passante per ABC viene trasformata in una curva riducibile. La trasformata stretta, ottenute eliminando gli scoppiamenti di A, B e C è una retta (rossa). **

Ci chiediamo: se C è una curva di grado n, cosa possiamo dire del grado della trasformata[38] C'? Supponiamo che C non passi per i vertici A, B e C. Il grado di C' è uguale al numero di intersezioni con una retta generica. Essendo la trasformazione di Bellavitis-Hirst biiettiva fuori dal triangolo fondamentale, il numero delle intersezioni di C' con una retta generica è uguale al numero di intersezioni di C con la controimmagine della retta, che è una conica per i vertici, quindi, per il teorema di Bezout, è uguale a 2n. Se però la curva passa con molteplicità pari a tA, tB e tC per i vertici, alla trasformata stretta vanno tolti gli scoppiamenti dei vertici (rette) contati con la loro molteplicità. Quindi il grado è 2n-tA-tB-tC.

[37] Per una dimostrazione, si può far riferimento alla forma analitica della trasformazione, data in fondo al lavoro.
[38] Stiamo assumendo che la trasformata di una curva algebrica sia ancora una curva algebrica. Questo si può dimostrare molto semplicemente osservando che le equazioni della trasformazione sono algebriche (si tratta di intersecare rette che dipendono algebricamente dalle coordinate del punto che stiamo trasformando.

Veniamo finalmente alla trasformazione della circonferenza tangente internamente al triangolo equilatero di vertici A, B e C. Poiché tale curva (di grado 2) non contiene alcun vertice, il grado della trasformazione è uguale a 4. Poiché la circonferenza è tangente ad ognuna delle tre rette del triangolo fondamentale, la curva immagine avrà una cuspide in ognuno dei vertici, per le ragioni illustrate nella figura seguente.

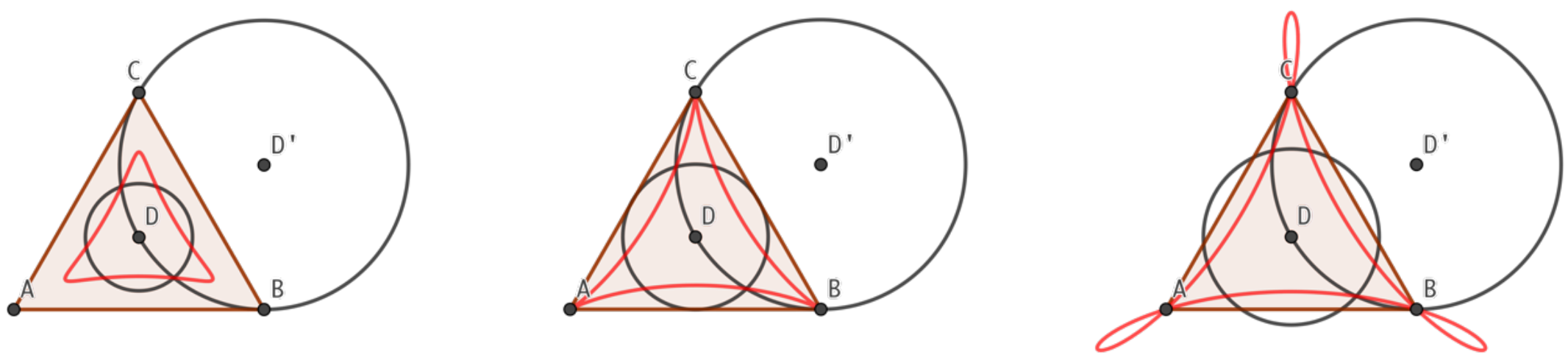

***Figura 41.*** *Nella prima figura, l'immagine del cerchio nero di centro D è la quartica rossa priva di punti singolari, in quanto il cerchio non interseca i lati del triangolo. Nella terza figura, l'immagine ha un nodo in ogni vertice in quanto interseca ogni lato in due punti. Nella figura centrale, la circonferenza si trova in una posizione limite tra le due precedenti. In ogni vertice, l'immagine ha quindi una cuspide, corrispondente alla tangenza della circonferenza al lato corrispondente. **

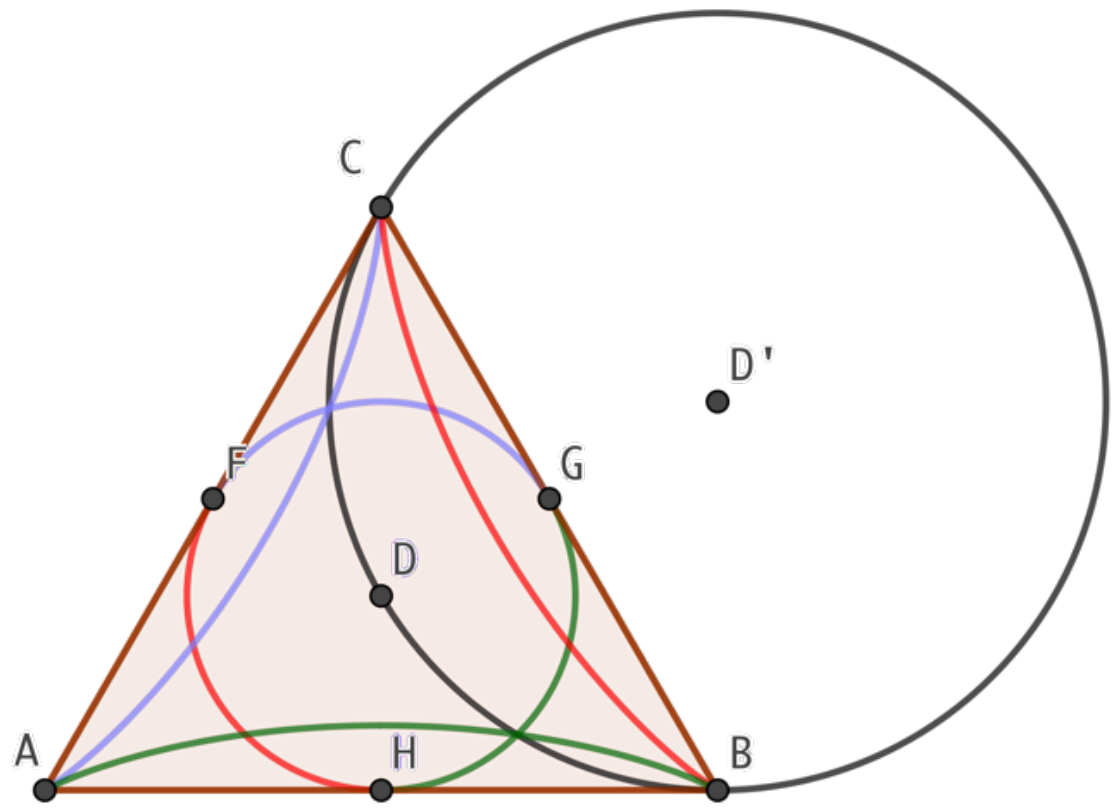

***Figura 42.*** *Il cerchio tangente internamente al triangolo ABC viene trasformato un una quartica con tre punti doppi cuspidali, corrispondenti ai tre punti di tangenza. **

E veniamo alla caratterizzazione data da Hirst della ipocicloide di Steiner come inviluppo. La circonferenza tangente internamente al triangolo si può definire come l'inviluppo delle sue rette tangenti. Allora l'immagine con l'inversione di Bellavitis-Hirst si può definire come l'inviluppo delle curve immagine delle rette tangenti. Esse sono coniche per i tre vertici, che sono iperboli con l'ulteriore proprietà che i loro asintoti formano un angolo uguale a $\pi/3$.[39]

[39] È possibile anche convincersi, con GeoGebra, che l'angolo tra gli asintoti delle iperboli tangenti è costante ed è uguale a $\pi/3$. Basta intersecare un cerchio di raggio abbastanza grande e centrato nel centro dell'iperbole, intersecare il cerchio con l'iperbole e calcolare l'angolo tra le diagonali. Il limite al crescere del raggio di tale angolo è il valore cercato. Per una costruzione, con riga e compasso degli asintoti di un'iperbole assegnata cfr. (Ignoto s.d.)

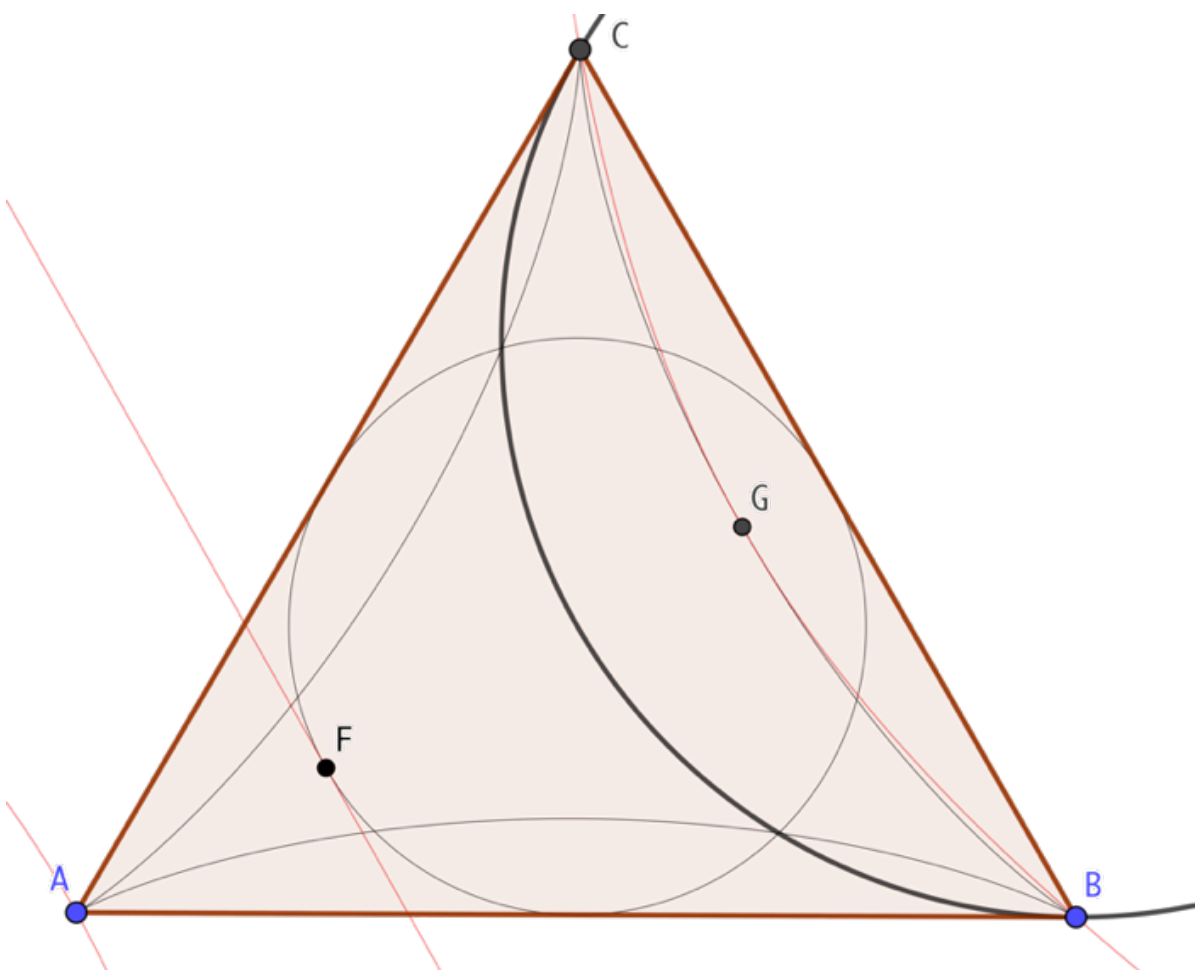

***Figura 43.*** *Sia F un punto della circonferenza $\mathcal{C}$ inscritta nel triangolo ABC e sia G la sua immagine con l'inversione di Bellavitis-Hirst. G appartiene alla quartica $\mathcal{Q}$ immagine di $\mathcal{C}$. La retta tangente in F a $\mathcal{C}$ viene trasformata nella conica tangente a $\mathcal{Q}$ in G passante per i punti A, B e C. **

Poiché il cerchio circoscritto al triangolo ha lo stesso centro del cerchio tangente internamente, i due cerchi sono tangenti nei punti ciclici. Allora, l'ipocicloide di Steiner è tangente all'immagine del cerchio circoscritto. Essa è una retta. Non è difficile usare GeoGebra per identificare questa retta. Infatti, al variare del punto D sulla circonferenza circoscritta (cfr. **Figura 44**), la congiungente ZD e la polare di D rispetto al cerchio di inversione restano parallele e quindi l'immagine del cerchio circoscritto è la retta all'infinito.

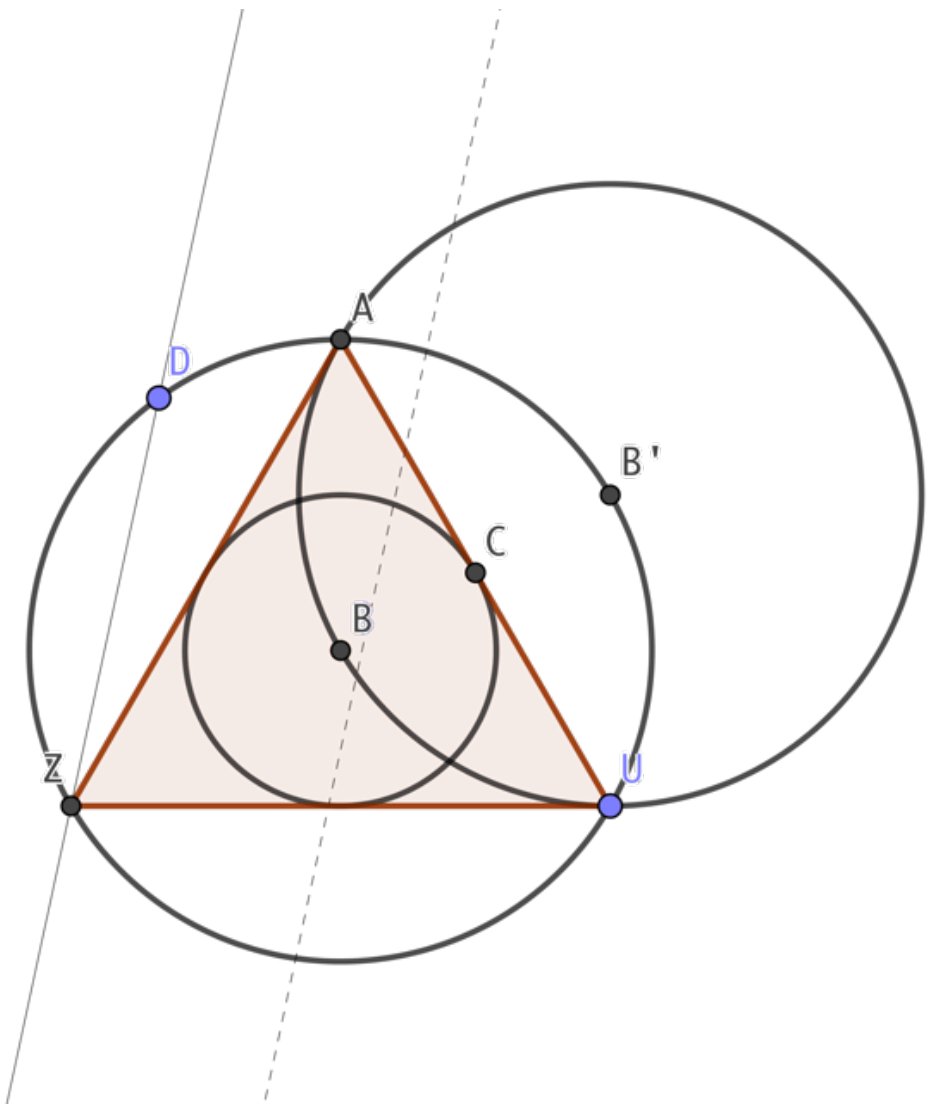

***Figura 44.*** *Per ogni punto D sulla circonferenza per Z, A e U la congiungente ZD è parallela alla polare di D rispetto alla circonferenza fondamentale.*

Ci sembrano opportune, anche per collegare questo argomento con la geometria algebrica moderna, alcune osservazioni sulla forma analitica della trasformazione. In fondo, va ribadito che l'obiettivo non è quello di eliminare le formule dalla matematica, ma di collegarle strettamente alla visualizzazione geometrica. Se scegliamo un sistema di coordinate cartesiane ortogonali monometriche in cui A abbia coordinate (0,0) e B abbia coordinate (1,0), allora $C(\frac{1}{2},\frac{\sqrt{3}}{2})$. Il cerchio di inversione ha quindi equazione

$$1-2x+x^2-\frac{2}{\sqrt{3}}y+y^2=0$$

La polare del punto $(\xi,\eta)$ è la retta

$$1-\frac{\eta}{\sqrt{3}}+x(\xi-1)-\xi+\left(\eta-\frac{1}{\sqrt{3}}\right)y=0$$

Intersecando la polare con la retta congiungente il punto A con il punto $(\xi,\eta)$ otteniamo la forma analitica della trasformazione

$$\left(-\frac{3x-3x^2-\sqrt{3}xy}{-3x+3x^2-\sqrt{3}y+3y^2},\frac{y(-3+3x+\sqrt{3}y)}{-3x+3x^2-\sqrt{3}y+3y^2}\right)$$

I due numeratori e il denominatore comune definiscono tre coniche che passano per i tre vertici del triangolo fondamentale.

Si poteva partire da qui per definire la trasformazione, che appare in questa forma come un caso particolare di trasformazione quadratica. Questo approccio si ricollega alla maniera standard di introdurre le trasformazioni quadratiche come caso particolare delle trasformazioni birazionali, che furono introdotte per la prima volta, in tutta la loro generalità, dallo stesso Cremona in due famosi lavori che stanno alla base della geometria algebrica.

# 6. Conclusioni

In conclusione, vogliamo ritornare sull'obiettivo che avevamo indicato nell'introduzione, cioè quello di illustrare come un software di Geometria Dinamica possa aiutare a veicolare l'insegnamento di contenuti superiori secondo modalità appropriate al completamento della formazione matematica del futuro insegnante. Prendiamo per esempio il caso delle trasformazioni birazionali, che abbiamo trattato nella sezione 5. Si tratta di uno dei concetti fondamentali della geometria algebrica il cui oggetto, secondo le vedute del programma di Erlangen, è proprio lo studio delle proprietà invarianti per queste trasformazioni. L'eventuale conoscenza di questi argomenti per il futuro insegnante proviene dai corsi di Geometria superiore, che non sono pensati per le loro specifiche esigenze ma per quelle dei futuri ricercatori. Non viene ad esempio messo sufficientemente in luce, a nostro parere, il rapporto delle più semplici tra queste trasformazioni (quelle quadratiche) con la geometria elementare, con l'inversione e con le geometrie non euclidee perché, dal punto di vista della formazione dei ricercatori, non è opportuno soffermarsi troppo su queste questioni che non sono immediatamente legate alla ricerca. Si tratta invece di collegamenti estremamente fecondi per il futuro insegnante, sia per collocare i contenuti del suo futuro insegnamento in un quadro più ampio, sia come stimolo per organizzare attività complementari, come l'esplorazione dell'inversione circolare e della trasformazione di Bellavitis-Hirst. Quando attività laboratoriali di questo genere vengono proposte a scuola, si richiede all'insegnante una conoscenza più ampia di quella che viene trasmessa allo studente ma diversa da quella trasmessa nei corsi istituzionali di Geometria superiore; crediamo che siano necessarie piuttosto attività laboratoriali del tipo di quelle descritte nelle sezioni 3 e 4. La conoscenza dei collegamenti di un argomento con le altre parti della matematica non è centrale infatti nella formazione del ricercatore, che ha tutto il tempo di apprezzarla e coltivarla nell'intero arco della sua vita professionale, mentre lo è per il futuro insegnante, che proprio nei collegamenti può trovare armi efficaci per motivare gli studenti allo studio della matematica e per chiarire i punti oscuri del proprio insegnamento. I collegamenti sono anche importanti per riuscire a impostare in maniera efficace le attività laboratoriali. Questo approccio ai contenuti matematici superiori che privilegia i collegamenti rispetto agli approfondimenti, che abbiamo definito "estensivo" trova nel software di geometria dinamica, almeno per quanto riguarda i contenuti geometrici, un alleato ideale in quanto permette, con riferimento alla dicotomia analisi/sintesi che abbiamo ricordato a p. 14Il tesoro dell'Analisi, di *analizzare* i contenuti e i collegamenti in maniera efficace senza la necessità di s*intetizzare* la teoria in maniera completa, ma aiutando a scoprire  gli snodi cruciali ai quali è possibile limitare la *sintesi*.

## Bibliografia

AAVV. s.d. *GeoGebra.* https://www.geogebra.org/.

Arcozzi, Nicola. 2012. «Beltrami's models of non-eucidean geometry.» In *Mathematics in Bologna 1861-1960*, di Salvatore Coen, 1-30. Basel: Springer.

Arzarello, F. Bartolini Bussi, M.G., Leung, A.Y.L., Mariotti, M.A, e I. Stevenson. 2012. «Experimental approaches to theoretical thinking: Artefacts and proof.» In *Proof and Proving in Mathematics Education*, di G. Hanna e M. de Villier, 1-10. New York: Springer.

Bartolini Bussi, M. G., e Michela Maschietto. 2006. *Macchine matematiche: dalla storia alla scuola.* Milano: Springer.

Bellavitis, Giusto. 1838. «Saggio di Geometria derivata.» *Nuovi Saggi dell'i.r. Accademia di Scienze lettere ed Arti di Padova* 243-289.

Beltrami, Eugenio. 1862. «Intorno alle coniche dei nove punti ed ad alcune quistioni che ne dipendono.» *Memorie dell'Accademia delle Scienze dell'Istituto di Bologna* Ser. II, vol. 2: 361-395.

Capone, Roberto, Enrico Rogora, e Francesco Saverio Tortoriello. 2017. «La matematica come collante culturale nell'insegnamento.» *Matematica, Cultura e Società* 2 (1): 293-304.

Cremona, Luigi. 1865. «Sur l'hypocycloïde à trois rebroussements.» *J. für Math.* 64: 101-123.

De Marchis, Martina, e Enrico Rogora. 2017. «Attualità delle riflessioni di Guido Castelnuovo sulla Formazione dell'insegnante di Matematica.» *Periodico di Matematiche* 71-79.

De Marchis, Martina, Marta Menghini, e Enrico Rogora. 2020. «The importance of "Extensive teaching" in the education of prospective teachers of Mathematics.» *19th Conference on Applied Mathematics APLIMAT 2020.* Bratislava : Slovak University of Technology. 344-353.

Drijvers, P., C. Kieran, e M. A. Mariotti. 2010. «Integrating Technology into Mathematics Education: Theoretical Perspectives.» In *Mathematics Education and Technology-Rethinking the Terrain*, di Hoyles C. e Lagrange J. B., 89-132. New York: Springer.

Euclide. s.d. *Elements.* https://mathcs.clarku.edu/~djoyce/java/elements/elements.html.

Hilbert, David. 2009. *Fondamenti della Geometria.* Milano: Franco Angeli.

Hirst, Thomas Archer. 1865. «Sull'inversione quadrica delle curve piane.» *Annali di Matematica Pura ed Applicata* VII (2): 49-65.

Ignoto, Autore. s.d. «Construct the focii and asymptotes with compass and ruler, given a non rectangular hyperbola.» *Mathematics Stack Exchange.* Consultato il giorno ottobre 25, 2020. https://math.stackexchange.com/questions/3784473/construct-the-focii-and-asymptotes-with-compass-and-ruler-given-a-non-rectangul?newreg=94489116d2344398a3c177ec7243b487.

Jahnke, H. N., A. Arcavi, E. Barbin, O. Bekken, F. Furinghetti, A. Idrissi, C. M. S. da Silva, e C. Weeks. 2000. «The use of original sources in mathematics classroom.» In *History in Mathematics Education*, di J. Fauvel e J. van Maanen, 291-328. Dordrecht: Kluver.

Laborde, C., C. Kynigos, K. Hollebrands, e Sträßer R. 2006. «Teaching and learning geometry with technology.» In *Handbook of research on the psychology of mathematics education: Past, present and future.*, di A. Gutierrez e P. Boero, 275-304. Rotterdam: Sense Publishers.

Lakatos, Imre. 1979. *Dimostrazioni e confutazioni. La logica della scoperta matematica.* Milano: Feltrinelli.

Leung A., Baccaglini-Frank A. 2017. *Digital Technologies in Designing Mathematics Education Tasks.* New York: Springer.

Mariotti, M. A., D. Paola, O. Robutti, e D. Venturi. 2004. *Quaderno interattivo di geometria.* Prod. Media Direct. Bassano del Grappa.
Mariotti, M. A., e A. Baccaglini-Frank. 2018. «Developing the Mathematical Eye Through Problem Solving in a Dynamic Geometry Environment.» In *Broadening the Scope of Research on Mathematical Problem Solving. Research in Mathematics Education*, di N. Amado, S. Carriera e K. Jones, 153-176. New York: Springer.
Mariotti, Maria Alessandra, e Andrea Maffia. 2018. «Dall'utilizzo degli artefatti ai significati matematici: il ruolo dell'insegnante nel processo di mediazione semiotica.» *idattica della matematica. Dalla ricerca alle pratiche d'aula.* (4): 50-64.
Mumma, John, e Marco Panza. 2012. «Diagrams in Mathematics:History and Philosophy.» *Syntheses* 1-5.
Netz, R. 1999. *The Shaping of Deduction in Greek Mathematics. A STudy in Cognitive History.* Cambridge University Press: Cambridge.
Nurzia, Laura. 1999. *La corrispondenza di Luigi Cremona (1830-1903), vol. IV.* Milano: PRISTEM Università Bocconi.
Polya, George. 1971. *La scoperta della matematica. Capire, imparare e insegnare a risolvere i problemi.* Milano: Feltrinelli.
Poncelet, J. V. 1822. *Traité des Propriétés projectives des Figures.* Paris: Bachelier.
Poncelet, J.V. 1822. *Traité des propriétés projectives des figures.* Paris: Gauthier-Villars.
Raspanti, Maria Anna. 2016. «Dall'inversione circolare all'inversione quadratica: aspetti storici e potenzialità didattiche.» *Bollettino Mathesis* xxx-yyy.
Raspanti, Maria Anna. s.d. «Giusto Bellavitis e la sua geometria di derivazione.» *Bollettino di Storia delle Scienze Matematiche* xxx-yyy.
—. 2017. *Tesi di dottorato.* Torino: Università.
Rogora, Enrico. 2020. «Materiali preparati con GeoGebra.» *Matstor.* 28 Settembre. Consultato il giorno Settembre 28, 2020. http://matstor.wikidot.com/materialiggb.
Rogora, Enrico, e Francesco Saverio Tortoriello. 2020. «Interdisciplinarity for learning and teaching mathematics.» *Bolema* xxx-yyy.
Rogora, Enrico, e Saverio Tortoriello. 2018. «Matematica e cultura umanistica.» *Archimede* (2): 82-88.
Russo, Lucio, e Giuseppina Pirro. 2017. *Euclide: il I libro degli Elementi.* Roma: Carocci.
Saito, Ken. 2012. «Traditions of the diagram, tradition of the text: a case study.» *Synthèse* 186 (1): 7-20.
Sapienza, Università di Roma. s.d. *Moodle Sapienza.* Consultato il giorno ottobre 26, 2020. https://elearning.uniroma1.it/.
UMI. 2005. *Riflessioni sul Laboratorio di Matematica.* https://umi.dm.unibo.it/materiali-umi-ciim/trasversali/riflessioni-sul-laboratorio-di-matematica/.
Vaccaro, A. 2020. «Historical origins of the nine-point conic. The contribution of Eugenio Beltrami.» *Historia* 51: 26-48.
Weisstein, Eric W. s.d. «Deltoid.» *Math World - A Wolfram Web Resource.* Consultato il giorno ottobre 26, 2020. https://mathworld.wolfram.com/Deltoid.html.